\documentclass[11pt]{article}
\usepackage[centertags]{amsmath}
\usepackage{amsfonts}
\usepackage{amssymb}
\usepackage{amsthm}
\usepackage{newlfont}
\usepackage[pdftex]{graphicx}
\usepackage{epstopdf}
\usepackage{color}

\newtheorem{Theorem}{Theorem}[section]
\newtheorem{Definition}[Theorem]{Definition}
\newtheorem{Proposition}[Theorem]{Proposition}
\newtheorem{Lemma}[Theorem]{Lemma}

\newtheorem{Remark}[Theorem]{Remark}
\newtheorem{Example}[Theorem]{Example}

\newtheorem{Hypothesis}{Hypothesis}

\numberwithin{equation}{section}

\begin{document}
\renewcommand{\figurename}{Fig.1}

\def\le{\left}
\def\r{\right}
\def\cost{\mbox{const}}
\def\a{\alpha}
\def\d{\delta}
\def\ph{\varphi}
\def\e{\epsilon}
\def\la{\lambda}
\def\si{\sigma}
\def\La{\Lambda}
\def\B{{\cal B}}
\def\A{{\mathcal A}}
\def\L{{\mathcal L}}
\def\O{{\mathcal O}}
\def\bO{\overline{{\mathcal O}}}
\def\F{{\mathcal F}}
\def\K{{\mathcal K}}
\def\H{{\mathcal H}}
\def\D{{\mathcal D}}
\def\C{{\mathcal C}}
\def\M{{\mathcal M}}
\def\N{{\mathcal N}}
\def\G{{\mathcal G}}
\def\T{{\mathcal T}}
\def\R{{\mathcal R}}
\def\I{{\mathcal I}}

\def\bw{\overline{W}}
\def\phin{\|\varphi\|_{0}}
\def\s0t{\sup_{t \in [0,T]}}
\def\lt{\lim_{t\rightarrow 0}}
\def\iot{\int_{0}^{t}}
\def\ioi{\int_0^{+\infty}}
\def\ds{\displaystyle}
\def\pag{\vfill\eject}
\def\fine{\par\vfill\supereject\end}
\def\acapo{\hfill\break}

\def\Amu{{A_\mu}}
\def\Qmu{{Q_\mu}}
\def\Smu{{S_\mu}}
\def\H{{\mathcal{H}}}
\def\Im{{\textnormal{Im }}}
\def\Tr{{\textnormal{Tr}}}
\def\P{{\mathbb{P}}}

\def\beq{\begin{equation}}
\def\eeq{\end{equation}}
\def\barr{\begin{array}}
\def\earr{\end{array}}
\def\vs{\vspace{.1mm}   \\}
\def\rd{\reals\,^{d}}
\def\rn{\reals\,^{n}}
\def\rr{\reals\,^{r}}
\def\bD{\overline{{\mathcal D}}}
\newcommand{\dimo}{\hfill \break {\bf Proof - }}
\newcommand{\nat}{\mathbb N}
\newcommand{\E}{\mathbb E}
\newcommand{\Pro}{\mathbb P}
\newcommand{\com}{{\scriptstyle \circ}}
\newcommand{\reals}{\mathbb R}

\title{On the Smoluchowski-Kramers approximation for SPDEs and its interplay with  large deviations and long time behavior}

\author{Sandra Cerrai\thanks{Partially supported by the NSF grant DMS 1407615.}\\
\normalsize University of Maryland, College Park\\ United States
\and
Mark Freidlin\thanks{Partially supported by the NSF grant DMS 1411866.}\\
\normalsize University of Maryland, College Park\\ United States
\and
Michael Salins\\
\normalsize Boston University, Boston\\ United States}
\date{}

\maketitle

\begin{abstract}
We discuss here the validity of the small mass limit (the  so-called Smoluchowski-Kramers approximation) on a fixed time interval for a class of semi-linear stochastic wave equations, both in the case of the presence of a constant friction term and in the case of the presence of a constant magnetic field. We also consider the small mass limit in an infinite time interval and we see how the approximation is stable in terms of the invariant measure and of the large deviation estimates and the exit problem from a bounded domain of the space of square integrable functions.
\end{abstract}

\section{Introduction}
\label{sec1}

The motion of a particle of a mass $\mu$ in the field
$b(q)+\sigma(q)\dot{W}$, with a constant damping proportional to the speed,
 is described, according to the
Newton law, by the Langevin equation

\begin{equation}
\label{ieq1}
 \mu\ddot{q}^{\mu}_{t}=b(q^{\mu}_{t})+\sigma(q^{\mu}_{t})\dot{W}_{t}-
 \dot{q}_{t}, \ \ \ \ \
q_{0}^{\mu}=q \in \reals^n, \ \ \  \dot{q}^{\mu}_{0}=p \in \reals^n,
\end{equation}
(for the sake of simplicity the friction coefficient is taken equal to 1).

Here $b(q)$ is the deterministic component of the force and
$\sigma(q)\dot{W}_{t}$, where $\dot{W}_{t}$ is the standard
Gaussian white noise in $\reals^{n}$ and $\sigma(q)$ is an $n\times
n$-matrix, is the stochastic part. It is  known that, for
$0<\mu \ll1$, the random position $q_{t}^{\mu}$ of the particle can be approximated by the solution of
the first order equation

\begin{equation}
\label{ieq2}
 \dot{q}_{t}=b(q_{t})+\sigma(q_{t})\dot{W}_{t},\ \ \ \
q_{0}=q \in \reals^n,
\end{equation}
in the sense that
\begin{equation}
\label{ieq3} \lim_{\mu \downarrow 0}\, P\le(\max_{0\leq t \leq T}
|q_{t}^{\mu} -q_{t} |>\delta\r)=0,
\end{equation}
for any $0 \leq T <\infty$ and $\delta>0$ fixed. Statement \eqref{ieq3}
is called then {\em Smoluchowski-Kramers approximation} of $q_{t}^{\mu}$ by
$q_{t}$ (see to this purpose \cite{kra,smolu,freidlin}). This statement justifies the description of the motion of
a small particle by the first order equation \eqref{ieq2} instead
of the second order equation \eqref{ieq1}. Several authors have considered generalizations of this phenomenon in the presence of a magnetic field (see \cite{smolu3,jjl}) and for a non-constant   friction, both in the case it is strictly positive, as in \cite{fh,hvw,hvw12}, and in the case it is possibly vanishing, as in  \cite{fhw}. Large deviations and exit problems in the small noise regime  have also been studied (see \cite{cf}).

\bigskip

In the present paper, we will review some results about the Smoluchowski-Kramers approximation for systems with an infinite number of degrees of freedom,  that the three authors have obtained  in a series of papers written in last several years (cfr. \cite{smolu2,smolu1,csal,cs-annals,cs-asy}).

Let $\mathcal{O}$ be a bounded smooth domain of $\mathbb{R}^d$, with $d\geq 0$. We are dealing here with  the following stochastic semi-linear damped wave equation on $\mathcal{O}$ \begin{equation}
\label{wave-intro-rev}
\le\{\begin{array}{l} \ds{\mu \frac{\partial^2
u}{\partial t^2}(t,x)+\frac{\partial u}{\partial t}(t,x)=\Delta u(t,x)+b(x,u(t,x))+ g(x,u(t,x))
\frac{\partial w^Q}{\partial t}(t,x),\ x \in\,\mathcal{O},\ t\geq 0,}\\
\vs \ds{u(0,x)=u_0(x),\ \ \ \ \ \frac{\partial u}{\partial
t}(0,x)=v_0(x),\ \ \ \ \ \ u(t,x)=0,\ \ \ \ x \in \partial \mathcal{O}.} \end{array}\r.
\end{equation}

Equation \eqref{wave-intro-rev} models the displacement of an elastic material with mass density $\mu>0$ in the region $\mathcal{O}$, which is exposed to deterministic and random forces. The term $-\partial u/\partial t$  models the damping, the Laplacian $\Delta$ models the forces that neighboring particles exert on each other and  the non-linearity $b$ models some deterministic forcing. State-dependent stochastic perturbations are modeled by the term $g\,\partial w^Q/\partial t$, for some Wiener process $w^Q$, that is white in time and $Q$-correlated in space (see Section \ref{se3} for all definitions), and for some nonlinear coefficient $g$.

Many authors have studied stochastic wave equations under various assumptions on the non-linear coefficients $b$ and $g$, on the correlation $Q$ of the noise and on the domain  $\mathcal{O}$ and the boundary conditions (see for example \cite{carmona,dalang,dpb,kaza,millet1,millet2,ondre,ondre2,peszat,peszat2}). We will specify our hypotheses in the following sections.

\medskip
In this paper we explore the  asymptotic behavior  of the solution to \eqref{wave-intro-rev} as the mass density of the material $\mu$ vanishes. This is the infinite dimensional analogue to the Smoluchowski-Kramers approximation described in \eqref{ieq3}. In Sections \ref{sec3} and \ref{sec4}, where the case of additive and multiplicative noise are considered, respectively, we review the results of \cite{smolu2,smolu1} to show that as, $\mu \to 0$, the solution of equation  \eqref{wave-intro-rev} converges to the solution of the stochastic heat equation
\begin{equation}
  \label{heat-intro-rev}
  \left\{
  \begin{array}{l}
   \ds{\frac{\partial u}{\partial t}(t,x) = \Delta u(t,x) + b(x,u(t,x)) + g(x,u(t,x))\frac{\partial w^Q}{\partial t}(t,x),\ x \in \mathcal{O},\ t \geq 0,}\\
   \vs
   \ds{u(0,x) = u_0(x), \ \ \ u(t,x) = 0, \ \ x \in \partial \mathcal{O}. }
   \end{array}
   \right.
\end{equation}
Specifically, if we denote by $u^\mu$ and $u$ the solutions to \eqref{wave-intro-rev} and \eqref{heat-intro-rev}, respectively, we show that for any $T>0$ and $\d>0$ fixed,
\begin{equation} \label{SK-eq-intro}
  \lim_{\mu \to 0} \Pro \left( \sup_{0\leq t\leq T} \int\limits_{\mathcal{O}} |u^\mu(t,x) - u(t,x)|^2 dx > \d \right) = 0.
\end{equation}
This means that when the mass density of a material is small, the stochastic heat equation approximates the stochastic wave equation well, on any finite time interval.

The proof of this infinite dimensional Smoluchowski-Kramers approximation requires uniform estimates on the Sobolev regularity in space and the   H\"older  regularity in time for the solutions to the stochastic wave equation \eqref{wave-intro-rev}. As known, such uniform bounds are used  to establish tightness in an appropriate functional space. Once we obtain a weakly convergent subsequence, by the Prokhorov theorem, the identification of the unique limit, and hence the convergence of the whole sequence, is obtained by using a non-trivial integration-by-parts formula for SPDEs.

\medskip

In Section \ref{sec5} we review the results of \cite{csal}, which concern the Smoluchowski-Kramers approximation for an electrically charged material in the presence of a uniform magnetic field. In this case, the displacement of the material is modeled by the equation
\begin{equation}
  \label{wave-magnetic-intro}
  \left\{
  \begin{array}{l}
   \ds{\mu \frac{\partial^2 u}{\partial t^2}(t,x) = \Delta u(t,x) + b(u(t,x),x,t) + \vec{m} \times \frac{\partial u}{\partial t} + g(u(t,x),x,t)\frac{\partial w^Q}{\partial t}(t,x), }\\
   \vs
   \ds{u(0,x) = u_0(x), \ \ \frac{\partial u}{\partial t}(0,x) + v_0(x) ,\ x \in \mathcal{O}, \ \ u(t,x) = 0, \ x \in \partial \mathcal{O},}
  \end{array}
  \right.
\end{equation}
where $\vec{m} = (0,0,m)$ is a constant vector field that is perpendicular to the plane of motion of the material and that models a magnetic field (notice that here the symbol $\times$ denotes the usual vector product in $\mathbb{R}^3$). The material is assumed to be electrically charged and therefore the above equation describes the movement of an elastic material with constant mass density $\mu>0$, that is exposed to the electric field, some deterministic forcing $b$, and a state dependent stochastic forcing $g\,\partial w^Q/\partial t$.

One might hope, based on the results from Sections \ref{sec3} and \ref{sec4}, that, for any $T>0$ and $\delta>0$ fixed, a limit analogous to \eqref{SK-eq-intro} holds,
where $u_\mu$ is the solution to equation \eqref{wave-magnetic-intro} and $u$ is the solution to  the equation
\begin{equation}
  \label{magnetic-limit-intro}
  \left\{
  \begin{array}{l}
    \ds{\frac{\partial u}{\partial t}(t,x) = J_0^{-1} \left[\Delta u(t,x) + b(u(t,x),x,t) + g(u(t,x),x,t)\, \frac{\partial w^Q}{\partial t}(t,x) \right],}\\
    \vs
    \ds{u(t,x) = 0, \ \ \ x \in \partial \mathcal{O}, \ \ \ \ \ \ \ \ \ u(0,x) =u_0(x),\ \ \ x \in\,\mathcal{O},}
  \end{array}
  \right.
\end{equation}
where $J_0^{-1}$ is the inverse of the matrix
\[J_0= \begin{pmatrix}
  0 & 1 \\ -1 & 0
\end{pmatrix}.\]
Unfortunately, because of the presence of the stochastic  term, the small-mass limit \eqref{SK-eq-intro} is not true anymore. A similar situation was explored in \cite{smolu3}, where it was shown that, in the case of a system with a finite number of degrees of freedom,   the difference $u_\mu-u$ does not converge to zero, as $\mu$ tends to zero.   Thus, in Section \ref{sec5} we study the problem of the small mass limit by  adding a small friction $-\e\,\partial u/\partial t$ to equation \eqref{wave-magnetic-intro} and we show that,    for any fixed $\e>0$, the  solution to \eqref{wave-magnetic-intro}, with fixed positive friction, converges to the solution of the system of semi-linear stochastic heat-equations formally obtained by setting $\mu = 0.$

Here, we use a slightly different line of argument than in Sections \ref{sec3} and \ref{sec4}, which allows us to prove the validity of the convergence  in $L^p(\Omega;C([0,T];L^2(\mathcal{O})))$. This means that, if  $u_\mu^\e$ is the solution to the wave equation exposed to a magnetic field and an $\e$ friction, and if $u_\e$ is the solution to the equation
\[ \left\{
  \begin{array}{l}
    \ds{\frac{\partial u}{\partial t}(t,x) = \le(\e\,I+J_0\r)^{-1} \left[\Delta u(t,x) + b(u(t,x),x,t) + g(u(t,x),x,t)\, \frac{\partial w^Q}{\partial t}(t,x) \right],}\\
    \vs
    \ds{u(t,x) = 0, \ \ \ x \in \partial \mathcal{O}, \ \ \ \ \ \ \ \ \ u(0,x) =u_0(x),\ \ \ x \in\,\mathcal{O},}
  \end{array}
  \right.\]
then for any $p\geq 1$ and $T\geq0$ fixed,
\begin{equation}
  \lim_{\mu \to 0} \E\, \le(\sup_{t \leq T} \int\limits_{\mathcal{O}}\left|u_\mu^\e(t,x) - u_\e(t,x)\right|^2dx \r)^{\frac p2} = 0.
\end{equation}
Notice that the $L^p$ convergence can also be proven in the case  there is no magnetic field, as long as $b$ and $g$ are suitably regular.

\bigskip

In Sections \ref{sec6}, \ref{sec7}, and \ref{sec8}, we investigate the multi-scale interactions between the small mass limit  and  long-time behaviors of the stochastic wave equation. In Sections \ref{sec3} to \ref{sec5}, we demonstrate that the solutions to \eqref{wave-intro-rev} converge to the solution to \eqref{heat-intro-rev} in the topology of $C([0,T];L^2(\mathcal{O}))$. In fact \eqref{SK-eq-intro}  is only true for fixed finite $T>0$ and, on longer time scales, the solutions will deviate arbitrarily far apart in a pathwise sense. But, despite the fact that the Smoluchowski-Kramers approximation is not valid in a pathwise sense on long-time scales, there are still many important long-time characteristics of the small-mass wave equation that are approximated by the heat equation.

\medskip

In Section \ref{sec6} we review some results from \cite{smolu2} concerning the relationship between the Smoluchowski-Kramers approximation and invariant measures for equations \eqref{wave-intro-rev} and \eqref{heat-intro-rev}, in the case of gradient systems. We recall here that by gradient system we mean  a system where the noise is additive (that is $g\equiv I$) and, if $B: L^2(\mathcal{O}) \to L^2(\mathcal{O})$ is defined as
\[B(h)(x) = b(h(x)),\ \ \ \ x \in\,\mathcal{O},\]
then
there exists a potential functional $F:L^2(\mathcal{O}) \to \mathbb{R}$,  such that for any $h \in L^2(\mathcal{O})$
\[B(h) = -Q^2 DF(h),\]
where $DF$ is the Fr\'echet derivative of $F$  in $L^2(\mathcal{O})$ and $Q$ is the correlation of the noise.

Due to the fact that  \eqref{wave-intro-rev} is a gradient system, by using suitable finite dimensional approximations, we show that the invariant measure for the pair $(u^\mu,\partial u^\mu/\partial t)$ solving \eqref{wave-intro-rev} in the space $L^2(\mathcal{O})\times H^{-1}(\mathcal{O})$ is given by the probability measure
\[\nu_\mu(du,dv):=\frac 1{Z}\,e^{\,-2 F(u)}\,\mathcal{N}\le(0,
(-\Delta)^{-1}Q^2/2\r)(du)\times \mathcal{N}\le(0,
(-\Delta)^{-1}Q^2/2\mu\r)(dv),\ \ \ \ \mu>0,\]
where $Z$ is a  normalizing constant, independent of $\mu>0$. This means that we have  an explicit expression of the density of the invariant measure with respect to a suitable Gaussian measure on the product space $L^2(\mathcal{O})\times H^{-1}(\mathcal{O})$. In particular, due to the special form of $\nu_\mu$,  its first marginal $\Pi_1 \nu_\mu$  does not depend on $\mu>0$ and coincides with the probability measure
\[\nu(du):=\frac 1{Z}\,e^{\,-2 F(u)}\,\mathcal{N}\le(0,
(-\Delta)^{-1}Q^2/2\r)(du),\]
 which is the invariant measure of system \eqref{heat-intro-rev}. In particular, we have that the Smoluchowski-Kramers approximation is valid in the sense that the stochastic heat equation and the damped stochastic wave equation have the same long-time behavior.

The case of non-gradient systems is still open, and of course we cannot expect to have any explicit expression for the invariant measures of systems \eqref{wave-intro-rev} and \eqref{heat-intro-rev}. Nevertheless,  we expect that if $\nu_\mu$ and $\nu$ denote the invariant measures of \eqref{wave-intro-rev} and \eqref{heat-intro-rev}, respectively, then some sort of convergence for $\Pi_1 \nu_\mu$ to $\nu$ holds, in the small mass $\mu$ limit.

\medskip

In Sections \ref{sec7} and \ref{sec8} we review results from \cite{cs-asy} and \cite{cs-annals} about the relationship between the small mass and the small noise asymptotics. More precisely, we consider for any $\e>0$ the stochastic damped wave equation
\begin{equation} \label{wave-intro-small-noise}
  \mu\frac{\partial^2 u}{\partial t^2}(t,x) + \frac{\partial u}{\partial t}(t,x) = \Delta u(t,x) + b(x,u(t,x)) + \sqrt{\e}\,\frac{\partial w^Q}{\partial t}(t,x),
\end{equation}
with appropriate initial and boundary conditions and the corresponding heat equation
\begin{equation} \label{heat-intro-small-noise}
  \frac{\partial u}{\partial t}(t,x) = \Delta u(t,x) + b(x,u(t,x)) + \sqrt{\e}\,\frac{\partial w^Q}{\partial t}(t,x),
\end{equation}
with the same initial and boundary conditions. We  are interested in the multiscale behavior of system \eqref{wave-intro-small-noise},  as both $\mu$ and $\e$ go to zero.

From the earlier results, it is clear that if we first take the limit as  $\mu \to 0$ and then as $\e \to 0$,  the large deviation principle for the heat equation should describe the behavior of the system. Our purpose here is to  show that when the order in the two limits is reversed and we first take the limit as $\e \to 0$ and then $\mu \to 0$, the large deviations principle for the heat equation is still appropriate for studying the long-time behaviors of the wave equation. In particular, this result provides a rigorous mathematical justification of what is done in  applications, when, in order to study rare events and transitions between metastable states for the more complicated system \eqref{wave-intro-small-noise}, as well as exit times from basins of attraction and the corresponding exit places,  the relevant quantities associated with the large deviations for system \eqref{heat-intro-small-noise} are considered.

Because of the dissipation introduced by the friction term, under suitable conditions on $b$, the  solution to the unperturbed version of  \eqref{wave-intro-small-noise} (that is with $\e=0$) will converge to 0 as $t \to +\infty$. When the system is exposed to small random perturbations (that is when $0< \e\ll 1$), the solution will deviate from this equilibrium point on long time scales. It is thus of interest to study exit times of the form
\[\tau^{\mu,\e} = \inf\{t>0: u^\mu_\e(t,\cdot) \not \in D\},\]
where $D \subset L^2(\mathcal{O})$ contains the equilibrium solution $0$.

By extending to this infinite dimensional setting well known results in the theory of large deviations for finite dimensional systems (see \cite{dz,fw}), in Section \ref{sec8} we show that the logarithmic exit time asymptotics, as well as the logarithmic expected value of the exit time and the exit position $u^\mu_\e(\tau^{\mu,\e},\cdot)$, can be characterized by  the quasi-potential $\bar{V}_\mu$.
More precisely, we show that, for small enough fixed $\mu>0$ and appropriate initial conditions
\begin{equation}
\label{intro22}
\lim_{\e\to 0}\,\e\log\,\E\tau^{\mu,\e}=\inf_{u \in\,\partial D} \bar{V}_\mu(u),
\end{equation}
and
\begin{equation} \label{intro23}
  \lim_{\e \to 0}\,\e\log\,(\tau^{\mu,\e}) = \inf_{u \in\,\partial D} \bar{V}_\mu(u),\ \ \  \text{ in probability}.
\end{equation}
We also  prove that if $N \subset \partial D$ has the property that $\ds{\inf_{u \in N} \bar{V}_\mu(u) > \inf_{u \in \partial D} \bar{V}_\mu(u) }$, then
\begin{equation} \label{intro24}
  \lim_{\e\to0}\, \Pro \left(u^\mu_{\e}(\tau^{\mu,\epsilon}) \in N  \right) =0.
\end{equation}

Thus, due to the role played by the quasi-potential in the description of these important asymptotic features of the system, our purpose here is  to compare the quasi-potential $V^\mu(u,v)$ associated with \eqref{wave-intro-small-noise}, with the quasi-potential $V(u)$ associated with \eqref{heat-intro-small-noise}, and  to show that for any closed set $N \subset L^2(\mathcal{O})$ it holds
\begin{equation}
\label{intro10}
\lim_{\mu\to 0}\,\inf_{u \in N}\,\bar{V}_\mu(u):=\lim_{\mu\to 0}\,\inf_{u \in N}\,\inf_{v \in\,H^{-1}(\mathcal{O})}V^\mu(u,v)=\inf_{u \in N}\,V(u).
\end{equation}

This means that, in the description of the large deviation principle,  taking first the limit as  $\e\downarrow 0$ (large deviation) and then taking the limit as $\mu\downarrow 0$ (Smoluchowski-Kramers approximation) is the same as first taking the limit as  $\mu\downarrow 0$ and then   as  $\e\downarrow 0$.

\medskip

In Section \ref{sec7}, we  address this problem in  the particular case system \eqref{wave-intro-small-noise} is of gradient type. As for the invariant measures studied in Section \ref{sec6},  in the case of gradient systems all relevant quantities associated with the large deviation  can be explicitly computed. In particular, we show that
 for any $\mu>0$
 \begin{equation}
 \label{intro11}
V^\mu(u,v) = \left|(-\Delta)^{1/2}Q^{-1} u \right|_{L^2(\mathcal{O})}^2 + 2F(u) + \mu \left|Q^{-1} v \right|_{L^2(\mathcal{O})}^2,
  \end{equation}
  for any $(u,v) \in\,\textnormal{Dom}((-\Delta)^{1/2}Q^{-1})\times \textnormal{Dom}(Q^{-1}).$
Therefore, as
\[V(u) = \left|(-\Delta)^{1/2} Q^{-1} u \right|_{L^2(\mathcal{O})}^2 + 2F(u),\ \ \ \ u \in\,\textnormal{Dom}((-\Delta)^{1/2}Q^{-1}), \]
from \eqref{intro11} we have   that for any $\mu>0$,
\begin{equation}
\label{intro20}   \bar{V}_\mu(u):= \inf_{v \in H^{-1}(\mathcal{O})} V^\mu(u,v) = V^\mu(u,0) = V(u),\ \ \ \ u \in\,\textnormal{Dom}((-\Delta)^{1/2}Q^{-1}).
\end{equation}
In particular, this means that $\bar{V}_\mu(u)$ does not just coincide with $V(u)$ at the limit,  but for any fixed $\mu>0$.

\medskip

In the general non-gradient case considered in Section \ref{sec8}, the situation is considerably more delicate and we cannot expect anything explicit as in \eqref{intro11}.  The lack of an explicit expression for $V^\mu(u,v)$ and $V(u)$ makes the proof of  \eqref{intro10} much more difficult and requires the introduction of new arguments and techniques.

The first key idea in order to prove \eqref{intro10} is to characterize $V^\mu(u,v)$ as the minimum value for a suitable functional. We recall that the quasi-potential $V^\mu(u,v)$ is defined as the minimum energy required to the system to go from the asymptotically stable equilibrium $0$ to the point $(u,v) \in\,L^2(\mathcal{O})\times H^{-1}(\mathcal{O})$, in any time interval. Namely
\[V^\mu(u,v) = \inf \left\{ I^\mu_{0,T}(z) \ ;\  z(0) = 0,\  z(T) =(u,v),\  T > 0 \right\},\]
where
\[I^{\mu}_{0,T}(z) =\frac{1}{2}  \inf \left\{ | \psi|_{L^2((0,T);L^2(\mathcal{O})}^2\,:\, z=z^\mu_{\psi}\right\},
\] is the large deviation action functional
 and $z^\mu_{\psi} = (u^\mu_\psi, \partial u^\mu_\psi/\partial t)$ is a mild solution of the skeleton equation associated with equation  \eqref{wave-intro-small-noise}, with control $\psi \in\,L^2((0,T);L^2(\mathcal{O}))$,
\begin{equation}
\label{intro12}
\mu \frac{\partial^2 u^\mu_\psi}{\partial t^2}(t) = \Delta u^\mu_\psi(t) - \frac{\partial u^\mu_\psi}{\partial t}(t) + B(u^\mu_\psi(t)) + Q \psi(t),\ \ \ \ t \in\,[0,T].
\end{equation}
By working thoroughly with the skeleton equation \eqref{intro12}, we  show  that,
for small enough $\mu>0$,
\begin{equation}
\label{intro13}
  V^\mu(u,v) = \min \left\{ I^\mu_{-\infty,0}(z): \lim_{t \to -\infty} |z(t)|_{L^2(\mathcal{O})\times H^{-1}(\mathcal{O})} =0,\  z(0)=(u,v) \right\}\end{equation}
  where the minimum is taken over all $z \in C((-\infty,0];L^2(\mathcal{O})\times H^{-1}(\mathcal{O}))$.
In particular, we get that the level sets of $V^\mu$ and $\bar{V}_\mu$ are compact in $L^2(\mathcal{O})\times H^{-1}(\mathcal{O})$ and $L^2(\mathcal{O})$, respectively. Moreover, we show that both $V^\mu$ and $\bar{V}_\mu$ are well defined and continuous in suitable Sobolev spaces of functions.

The second key idea is based on the fact  that,  as in \cite{cf} where the finite dimensional case is studied, for all functions $z \in\,C((-\infty,0];L^2(\mathcal{O})\times H^{-1}(\mathcal{O}))$ that are regular enough, if we denote $\varphi(t)=\Pi_1 z(t)$, we have
\begin{equation} \label{intro14}
\begin{array}{l}
\ds{  I^\mu_{-\infty} (z) = I_{-\infty}(\varphi) + \frac{\mu^2}{2} \int_{-\infty}^0 \left| Q^{-1} \frac{\partial^2\varphi}{\partial t^2}(t) \right|_{L^2(\mathcal{O})}^2 dt}\\
\vs
\ds{  + \mu \int_{-\infty}^0 \left< Q^{-1} \frac{\partial^2\varphi}{\partial t^2}(t), Q^{-1} \left( \frac{\partial\varphi}{\partial t}(t) - A \varphi(t) - B(\varphi(t)) \right) \right>_{L^2(\mathcal{O})} dt}\\
\vs
\ds{=:I_{-\infty}(\varphi)+J^\mu_{-\infty}(z).}
\end{array}
\end{equation}
Thus, if $\bar{z}^\mu$ is the minimizer of $\bar{V}_\mu(u)$, whose existence is guaranteed by \eqref{intro13}, and if $\bar{z}^\mu$ has enough regularity to guarantee that all terms in \eqref{intro14} are meaningful, we obtain
\begin{equation}
\label{intro16}
\bar{V}_\mu(u)=I_{-\infty}(\bar{\varphi}_\mu)+J^\mu_{-\infty}(\bar{z}^\mu)\geq V(u)+J^\mu_{-\infty}(\bar{z}^\mu).\end{equation}
In the same way, if $\bar{\varphi}$ is a minimizer for $V(u)$ and is regular enough, then
\begin{equation}
\label{intro17}
\bar{V}_\mu(u)\leq I^\mu_{-\infty}(\bar{\varphi},\partial \bar{\varphi}/\partial t)=V(u)+J^\mu_{-\infty}((\bar{\varphi},\partial \bar{\varphi}/\partial t)).\end{equation}
If we could prove that
\begin{equation}
\label{intro18}
\liminf_{\mu\to 0}J^\mu_{-\infty}(\bar{z}^\mu)=\limsup_{\mu\to 0}J^\mu_{-\infty}((\bar{\varphi},\partial \bar{\varphi}/\partial t))=0,\end{equation}
from \eqref{intro16} and \eqref{intro17} we could conclude that \eqref{intro10} holds true. But unfortunately, neither $\bar{z}^\mu$ nor $\bar{\varphi}$ have the required regularity to justify \eqref{intro18}. Thus, we have to proceed with suitable approximations, which, among other things,  require us to prove the continuity of the mappings $\bar{V}_\mu:\textnormal{Dom}((-\Delta)^{1/2}Q^{-1})\to \reals$, uniformly with respect to  $\mu \in\,(0,1]$.

\medskip

In the second part of Section \ref{sec8}, we apply \eqref{intro10} to the study of the exit time and of the exit place of $u^\mu_\e$ from a given domain in $L^2(D)$ . If
\[\tau^\epsilon = \inf\{ t>0: u_\e(t) \not \in D\}\] is the exit time from $D$ for the solution of \eqref{heat-intro-small-noise}, and $V(u)$ is the quasi-potential associated with this system, the exit time and exit place results for the first-order system are analogous to \eqref{intro22}, \eqref{intro23}, and \eqref{intro24}.

As a consequence of \eqref{intro20}, in the gradient case, \eqref{intro22}, and \eqref{intro23} imply that,
for any fixed $\mu>0$, the exit time and exit place asymptotics of \eqref{heat-intro-small-noise} match those of \eqref{wave-intro-small-noise}. In particular, for any $\mu>0$
\begin{equation}
\label{i.1}\lim_{\e\to 0}\,\e\log\,\E\tau^{\mu,\e}=\inf_{u \in\,\partial D} V(u) = \lim_{\e \to 0} \e \log\,\E \tau^\epsilon,\end{equation}
and
\begin{equation}
\label{i.2}
\lim_{\e \to 0} \e \log \tau^{\mu,\e} = \inf_{u \in\,\partial D} V(u) = \lim_{\e \to 0} \e \log\, \tau^\e,\ \ \ \  \text{ in probability}.\end{equation}

In the general non-gradient case, we cannot have \eqref{i.1} and \eqref{i.2}. Nevertheless,  in view of \eqref{intro10}, the exit time and exit place asymptotics of \eqref{wave-intro-small-noise} can be approximated by $V$. Namely
 \[\lim_{\mu \to 0}\lim_{\e\to 0}\,\e\log\,\E\tau^{\mu,\e}=\inf_{u \in\,\partial D} V(u) = \lim_{\e \to 0} \e \log\,\E \tau^\epsilon,\]
and
 \[\lim_{\mu \to 0}\lim_{\e \to 0} \e \log \tau^{\mu,\e} = \inf_{u \in\,\partial D} V(u) = \lim_{\e \to 0} \e \log\, \tau^\e,\ \ \ \  \text{ in probability}.\]

\section{Notations, assumptions and a few preliminary results}
\label{sec2}

We denote by $\mathcal{O}$  a bounded open subset of $\rd$, with $d\geq
1$, and we assume that $\mathcal{O}$ is of class $C^3$. In what
follows we shall denote by $\{e_k\}_{k \in\,\nat}$ the complete
orthonormal basis which diagonalizes the Laplace operator $\Delta$,
endowed with Dirichlet boundary conditions on $\partial\,\mathcal{O}$.
Moreover we shall denote by $\{-\a_k\}_{k \in\,\nat}$ the
corresponding sequence of eigenvalues.

As we are assuming the boundary of $\mathcal{O}$ to be smooth, for
any $\delta \in\,(0,1)$ we have
\[|e_k|_{C^{2+\delta}(\bar{\mathcal{O}})}\leq c\,|\Delta
e_k|_{C^{\delta}(\bar{\mathcal{O}})}=c\,\a_k
\,|e_k|_{C^{\delta}(\bar{\mathcal{O}})},\] (e.g. see
\cite[Theorem 6.3.2]{krylov}) and then, by interpolation, we get
\[|e_k|_{C^{\delta}(\bar{\mathcal{O}})}\leq c\,|e_k|_\infty^{2/2+\d}|e_k|_{C^{2+\delta}(\bar{\mathcal{O}})}^{\delta/2+\delta}\leq c\,\a_k^{\delta/2+\delta}\,|e_k|_\infty^{2/2+\d}
|e_k|_{C^{\delta}(\bar{\mathcal{O}})}^{\delta/2+\delta},\]
(for a proof see \cite[Lemma 6.3.1]{krylov}).
This implies
\[
|e_k|_{C^{\delta}(\bar{\mathcal{O}})}\leq
c\,\a_k^{\delta/2}\,|e_k|_\infty,\ \ \ \ \d \in\,(0,1).
\]

For any $\delta \in\,\reals$, we denote by $H^\delta(\mathcal{O})$
the completion of $C^\infty_0(\mathcal{O})$ with respect to the
norm
\[|h|_{H^\delta(\mathcal{O})}^2=\sum_{i=1}^\infty \a_i^\delta
\le<h,e_i\r>_{L^2(\mathcal{O})}^2.\]
Here and in what follows, for each $h
\in\,H^\delta(\mathcal{O})$ we shall denote by $h_k$ the $k$-th Fourier
coefficient of $h$, that is
\[h_k=\le<h,e_k\r>_{L^2(\mathcal{O})}.\]
$H^\delta(\mathcal{O})$ is a Hilbert space, endowed with the
scalar product
\[\le<h,k\r>_{H^\delta(\mathcal{O})}=\sum_{i=1}^\infty \a_i^\delta
h_i k_i,\ \ \ \ h,k \in\,H^\delta(\mathcal{O}).\] Next, for
any $\delta \in\,\reals$ we denote by $\mathcal{H}_\delta$ the
Hilbert space $H^{\delta}(\mathcal{O})\times
H^{\delta-1}(\mathcal{O})$, endowed with the natural scalar
product and norm inherited from each component. In what follows, we shall denote $L^2(\mathcal{O})=H^0(\mathcal{O})=:H$ and $\mathcal{H}_0=:\mathcal{H}$.

\medskip

For any $\mu>0$ and $\delta \in\,\reals$, we define  the unbounded operator $A_{\mu}$ by setting
\[A_{\mu}(h,k)=\frac 1\mu\le(\mu k, \Delta h-k\r),\ \ \ \ \
(h,k) \in\,D(A_{\mu}):=\mathcal{H}_{\delta+1}.\] The operators $A_\mu$ defined on different $\mathcal{H}_\delta$
are all consistent. It is known that $A_{\mu}$ is the generator of
a  group of bounded linear transformations $\{S_{\mu}(t)\}_{t
\in\,\reals}$ on $\mathcal{H}_\delta$ which is  strongly
continuous (for a proof see e.g. \cite[section 7.4]{pazy}). This means that   for any  $(u_0,v_0) \in\,\mathcal{H}_\delta$ and for any
$\mu>0$, $S_\mu(t)(u_0,v_0)$ is the solution  of the deterministic
linear system
\[\le\{\begin{array}{l} \ds{\frac{\partial u}{\partial
t}(t,x)=v(t,x),\ \ \ \ \mu \frac{\partial v}{\partial
t}(t,x)=\Delta
u(t,x)-v(t,x),\ \ \ \ \ t>0,\ \ x \in\,\mathcal{O},}\\
\vs \ds{u(0)=u_0,\ \ \ v(0)=v_0,\ \ \ \ \ u(t,x)=0,\ \ \ \ \ t\geq
0,\ \ x \in\,\partial \mathcal{O},}
\end{array}\r.\]
which can be written as the following abstract evolution problem
in $\mathcal{H}_\delta$
\[\frac {dz}{dt}(t)=A_\mu z(t),\ \ \ \ \ z(0)=(u_0,v_0),\]
where $z(t):=(u(t),v(t))$.

Next, we notice that the adjoint operator to $A_\mu$ is given by
\[A^\star_\mu(h,k)=\frac 1\mu\le(-k, -\mu \Delta h-k\r),\ \ \ \ \
(h,k) \in\,D(A^\star_{\mu}):=\mathcal{H}_{\delta+1}.\] In what
follows we shall denote by $\{S^\star_\mu(t)\}_{\{t\geq 0\}}$ the
semigroup generated by $A^\star_\mu$.

If we set
$\Pi_1(u,v):=u$ and $\Pi_2(u,v):=v$  we have that
\[\Pi_1\,S_\mu(t)(u,v)=\sum_{k=1}^\infty f^\mu_k(t;u,v) e_k,\ \ \ \ \
\Pi_2\,S_\mu(t)(u,v)=\sum_{k=1}^\infty g^\mu_k(t;u,v) e_k,\] where the
pair $(f^\mu_k(t;u,v),g^\mu_k(t;u,v))$ is for each $k \in\,\nat$ and
$\mu>0$ the solution of the system
\begin{equation}
\label{sistema} \le\{\begin{array}{l} \ds{f^\prime(t)=g(t),\ \ \ \ \ f(0)=u_{k}}\\
\vs \ds{\mu\,g^\prime(t)=-\a_k \, f(t)-g(t),\ \ \ \ \ g(0)=v_{k}.}
\end{array}\r.\end{equation}

In fact, both
$f^\mu_k$ and $g^\mu_k$ can be explicitly computed.

\begin{Proposition}\cite[Prop 2.2]{smolu2}
\label{propy} For any $\mu>0$ and $k \in\,\nat$,  let us define
\[ \gamma^\mu_k:= \frac 1{2\mu}\sqrt{1-4\a_k \mu}.\]
Then, if $4\a_k \mu\neq 1$ we have
\begin{equation}
\label{pop1}
\begin{array}{l}
\ds{f^\mu_k(t;u,v)=\frac 12\exp\le(-\frac
t{2\mu}\r)\le(\le[\le(1+\frac 1{2\mu
\gamma^\mu_k}\r)\exp\le(\gamma^\mu_k t\r)+
\le(1-\frac 1{2\mu\gamma^\mu_k}\r)\exp\le(-\gamma^\mu_k t\r)\r]\,u_k\r.}\\
\vs \ds{\le.+\frac 1{ \gamma^\mu_k}\le[\exp\le(\gamma^\mu_k
t\r)-\exp\le(-\gamma^\mu_k t\r)\r]\,v_k\r),}
\end{array}
\end{equation}
and
\[\begin{array}{l} \ds{ g^\mu_k(t;u,v)=\frac
12\exp\le(-\frac t{2\mu}\r)\le( -\frac {\a_k}{\mu
\gamma^\mu_k}\le[\exp\le(\gamma^\mu_k
t\r)-\exp\le(-\gamma^\mu_k t\r)\r]\,u_k\r.}\\
\vs \ds{\le.+\le[\le(1-\frac
1{2\mu\gamma^\mu_k}\r)\exp\le(\gamma^\mu_k t\r) +\le(1+\frac
1{2\mu\gamma^\mu_k}\r)\exp\le(-\gamma^\mu_k t\r)\r]\,v_k\r).}
\end{array}
\]
Moreover, if $4\a_k \mu=1$ we have
\[ \begin{array}{l} \ds{f^\mu_k(t;u,v)=\exp\le(-\frac
t{2\mu}\r)\le[\le(\frac t{2\mu}+1\r)\,u_{k}+ t\,v_{k}\r],}
\end{array}
\]
and
\[ g^\mu_k(t;u,v)=\exp\le(-\frac t{2\mu}\r)\le[-\frac
{t}{4\mu^2}\,u_{k}+\le(1-\frac t{2\mu}\r)v_{k}\r].
\]

\end{Proposition}

Moreover, in view of the explicit formula \eqref{pop1}, it is possible to prove the following bounds for $f^\mu_k(t;u,v)$.
\begin{Lemma}\cite[(3.12) and (3.13) in Lemma 3.2]{smolu2}
\label{rev-l1}
For every $k \in\,\mathbb{N}$, $\mu>0$  and $t\geq 0$, we have
\[
\int_0^t|f^\mu_k(s;0,e_k)|^2\,ds\leq c\,\frac{\mu^2}{\a_k}.
\]
Moreover, for any $k \in\,\mathbb{N}$, $\mu>0$, $\theta \in\,[0,1]$  and $t>s$, we have
\[
\int_0^s|f^\mu_k(t-r;0,e_k)-f^\mu_k(s-r;0,e_k)|^2\,dr\leq c\,\frac{\mu^2}{\a_k^{1-\theta}}(t-s)^\theta.
\]
\end{Lemma}

An important consequence of Proposition \ref{propy} is the
following result on the asymptotic behavior of $S_\mu(t)$.

\begin{Proposition}\cite[Proposition 2.4]{smolu2}
\label{asy} For any fixed $\mu>0$ and for any $\delta \in\,\reals$
the semigroup $\{S_\mu(t)\}_{t\geq 0}$ is of negative type in
$\mathcal{H}_\delta$. This means that there exist some $\omega_\mu>0$ and
$M_\mu>0$  such that
\[ \|S_\mu(t)\|_{\L(\mathcal{H}_\delta)}\leq
M_\mu\,e^{-\omega_\mu t},\ \ \ \ t \geq 0.
\]
\end{Proposition}

In fact, we have the following representation for the dual semigroup $S^\star_\mu(t)$ in terms of the semigroup $S_\mu(t)$.
\begin{Proposition}\cite[Proposition 2.3]{smolu2}
\label{propybis} For any $\mu>0$ and  $(u,v)
\in\,\mathcal{H}_\delta$ we have

\[ S^\star_\mu(t)(u,v)=\le(\Pi_1\,S_\mu(t)\le(u,-
v/\mu\r),\Pi_2\,S_\mu(t)\le(-\mu u,v\r)\r),\ \ \ \ \ \ t\geq 0.
\]

\end{Proposition}

In the case $\mathcal{O}=[0,L]$, we can prove the following integral representation of $S_\mu(t)$ in terms of a suitable kernel $K_\mu(t,x,y)$.

\begin{Lemma}\cite[Lemma 2.2]{smolu1}
Assume $\mathcal{O}=[0,L]$ and fix $\mu>0$ and $\d<1/2$. For any $v \in\,H^{-\d}(0,L)$ it holds
\[ \Pi_1 S_\mu(t)(0,v)(x)=\int_0^LK_\mu(t,x,y)
v(y)\,dy,\ \ \ \ \ \ (t,x) \in\,[0,\infty)\times [0,L],
\]
 where $K_\mu:[0,\infty)\times [0,L]^2\to \mathbb{R}$
is  defined by
\[
K_\mu(t,x,y):=\sum_{k=1}^{\infty}f^\mu_k(t;0,e_k)e_k(x)
e_k(y).
\]

\end{Lemma}

The kernel
$K_\mu(t,x,y)$ satisfies some regularity properties both in the time and in the space variables, as shown in the next Lemma..

\begin{Lemma}\cite[Lemma 2.3]{smolu1}
\label{para4} Let $\d<1/2$. Then, for any $\rho<1/2-\d$ there
exists some  constant $c_\rho>0$ such that for any $0\leq r\leq t$
and $x,y \in\,[0,L]$
    \[
    \int_0^t\le|K_\mu(s,x,\cdot)-K_\mu(s,y,\cdot)\r|^2_{H^\d}\,ds\leq
    c_\rho \,\mu^2\,|x-y|^{2\rho}
    \]
    and
    \[
    \int_0^r\le|K_\mu(t-s,x,\cdot)-K_\mu(r-s,x,\cdot)\r|^2_{H^\d}\,ds+
    \int_r^t\le|K_\mu(t-s,x,\cdot)\r|^2_{H^\d}\,ds\leq c_\rho
    \,\mu^2\,|t-r|^{\rho}.
    \]
\end{Lemma}

\section{The approximation in the case of additive noise}
\label{sec3}

We are here concerned with the  following stochastic damped semilinear equation
\begin{equation}
\label{quasilinear} \le\{\begin{array}{l} \ds{\mu \frac{\partial^2
u}{\partial t^2}(t,x)=\Delta u(t,x)-\frac{\partial
u}{\partial t}(t,x)+b(x,u(t,x))+ \frac{\partial w^Q}{\partial t}(t,x),\ \ \ \ \ t>0,\ \ x \in\,\mathcal{O},}\\
\vs \ds{u(0,x)=u_0(x),\ \ \ \ \frac{\partial u}{\partial t}(0,x)=v_0(x),\ \
\ \ \ \ u(t,x)=0, \ \ \ \ t\geq 0,\ \ x \in\,\partial
\mathcal{O}.}
\end{array}\r.
\end{equation}
Our aim is proving that the solution $u^\mu(t)$ converges to the
solution of the  stochastic semi-linear heat equation
\begin{equation}
\label{heat} \le\{\begin{array}{l} \ds{\frac{\partial u}{\partial
t}(t,x)=\Delta u(t,x)+b(x,u(t,x))+ \frac{\partial w^Q}{\partial
t}(t,x),\ \ \ \ \
t>0,\ \ x \in\,\mathcal{O},}\\
\vs \ds{u(0,x)=u_0(x),\ \ \ \ \ \ \ \ u(t,x)=0, \ \ \ \ \ t\geq 0,\ \ x
\in\,\partial \mathcal{O},} \end{array}\r.
\end{equation}
as the parameter $\mu$ converges to zero.

\medskip

Here and in what follows, $w^Q(t,x)$ is a cylindrical Wiener process. We shall assume that for any $h,k \in\,H$ and $t,s\geq 0$
\begin{equation}
\label{covariance}
\E\,\le<w^Q(t),h\r>_{H}\le<w^Q(s),k\r>_{H}=(t\wedge
s) \le<Qh,k\r>_{H},\end{equation}
for some operator $Q \in\,\mathcal{L}(H)$. We will assume that $Q$ satisfies the following condition.

\begin{Hypothesis}
\label{H1} The bounded linear operator $Q:H\to
H$ is diagonal with respect to the basis
$\{e_k\}_{k \in\,\nat}$. If   $\{\la_k\}_{k \in\,\nat}$ denotes
the corresponding sequence of eigenvalues,  there exists a
constant $\theta \in\,(0,1)$ such that
\begin{equation}
\label{gamma} \sum_{k=1}^\infty
\frac{\la^2_k}{\a_k^{1-\theta}}\,|e_k|_\infty<\infty.
\end{equation}
\end{Hypothesis}

\begin{Remark}
\em{\begin{enumerate}

\item In several cases, as for example in the case of space dimension $d=1$ or in  the case $\mathcal{O}$ is  a hypercube and the Laplace operator $\Delta$ is endowed with Dirichlet boundary conditions,
the eigenfunctions $e_k$ are equi-bounded in the sup-norm and then
condition \eqref{gamma} becomes \[\sum_{k=1}^\infty
\frac{\la_k^2}{\a_k^{1-\theta}}<\infty.\] In general it holds
\[|e_k|_\infty\leq c\,k^\a,\ \ \ \ \ \ k \in\,\nat,\]
for some $\a\geq 0$. Thus,   condition \eqref{gamma}  is fulfilled
if
\[\sum_{k=1}^\infty \frac{\la_k^2\,k^\a}{\a_k^{1-\theta}} <\infty.\]

\item In case
\[\a_k\sim k^{2/d},\ \ \ \ \ \ k \in\,\nat,\]
(and this is true for several reasonable domains) \eqref{gamma} becomes
\[\sum_{k=1}^\infty \frac{\lambda_k^2}{k^{2(1-\theta)/d}}|e_k|_\infty<+\infty.\]
Thus, in dimension $d=1$ condition \eqref{gamma} is fulfilled by
a white noise in space and time. As soon as one goes to higher dimension, this of
course is no longer possible. It is important to stress, however,  that if the
sup-norms of the eigenfunctions $e_k$ are equi-bounded, condition \eqref{gamma}  does not require a noise with trace-class covariance, no matter how large the space dimension is.

\end{enumerate}}
\end{Remark}

Concerning the nonlinearity $b$, we shall assume the following condition

\begin{Hypothesis}
\label{H2} The mapping $b:\bar{\mathcal{O}}\times \reals\to
\reals$  is measurable and
\[\sup_{x \in\,\bar{\mathcal{O}}}|b(x,\si)-b(x,\rho)|\leq
L\,|\si-\rho|,\ \ \ \ \ \ \si,\,\rho \in\,\reals,\] for some
positive constant $L$. Moreover
\[\sup_{x \in\,\bar{\mathcal{O}}}|b(x,0)|=:b_0<\infty.\]
\end{Hypothesis}

\subsection{Estimates for the stochastic convolution}
\label{subsec3.1}

For each $\mu>0$, we consider  the linear problem
\begin{equation}
\label{linear} \le\{\begin{array}{l} \ds{\mu \frac{\partial^2
\eta}{\partial t^2}(t,x)=\Delta \eta(t,x)-\frac{\partial
\eta}{\partial t}(t,x)+ \frac{\partial w^Q}{\partial t}(t,x),\ \ \ \ \ t>0,\ \ x \in\,\mathcal{O},}\\
\vs \ds{\eta(0)=0,\ \ \ \ \ \frac{\partial \eta}{\partial
t}(0)=0,\ \ \ \ \ \ \eta(t,x)=0, \ \ \ \ t\geq 0,\ \ x
\in\,\partial \mathcal{O}.} \end{array}\r.
\end{equation}

It is well known  that if for some $\theta \in\,\reals$
condition \eqref{gamma} holds, then for any $\mu>0$ there exists a
unique solution $\eta^\mu$ to problem \eqref{linear} such that for
any  $T>0$ and $p\geq 1$
\[
 \eta^\mu
\in\,L^p(\Omega;C([0,T];H^\theta(\mathcal{O}))),\ \ \ \ \ \
\frac{\partial \eta^\mu}{\partial t}
\in\,L^p(\Omega;C([0,T];H^{\theta-1}(\mathcal{O})))\]
(for a proof we refer for example to \cite{dpz1} and \cite{dpb}).

Our aim here is proving that if the constant $\theta$ above is
strictly positive (as in Hypothesis \ref{H1}), then for any
$\d<\theta/2$ the process $\eta^\mu$ has a version which is
$\delta$-H\"older continuous with respect to $t\geq 0$ and $\xi
\in\,\bar{\mathcal{O}}$ and the momenta of the $\d$-H\"older norms
of $\eta^\mu$ are equi-bounded with respect to $\mu>0$. Namely we
prove the following result.

\begin{Proposition}\cite[Proposition 3.1]{smolu2}
\label{teorema1}
 Assume that Hypothesis \ref{H1} is satisfied.
 Then  for any $\mu>0$ and $\d <\theta/2$ the
process $\eta^\mu$ has a version (which we still denote by
$\eta^\mu$) which is $\delta$-H\"older continuous with respect to
$(t,x) \in\,[0,T]\times \bar{\mathcal{O}}$, for any $T>0$.

Moreover, for any $p\geq 1$
\[ \sup_{\mu>0}\E\,|\eta^\mu|^p_{C^\d([0,T]\times
\bar{\mathcal{O}})}=:c_{T,p}<\infty.
\]
\end{Proposition}

Due to \eqref{covariance} and Hypothesis \ref{H1}, the  cylindrical Wiener process  $w^Q(t,x)$ can be written as
\[w^Q(t,x)=\sum_{k=1}^\infty Q e_k\,\beta_k(t)=
\sum_{k=1}^\infty \la_k e_k(x)\,\beta_k(t),\ \ \ \ \ t\geq 0,\ \ \ \ x \in\,\mathcal{O},\] where
$\{e_k\}_{k \in\,\nat}$ is the complete orthonormal basis in
$H$ which diagonalizes $\Delta$ and
$\{\beta_k(t)\}_{k \in\,\nat}$ is a sequence of mutually
independent standard Brownian motions all defined on some
stochastic basis $(\Omega, \mathcal{F}, \mathcal{F}_t,
\mathbb{P})$. This means that   for all $(t,x) \in\,[0,\infty)\times
\bar{\mathcal{O}}$ we have
\[ \eta^\mu(t,x)=\sum_{k=1}^\infty
\eta_k^\mu(t)e_k(x),\]
 where, for each $k \in\,\nat$,
$\eta^\mu_k(t)$ is the solution of the one dimensional problem
\begin{equation}
\label{unidimensional} \le\{\begin{array}{l} \ds{\mu \frac{d^{\,2}
\eta_k^{\mu}}{dt^2}(t)=-\a_k \eta^\mu_k(t)-
\frac{d \eta_k^{\mu}}{dt} (t)+\la_k \frac{d \beta_k}{dt}(t),}\\
\vs \ds{\eta_k^\mu(0)=0,\ \ \ \ \ \frac{d \eta_k^{\mu}}{dt}(0)=0.}
\end{array}\r.
\end{equation}
Notice that the second order equation \eqref{unidimensional} can
be rewritten more rigorously as the following system
\[\le\{
\begin{array}{l}
\ds{d\eta_k^\mu(t)=\theta_k^\mu(t)\,dt}\\
\vs \ds{\mu\, d\theta_k^{\mu}(t)=-\le(\a_k \eta^\mu_k(t)+
\theta_k^{\mu}(t)\r)\,dt+\la_k\, d\beta_k(t),}\\
\vs \ds{\eta_k^\mu(0)=0,\ \ \ \ \ \theta_k^{\mu}(0)=0.}
\end{array}\r.
\]
Then, by the variation of constants formula, it is immediate to
check that
\[ \eta_k^\mu(t)=\frac {\la_k}{\mu}\int_0^t
f^\mu_k(t-s;0,e_k)\,d\beta_k(s),
\]
and
\[ \frac{d \eta_k^\mu}{d t}(t)=\theta^\mu_k(t)=\frac
{\la_k}{\mu}\int_0^t g^\mu_k(t-s;0,e_k)\,d\beta_k(s),
\] with $f_k^\mu$ and $g^\mu_k$ defined as the solutions of
 system \eqref{sistema}, with initial conditions $f(0)=0$ and
 $g(0)=1$.

\bigskip

Our aim here  is obtaining estimates for $\eta^\mu(t)$ which are
independent of $\mu>0$. We start with a uniform estimate for the mean-square of
$\eta^\mu(t,x)-\eta^\mu(s,y)$, for each $t,s\geq 0$ and $x, y
\in\,\bar{\mathcal{O}}$.

\begin{Lemma}\cite[Lemma 3.2]{smolu2}
\label{lemma1}
Under Hypothesis \ref{H1} there exists a constant $c_1>0$ such
that
\begin{equation}
\label{fr7}
 \sup_{\mu>0} \E\,|\eta^\mu(t,x)-\eta^\mu(t,y)|^2\leq c_1\,|x-y|^{2\theta},
\end{equation}
for any $t\geq 0$ and $x,y \in\,\bar{\mathcal{O}}$.

Moreover, there exists a constant
$c_2>0$ such that
\begin{equation}
\label{fr8} \sup_{\mu>0}\E\,|\eta^\mu(t,x)-\eta^\mu(s,x)|^2\leq
c_2|t-s|^\theta,
\end{equation}
for any $t,s\geq 0$ and $x \in\,\bar{\mathcal{O}}$.
\end{Lemma}

Once we have the uniform estimates \eqref{fr7} and \eqref{fr8}, the
proof of Proposition \ref{teorema1} is straightforward.
Actually, since for any $t,s\geq
0$ and $x,y \in\,\bar{\mathcal{O}}$ the random variable
$\eta^\mu(t,x)-\eta^\mu(s,y)$ is Gaussian, according to Lemma
\ref{lemma1}  for any $p\geq 1$ there exists
a constant $c_p$ such that
\[ \E\,|\eta^\mu(t,x)-\eta^\mu(s,y)|^p\leq c_p
\le(|t-s|^2+|x-y|^2\r)^{\frac{\theta p}4}.\] Now, if we fix any
$\d \in\,(0,\theta/2)$, there exists $p_\d\geq 1$ such that
\[\frac \theta 2-\frac {d+1}{p_\d}>\d.\]
Thus, if for any $p\geq p_\d$ we define
\[\a:=\d+\frac{d+1}{p},\] we have
\[\begin{array}{l}
\ds{\E\int_{[0,T]^2\times
\mathcal{O}^2}\frac{|\eta^\mu(t,x)-\eta^\mu(s,y)|^{p}}
{\le(|t-s|^2+|x-y|^2\r)^{\frac{d+1+\a p}2}}\,dt\, ds\, dx\, dy}\\
\vs \ds{\leq c_{\d}\int_{[0,T]^2\times
\mathcal{O}^2}\le(|t-s|^2+|x-y|^2\r)^{\le(\frac{\theta
p}{4}-\frac{d+1+\a p}{2}\r)}\,dt\, ds\, dx\, dy}\\
\vs \ds{= c_{\d}\int_{[0,T]^2\times
\mathcal{O}^2}\le(|t-s|^2+|x-y|^2\r)^{\le(\frac{\theta
p}{4}-\frac{\d p}{2}-(d+1)\r)}\,dt\, ds\, dx\,
dy=:c^\prime_\d<\infty,}
\end{array}
\]
as
\[2\le(\frac{\theta
p}{4}-\frac{\d p}{2}-(d+1)\r)\geq p\le(\frac \theta
2-\frac{d+1}{p_\d}-\d\r)-(d+1)>-(d+1).\] This implies that
$\eta^\mu$ belongs to $W^{\a,p}([0,T]\times \mathcal{O})$,
$\mathbb{P}$-a.s., so that, due to the Sobolev embedding theorem,
there exists a version of $\eta^\mu$ belonging to
$C^\d([0,T]\times \bar{\mathcal{O}})$, $\mathbb{P}$-a.s. Moreover,
there exists some $c_{T,p}$ independent of $\mu>0$ such that
\[\E\,|\eta^\mu|_{C^\d([0,T]\times \bar{\mathcal{O}})}^p\leq c\,\E\,
|\eta^\mu|_{W^{\a,p_\d}([0,T]\times \bar{\mathcal{O}})}^p\leq
c_{T,p}.\]

\subsection{The convergence result}
\label{subsec3.2}

According to
Proposition \ref{propy} for any constant $c>0$ and any $k
\in\,\nat$
\[\lim_{t\to \infty}\sup_{|u_k|+|v_k|\leq
c}|f^\mu_k(t;u,v)|=\lim_{t\to \infty}\sup_{|u_k|+|v_k|\leq
c}|g^\mu_k(t;u,v)|=0,\] so that for any $k \in\,\nat$ we have
\[ \lim_{t\to
\infty}\le(\sup_{|(u,v)|_{\mathcal{H}_\delta}\leq 1}\a_k^\delta\,
[\Pi_1 S_\mu(t)(u,v)]_k^2 +\sup_{|(u,v)|_{\mathcal{H}_\delta}\leq
1}\a_k^{\delta-1}\, [\Pi_2 S_\mu(t)(u,v)]_k^2\r)=0.
\]

Now, if we multiply the second equation in \eqref{sistema} by
$g^\mu_k(t)$ we get
\[\mu \frac{d\,|g^\mu_k(\cdot;u,v)|^2}{dt}(t)+\a_k\frac{d\,|f^\mu_k(\cdot;u,v)|^2}{dt}(t)+2\,|g^\mu_k(t;u,v)|^2=0,\]
and hence, integrating with respect to $t\geq 0$ we get
\[\mu\,|g^\mu_k(t;u,v)|^2+\a_k\,|f^\mu_k(t;u,v)|^2+2\int_0^t
|g^\mu_k(s;u,v)|^2\,ds= \mu\,|v_{k}|^2+\a_k\,|u_{k}|^2.
\]
This implies that
 the families of functions
\[\le\{\Pi_1 S_\mu(\cdot)(u_0,v_0)\r\}_{\mu \in\,(0,1)}\subset L^\infty(0,\infty
;H_0^{1}(\mathcal{O})),\ \ \ \ \ \le\{\Pi_2
S_\mu(\cdot)(u_0,v_0)\r\}_{\mu \in\,(0,1)}\subset L^2(0,\infty
;H)\] are equi-bounded, and by the
Ascoli-Arzel\`a theorem we have
\[ \le\{\Pi_1
S_\mu(\cdot)(u_0,v_0)\r\}_{\mu \in\,(0,1)}\subset
C([0,T];H),\ \ \ \ \ \text{compactly},
\] for any $T>0$.

Now, for any $\mu>0$ and $\delta \in\,[0,1]$, we define the
operators
\begin{equation}
\label{pp1} B_\mu(h,k)(x):=\frac 1\mu(0,b(x,h(x))),\ \ \ (h,k)
\in\,\mathcal{H}_\delta,\ \ \ x \in\,\O,\end{equation}
 and
\begin{equation}
\label{pp2} \mathcal{Q}_\mu h=\frac 1\mu(0, Qh),\ \ \ \ \ h
\in\,H.\end{equation} Note that, since
$\delta\in\,[0,1]$, for any $z_1=(u_1,v_1)$ and $z_2=(u_2,v_2)
\in\,\mathcal{H}_\delta$
\[|B_\mu(z_1)-B_\mu(z_2)|_{\mathcal{H}_\delta}=\frac
1{\mu}\,|b(\cdot,u_1)-b(\cdot,u_2)|_{H^{\delta-1}(\mathcal{O})}\leq
\frac c{\mu}\,|b(\cdot,u_1)-b(\cdot,u_2)|_{H},\]
and then, thanks to Hypothesis \ref{H2}
\[ |B_\mu(z_1)-B_\mu(z_2)|_{\mathcal{H}_\delta}\leq \frac
{c\,L}{\mu}|u_1-u_2|_{H}\leq \frac
{c\,L}{\mu}|z_1-z_2|_{\mathcal{H}_\delta}.
\]

\begin{Definition}
\label{def4.1}
Let $\delta\in\,[0,1]$. A process $z^\mu(t)=(u^\mu(t),v^\mu(t))$, $t\geq 0$,  is a {\em mild solution}
of problem \eqref{quasilinear} in $L^2(\Omega;C([0,T];\mathcal{H}_\d)$, if
\[u^\mu \in\,L^2(\Omega;C([0,T];H^\delta(\mathcal{O}))),\ \ \ \
v^\mu=\frac{\partial u^\mu}{\partial t}
\in\,L^2(\Omega;C([0,T];H^{\delta-1}(\mathcal{O}))),\]  for any
$T>0$, and
\[z^\mu(t)=S_\mu(t)(u_0,v_0)+ \int_0^t
S_\mu(t-s)B_\mu(z^\mu(s))\,ds+\int_0^t
S_\mu(t-s)\,dw^{\mathcal{Q}_\mu}(s).\]
\end{Definition}

Note that with these notations,  the weak solution $\eta^\mu$ of
problem \eqref{linear} studied in Subsection \ref{subsec3.1} can be
written as
\[\eta^\mu(t)=\Pi_1\int_0^t
\,S_\mu(t-s)\,dw^{\mathcal{Q}_\mu}(s),\ \ \ \ t\geq 0,\]
 and hence, if
$z^\mu=(u^\mu,v^\mu)$ is a mild solution of problem \eqref{quasilinear}, we have
\[
u^\mu(t)=\Pi_1S_\mu(t)(u_0,v_0)+\Pi_1\int_0^t
S_\mu(t-s)B_\mu(u^\mu(s),v^\mu(s))\,ds+\eta^\mu(t),\ \ \ \ \ t\geq
0. \]

Due to the Lipschitz continuity of $b$, and hence of $B_\mu$, it is possible to prove the following well-posedness result for problem \eqref{quasilinear}, for every $\mu>0$ fixed
 (for a proof see e.g. \cite[Theorem
5.3.1]{dpz2}).

\begin{Proposition}
\label{esistenza} Assume Hypotheses  \ref{H1} and \ref{H2}. Then
for any $\mu>0$  and for any initial data $u_0
\in\,H^\theta(\mathcal{O})$ and $v_0
\in\,H^{\theta-1}(\mathcal{O})$, there exists a unique mild
solution $z^\mu(t)$  to problem \eqref{quasilinear}. Moreover,
 for any $T>0$ and $p\geq 1$ there exists $c_{p,\mu}(T)>0$ such that
\begin{equation}
\label{sl5} \E\,\sup_{t
\in\,[0,T]}\,|z^\mu(t)|_{\mathcal{H}_\theta}^p\leq
c_{p,\mu}(T)\le(1+|(u_0,v_0)|_{\mathcal{H}_\theta}^p\r).
\end{equation}

\end{Proposition}

Once one has the well-posedness of equation \eqref{quasilinear}, for every $\mu>0$, in order to prove the validity of the Smoluchowski-Kramers approximation we have to show that
 the family of probability measures
$\{\L(u^\mu)\}_{\mu \in\,[0,1]}$ is  tight on
$C([0,T];L^2(\mathcal{O}))$, for any $T>0$.

\begin{Proposition}\cite[Proposition 4.3]{smolu2}
\label{tight} Assume that $u_0 \in\,H^1(\mathcal{O})$ and
$v_0\in\,L^2(\mathcal{O})$. Then, under Hypotheses \ref{H1} and
\ref{H2} the family of probability measures $\{\L(u^\mu)\}_{\mu
\in\,[0,1]}$ is tight on the space $C([0,T];H)$, for any
$T>0$.
\end{Proposition}

This follows from Lemma \ref{lemma1}. Once one has the tightness of the family $\{\L(u^\mu)\}_{\mu
\in\,[0,1]}$ on $C([0,T];H)$, by the Prokhorov theorem for every sequence $\{\L(u^{\mu_n})\}_{n \in\,\nat}$ there exists a subsequence $\{\L(u^{\mu_{n_k}})\}_{k \in\,\nat}$ and a probability $\mathbb{Q}$ on $C([0,T];H)$ such that $\L(u^{\mu_{n_k}})$ weakly converges to $\mathbb{Q}$.

The final step of our proof consists in identifying the probability $\mathbb{Q}$ with $\mathcal{L}(u)$ and showing that in fact the family $u^\mu$ converges to $u$ in $C([0,T];H)$ in the sense of convergence in probability.

To this purpose, the following integration by parts formula holds.

\begin{Lemma} \cite[Lemma 4.4]{smolu2}
\label{lemma3} Assume Hypotheses \ref{H1} and \ref{H2} and fix
$u_0 \in\,H^\theta(\mathcal{O})$ and $v_0
\in\,H^{\theta-1}(\mathcal{O})$. Then for any $\mu>0$ and for any
$\varphi \in\,C^{2}([0,T]\times \bar{\mathcal{O}})$, such that
$\varphi\equiv 0$ on $\partial \mathcal{O}$, we have
\begin{equation}
\label{fr64}
\begin{array}{l}
\ds{\int_{\mathcal{O}}
u^\mu(t,x)\varphi(t,x)\,dx=\int_{\mathcal{O}}
u_0(x)\varphi(0,x)\,dx+\int_0^t
\int_{\mathcal{O}}u^\mu(s,x)\le[\frac{\partial \varphi}{\partial
t}(s,x)+\Delta \varphi(s,x)\r]\,ds\,dx}\\
\vs \ds{+\int_0^t
\int_{\mathcal{O}}b(x,u^\mu(s,x))\varphi(s,x)\,ds\,dx+\int_0^t
\int_{\mathcal{O}}\varphi(s,x)\,w^Q(ds,dx)+R_\mu(t),}
\end{array}
\end{equation}
where
\begin{equation}
\label{fr53} \begin{array}{l} \ds{R_\mu(t):=\mu\le(1-e^{-\frac
t\mu}\r)\int_{\mathcal{O}}v_0(x)\varphi(0,x)\,dx-\int_0^t
 e^{-\frac{t-s}\mu}\,M_\mu(s)\,ds}\\
 \vs
 \ds{-\int_0^t
 e^{-\frac{t-s}\mu}\le[\int_{\mathcal{O}}\le(u_0(x)\frac{\partial \varphi}{\partial t}
 (0,x)-u^\mu(s,x)\frac{\partial \varphi}{\partial t}(s,x)+
 \int_0^s u^\mu(r,x)\frac{\partial
 ^2\varphi}{\partial
 t^2}(r,x)\,dr\r)\,dx\r]\,ds}\\
 \vs
 \ds{
 -\int_0^t \int_{\mathcal{O}}e^{-\frac{t-s}\mu}\varphi(s,x)w^Q(ds,dx),}
 \end{array}
 \end{equation}
 and
 \[
M_\mu(t):=
\int_\mathcal{O}\le(u^\mu(t,x)\Delta\varphi(t,x)+b(x,u^\mu(t,x))\varphi(t,x)\r)\,dx.
\]
\end{Lemma}

Concerning the remainder term $R_\mu(t)$ defined in \eqref{fr53}
we have the following limiting result.

\begin{Lemma}
\label{lemma4} Under the same hypotheses of Lemma \ref{lemma3} we
have
\[\lim_{\mu \to 0}\E\,|R_\mu(t)|^2=0,\ \ \ \ \ \ t\geq 0.\]
\end{Lemma}

Due to the tightness in $C([0,1];H)$ of the
sequence $\{\L(u^\mu)\}_{\mu \in\,[0,1]}$,  the Skorokhod theorem
assures that for any two sequences $\{\mu_n\}_n$ and $\{\mu_m\}_m$
converging to zero there exist  subsequences $\{\mu_{n(k)}\}_{k
\in\,\nat}$ and $\{\mu_{m(k)}\}_{k \in\,\nat}$ and a sequence of
random elements
\[\{\rho_k\}_{k \in\,\nat}:=\le\{(u_1^k,u_2^k,\hat{w_k}^Q)\r\}_{k \in\,\nat},\]
in $\mathcal{C}:=C([0,T];H)^2\times
C([0,T];\mathcal{D}^\prime(\mathcal{O}))$, defined on some
probability space $(\hat{\Omega},\hat{\F},\hat{\Pro})$, such that
the law of $\rho_k$ coincides with the law of
$(u^{\mu_{n(k)}},u^{\mu_{m(k)}},w^Q)$, for each $k \in\,\nat$, and
$\rho_k$ converges $\hat{\Pro}$-a.s. to some random element
$\rho:=(u_1,u_2,\hat{w}^Q) \in\,\mathcal{C}$.

Now, if we show that $u_1=u_2$, we have that there exists some $u
\in\,C([0,T];H)$ such that $u^\mu$ converges to $u$
in probability. Actually, as observed by Gy\"ongy and Krylov in
\cite{gk}, if $E$ is any Polish space equipped with the Borel
$\si$-algebra, a sequence $\{\rho_n\}$ of $E$-valued random
variables converges in probability if and only if for every pair
of subsequences $\{\rho_m\}$ and $\{\rho_l\}$ there exists an
$E^2$-valued subsequence $w_k:=(\rho_{m(k)},\rho_{l(k)})$
converging weakly to a random variable $w$ supported on the
diagonal $\{(h,k) \in\,E^2\ :\ h=k\}$.

Note that both $u_1^k$ and $u_2^k$ solve equation
\eqref{quasilinear} with $w^Q$ replaced by $\hat{w_k}^Q$. Then
they both verify formula \eqref{fr64},   with $R_1^k$ and $R_2^k$
obtained replacing $u^\mu$ respectively with $u_1^k$ and $u_2^k$
and $w^Q$ with $\hat{w}^Q_k$. According to Lemma \ref{lemma4} we
have that both $R_1^k$ and $R_2^k$ converge to zero in
$L^2(\hat{\Omega})$, as $m_{n(k)}$ and $\mu_{m(k)}$ go to zero,
and then, possibly for a subsequence, they converge
$\hat{\Pro}$-a.s. to zero. Due to formula \eqref{fr64} this
implies
\[\begin{array}{l}
\ds{\int_{\mathcal{O}} u_i(t,x)\varphi(t,x)\,dx=\int_{\mathcal{O}}
u_0(x)\varphi(0,x)\,dx+\int_0^t
\int_{\mathcal{O}}u_i(t,x)\le[\frac{\partial \varphi}{\partial
t}(s,x)+\Delta \varphi(s,x)\r]\,ds\,dx}\\
\vs \ds{+\int_0^t
\int_{\mathcal{O}}b(x,u_i(s,x))\varphi(s,x)\,ds\,dx+\int_0^t
\int_{\mathcal{O}}\varphi(s,x)\,\hat{w}^Q(ds,dx),\ \ \ \ \ \
i=1,2,}
\end{array}
\] and then they coincide with the solution of the semi-linear
heat equation perturbed by the noise $\hat{w}^Q$, which is unique.

As we have recalled above, thanks to the Gy\"ongy-Krylov remark in
\cite{gk} this implies that $u^\mu$ converges in probability to
some random variable  $u \in\,C([0,T];H)$. But, by
using again formula \eqref{fr64} and Lemma \ref{lemma4} we have
that $u$ solves the heat equation \eqref{heat}.

We have just proved the main result of this section.

\begin{Theorem}
\label{teo4.6}
Assume Hypotheses  \ref{H1} and \ref{H2} and fix $u_0
\in\,H^1(\mathcal{O})$ and $v_0 \in\,L^2(\mathcal{O})$. Then, if
$u^\mu$ is the solution of the stochastic semi-linear damped wave
equation \eqref{quasilinear} and $u$ is the solution of the
stochastic semi-linear heat equation \eqref{heat}, for any $T>0$
and for any $\e>0$ we have
\[\lim_{\mu\to 0}
\,\Pro\le(|u^\mu-u|_{C([0,T];H)}>\e\r)=0.
\]
\end{Theorem}

\section{The approximation in the case of multiplicative noise}
\label{sec4}

We are dealing  here with the following  damped semi-linear wave
equation perturbed by multiplicative noise in the interval $[0,L]$
\begin{equation}
\label{ieq4} \le\{\begin{array}{l} \ds{\mu \frac{\partial^2
u}{\partial t^2}(t,x)+\frac{\partial u}{\partial t}(t,x)=\frac{
\partial ^2u}{\partial x^2}(t,x)+b(x,u(t,x))+ g(x,u(t,x))
\frac{\partial w}{\partial t}(t,x),\ \ x \in\,[0,L],}\\
\vs \ds{u(0,x)=u_0(x),\ \ \ \ \ \frac{\partial u}{\partial
t}(0,x)=v_0(x),\ \ \ \ \ \ u(t,0)=u(t,L)=0.} \end{array}\r.
\end{equation}
In the present section,  we assume that $w(t,x)$ is a cylindrical Wiener process, white in space and time. The   coefficients  $b$ and $g$ are measurable from
$[0,L]\times \reals$ with values in $\reals$ and $g$ is
Lipschitz-continuous in the second variable, uniformly with
respect the first one. The coefficient $b$ is either assumed to be
Lipschitz-continuous in the second variable, uniformly with
respect the first one, or satisfying some polynomial growth and
dissipation conditions in the spirit of the Klein-Gordon model
(and in this second case $g$ is assumed bounded).

As in  the case of Lipschitz continuous $b$ and
additive noise in any space dimension, considered in Section \ref{sec3},   we want to show that  for any $\e>0$ and $T>0$
\[ \lim_{\mu\to 0}\mathbb{P}\le(\,\sup_{t
\in\,[0,T]}\le|u^\mu(t)-u(t)\r|_{L^2(0,L)}>\e\,\r)=0,\]
where
$u(t)$ is  the solution of the parabolic problem
\begin{equation}
\label{ieq5} \le\{\begin{array}{l} \ds{\frac{\partial u}{\partial
t}(t,x)=\Delta u(t,x)+b(x,u(t,x))+ g(x,u(t,x))\frac{\partial
w}{\partial t}(t,x),\ \ \ \
x \in\,[0,L],}\\
\vs \ds{u(0,x)=u_0(x),\ \ \ \ \ \ u(t,0)=u(t,L)=0.} \end{array}\r.
\end{equation}

\subsection{The coefficients $b$ and $g$}
\label{subsec4.1}

Concerning the  coefficients $b$ and   $g$, as we already mentioned above,we shall
consider two different cases. Here is described the first one.

\begin{Hypothesis}
\label{H3} \begin{enumerate} \item The mapping  $b:[0,L]\times
\reals\to \reals$ is measurable  and  there exists $M>0$ such that
\[\sup_{x \in\,[0,L]}|b(x,\si)-b(x,\rho)|\leq
M\,|\si-\rho|,\] for any $\si,\rho \in\,\reals$. Moreover,
$\sup_{x \in\,[0,L]}|b(x,0)|=:b_0<\infty$.

\item The mapping  $g:[0,L]\times \reals\to \reals$ is measurable
 and  there exists  $M>0$ such
that
\[\sup_{x \in\,[0,L]}|g(x,\si)-g(x,\rho)|\leq
M\,|\si-\rho|,\] for any $\si,\rho \in\,\reals$. Moreover,
$\sup_{x \in\,[0,L]}|g(x,0)|=:g_0<\infty$.
\end{enumerate}

\end{Hypothesis}

In particular from the assumptions above we have that both $b$ and
$g$ have linear growth in the second variable, uniformly with
respect to the first. Namely
\[  \sup_{x \in\,[0,L]}|b(x,\si)|\leq c\,\le(1+|\si|\r),\
\ \ \ \ \sup_{x \in\,[0,L]}|g(x,\si)|\leq c\,\le(1+|\si|\r),
\]
for some constant $c>0$.

As we are assuming
that $b(x,\cdot)$ is Lipschitz continuous, uniformly with respect
to $x \in\,[0,L]$, we have seen that the mapping $B_\mu$ defined in \eqref{pp1} is Lipschitz continuous  from $\mathcal{H}_\d$ into itself.

Now, for any $\mu>0$ and $\d \in\,[0,1]$, we define
\[[G_\mu(u,v)h](x):=\frac 1\mu\,(0,g(x,u(x)))h(x),\ \ \ x
\in\,[0,L],\ \ \ (u,v) \in\,\mathcal{H}_\d,\ \ \ h
\in\,L^\infty(0,L).\]
   Due to Hypothesis \ref{H3},  the mapping
$G_\mu(\cdot)h:\mathcal{H}_\d\to \mathcal{H}_\d$ is Lipschitz
continuous, for any fixed $h \in\,L^\infty(0,L)$. Actually, as $\d
\in\,[0,1]$, we have
\begin{equation}
\label{ph5}
\begin{array}{l}
\ds{|G_\mu(z_1)h-G_\mu(z_2)h|_{\mathcal{H}_\d}=\frac
1\mu\,|(g(u_1)-g(u_2))h|_{H^{\d-1}}
\leq \frac
c\mu\,|(g(u_1)-g(u_2))h|_{H}}\\
\vs \ds{\leq \frac{c M}\mu\,|u_1-u_2|_H\,|h|_\infty\leq \frac{c
M}\mu\,|z_1-z_2|_{\mathcal{H}_\d}\,|h|_\infty.}
\end{array}
\end{equation}

The second case that we shall consider  is
described below.

\begin{Hypothesis}
\label{H4} \begin{enumerate} \item The mapping  $b:[0,L]\times
\reals \to \reals$ is measurable  and $b(x,\cdot):\reals\to
\reals$ is of class $C^1$, for almost all $x \in\,[0,L]$.
Moreover,

\begin{enumerate}
\item there exist $\la \in\,(1,3]$ and
$c_1>0$ such that for any $\si \in\ [0,L]$
\begin{equation}
\label{poli} \sup_{x \in\,[0,L]}|b(x,\si)|\leq
c_1\,\le(1+|\si|^\la\r),\ \ \ \ \ \sup_{x
\in\,[0,L]}\le|\partial_\si\,b (x,\si)\r|\leq
c_1\,\le(1+|\si|^{\la-1}\r);
\end{equation}
\item there exists $c_2>0$ such that for any $\si
\in\,\mathbb{R}$
  \begin{equation}
  \label{dissi}
\sup_{x \in\,[0,L]}\int_0^\si b(x,\rho)\,d\rho\leq \,
c_2\,\le(1-|\si|^{\la+1}\r);
  \end{equation}
\item for any $(x,\si)
\in\,[0,L]\times \mathbb{R}$
\begin{equation}
\label{dissibis} \partial_\si\, b(x,\si)\leq 0.
\end{equation}
\end{enumerate}

\item The mapping  $g:[0,L]\times \reals\to \reals$ is measurable and
there exists  $M>0$ such that
\[\sup_{x \in\,[0,L]}|g(x,\si)-g(x,\rho)|\leq
M\,|\si-\rho|,\] for any $\si,\rho \in\,\reals$. Moreover,
\begin{equation}
\label{glim} \sup_{(x,\si) \in\,[0,L]\times
\mathbb{R}}|g(x,\si)|=:g_0<\infty.
\end{equation}

\end{enumerate}
\end{Hypothesis}

\begin{Remark}
{\em \begin{enumerate} \item A typical example of a function $b$
fulfilling conditions \eqref{poli}, \eqref{dissi} and
\eqref{dissibis} is
\[b(x,\si)=-\a\,|\si|^{\la-1}\si,\]
 for any strictly positive
constant $\a$ (the Klein-Gordon equation).

\item Here we are assuming  the condition $\partial_\si\,b(x,\si)\leq 0$, for any
$(x,\si) \in\,[0,L]\times \reals$, just for simplicity of
notations. Actually we could also treat the case
\[\sup_{(x,\si) \in\,[0,L]\times
\reals}\partial_\si\,b(x,\si)\leq c,\] for some constant $c$, by
setting $b(x,\si)=b_1(x,\si)+b_2(x,\si)$, where
$b_1(x,\si)=b(x,\si)-c\,\si$ fulfills conditions \eqref{poli},
\eqref{dissi} and \eqref{dissibis}, and $b_2(x,\si)=c\,\si$ is a
Lipschitz perturbation.
\end{enumerate}}
\end{Remark}

In this second case the mapping $b(x,\cdot):\reals \to \reals$ is
no more Lipschitz continuous and does not necessarily have sublinear
 growth. From  \eqref{dissibis}
we obtain
\[\si b(x,\si)\leq \int_0^\si b(x,\rho)\,d\rho,\ \ \ \ (x,\si)
\in\,[0,L]\times \reals,\]
 and then according to \eqref{dissi} it
follows
 \begin{equation}
 \label{was5}
 \lim_{|\si|\to \infty} \si b(x,\si)=-\infty,\ \ \ \ x \in\,[0,L].
 \end{equation}
 Moreover, due to \eqref{dissi}, for any $(x,\si)
\in\,[0,L]\times \mathbb{R}$ we have
 \begin{equation}
 \label{was6}
 |\si|^{\la+1}\leq 1-\frac 1{c_2}
\int_0^\si b(x,\rho)\,d\rho=:-\beta(x,\si).
\end{equation}

Now, by proceeding as in \cite{millet2}, where the
 particular case of the Klein-Gordon equation is considered,
 we approximate $b$ by means of Lipschitz continuous
 mappings, by setting for any $n
 \in\,\nat$
 \[
 b_n(x,\si):=\le\{ \begin{array}{ll}
 \ds{b(x,n)+(\si-n)\,\partial_\si\, b(x,n),} & \ds{\text{if}\ \si\geq n,}\\
 & \vs
 \ds{b(x,\si),}
 & \ds{\text{if}\ |\si|\leq n,}\\
 & \vs
 \ds{b(x,-n)+(\si+n)\,\partial_\si\, b(x,-n),} & \ds{\text{if}\ \si\leq -n.}
 \end{array}\r.
 \]
Clearly,  $b_n(x,\cdot):\reals\to \reals$ is
 Lipschitz continuous, uniformly with respect to $x \in\,[0,L]$,
 and
 $b_n(x,\si)\equiv b(x,\si)$ on  $[0,L]\times [-n,n]$, so that
the mapping $B_{n,\mu}$ defined by
\[ B_{n,\mu}(u,v)(x)=\frac 1\mu\,(0, b_n(x,u(x))),\ \ \
\ \ x \in\,[0,L],\ \ \ (u,v) \in\,\mathcal{H}_\d,
\]
is Lipschitz continuous  from $\mathcal{H}_\d$ into itself, for
any $\d \in\,[0,1]$, and for any $(u,v) \in\,\mathcal{H}_\d$ it
holds
\[ |u|_\infty\leq n\Rightarrow
B_{n,\mu}(u,v)=B_\mu(u,v).
\]
Moreover,   for any $\si \in\,\reals$
 \[ \sup_{x \in\,[0,L]}|b_n(x,\si)|\leq
c\,\le(1+|\si|^\la\r),\ \ \ \ \ \sup_{x
\in\,[0,L]}\le|\partial_\si\,b_n (x,\si)\r|\leq
c\,\le(1+|\si|^{\la-1}\r),
\]
for some constant $c$ not depending on $n$.

Next, for any $n \in\,\mathbb{N}$ we  set
 \[
 {\beta}_n(x,\si):=\frac 1 c_2
\int_0^\si b_n(x,\rho)\,d\rho-1.
\]
Due to \eqref{dissibis}, \eqref{was5} and \eqref{was6} it
is possible to show that there exists some $n_0>0$ such that for
any $n\geq n_0$ and $(x,\si) \in\,[0,L]\times \mathbb{R}$
\[ -{\beta}_n(x,\si)\geq n^{\la+1}>0.
\]
Now, since $\la\leq 3$, it is possible to  adapt the  arguments
used in \cite[Lemma A1]{millet2} for the Klein-Gordon case and
 prove the following inequality.
\begin{Lemma}\cite[Lemma 2.4]{smolu1}
\label{was2} Assume that $b$ fulfills Hypothesis \ref{H4}-1. Then,
there exists $n_0>0$ such that for any $n\geq n_0$,
 $x
\in\,[0,L]$ and $\si, \rho \in\,\reals$ it holds
 \[
\le|\partial_\si\, b_n (x,\si+\rho)\r|^2\leq
c\,\le(1-{\beta}_n(x,\si)+|\rho|^{2(\la-1)}\r),\] for some
constant $c$ not depending on $n$. In particular, as $b_n$
coincides with $b$ on $[0,L]\times [-n,n]$, for any $n>0$, an
analogous inequality is fulfilled by $b$ and $\beta$ on
$[0,L]\times \reals$.

\end{Lemma}

As far as the diffusion coefficient $G_\mu$ is concerned, due to
\eqref{glim} in this second case the mapping
$G_\mu(\cdot)h:\mathcal{H}_\d\to \mathcal{H}_\d$ is bounded, for
any fixed $h \in\,L^\infty(0,L)$.

\subsection{Uniform bounds for the stochastic convolution}
\label{subsec4.2}

For any $\mu>0$ and $T>0$ and for any
$z\in\,L^p(\Omega;C([0,T];\mathcal{H}_\delta))$ we define
\[ \Gamma_\mu(z)(t):=\int_0^t
S_\mu(t-s)G_\mu(z(s))\,dw(s),\ \ \ \ t \in\,[0,T]. \]

Our aim in this section is proving some a-priori estimates for
$\Gamma_\mu(z)$, which are uniform with respect to $\mu \in(0,1]$.
These estimates represent the key point in the proof of the
tightness of the family of probability measures
$\{\mathcal{L}(u^\mu)\}_{\mu \in\,(0,1]}$ in $C([0,T];H)$ and, as a
consequence, of the Smoluchowski-Kramers approximation (see Section \ref{sec3}).

We first prove uniform bounds for $\Pi_1 \Gamma_\mu(z)$, as a
function of $t \in\,[0,T]$ with values in $H^\d(0,L)$.

\begin{Proposition}\cite[Proposition 3.1]{smolu1}
\label{spazio} Fix  $T>0$. Then for any $p>4$, $\d<1/2-2/p$ and $z
\in\,L^p(\Omega;L^p(0,T;\mathcal{H}))$   we have
\[\sup_{\mu>0}\,E\,|\Pi_1 \Gamma_\mu(z)|_{C([0,T];H^\d(0,L))}^p\leq
c_p(T)\,\le(1+\E\,\int_0^T|\Pi_1 z(t)|^{p}_{H}\,dt\r).
\]

\end{Proposition}

The estimates above  for the $H^\d(0,L)$-norm of $\Pi_1
\Gamma_\mu(z)$, are obtained only for $\d<1/2$. This means that
they do not provide any  bound in the space of continuous
functions. In what follows we shall state some pointwise uniform
bounds, both in the space and in the time variables, which will
lead us  to uniform bounds in the space of H\"older continuous
functions.

\begin{Lemma}\cite[Lemma 3.2]{smolu1}
\label{tempo} Assume that $g:[0,L]\times \mathbb{R}\to \mathbb{R}$
is a measurable mapping such that
\[\sup_{x \in\,[0,L]}\,\le|g(x,\si)\r|\leq
c\,\le(1+|\si|^\kappa\r),\ \ \ \ \si \in\,\mathbb{R},
\]
 for some $\kappa \in\,[0,1]$. Then, for any $\e
\in\,(0,1/2\kappa)$, $\rho \in\,(0,1/2-\kappa(1/2-\kappa \e))$ and
$p\geq 1$ there exists a constant $c=c(\e,\rho,p)$ such that
\[\begin{array}{l}
\ds{
\sup_{\mu>0}\,E\,|\Pi_1\,\Gamma_\mu(z)(t,x)-\Pi_1\,\Gamma_\mu(z)(s,y)|^p}\\
\vs
\ds{\leq
c\,\le(1+\E\,|\Pi_1 z|^{p
\kappa}_{C([0,T];H^{\kappa\e})}\r)\,\le(\,|t-s|^{\frac\rho
2}+|x-y|^{\rho}\r)^p,}
\end{array}
\]
for any $z \in\,L^p(\Omega;C([0,T];\mathcal{H}_{\kappa\e}))$, $x,y
\in\,[0,L]$ and $t,s \,\in\,[0,T]$.

\end{Lemma}

As a consequence of the Garcia-Rademich-Rumsey theorem, from the
previous lemma  we obtain the following result.

\begin{Proposition}\cite[Proposition 3.2]{smolu1}
\label{holder}
 Assume that $g:[0,L]\times \mathbb{R}\to
\mathbb{R}$ is a measurable mapping such that \[\sup_{x
\in\,[0,L]}|g(x,\si)|\leq c\,(1+|\si|^\kappa),\ \ \ \ \ \si
\in\,\mathbb{R},\] for some $\kappa \in\,[0,1]$. Then, for any
$\mu>0$, $\e \in\,(0,1/2\kappa)$, $\rho
\in\,(0,1/4-\kappa(1/2-\kappa \e)/2)$, $p\geq 1$ and $z
\in\,L^p(\Omega;C([0,T];\mathcal{H}_{\kappa\e}))$, the process
$\Pi_1\,\Gamma_\mu(z)$ has a version which is $\rho$-H\"older
continuous with respect to $(t,x) \in\,[0,T]\times [0,L]$.
Moreover, \[
\sup_{\mu>0}\,E\,|\Pi_1\,\Gamma_\mu(z)|^p_{C^\rho([0,T]\times
[0,L])}\leq c\,\le(1+\E\,|\Pi_1 z|^{p
\kappa}_{C([0,T];H^{\kappa\e})}\r),
\]
for some constat $c=c(\e,\rho,p)$.

\end{Proposition}

\begin{Remark}
{\em In particular,  if Hypothesis \ref{H3} is satisfied, then $g$
has linear growth ($\kappa=1$), so that as a consequence of
Proposition \ref{holder} for any $\e \in\,(0,1/2)$, $\rho<\e/2$
and $p\geq 1$ there exists some constant $c=c(\e,\rho,p)$ such
that
\[\sup_{\mu>0}\,E\,|\Pi_1\,\Gamma_\mu(z)|^p_{C^\rho([0,T]\times
[0,L])}\leq c\,\le(1+\E\,|\Pi_1 z|^{p}_{C([0,T];H^{\e})}\r),
\]
 for any  $z \in\,L^p(\Omega;C([0,T];\mathcal{H}_\e))$.

Analogously, if Hypothesis \ref{H4} is verified,  then $g$ is
bounded ($\kappa=0$), so that  for any  $\rho<1/4$ and $p\geq 1$
there exists some constant $c=c(\rho,p)$ such that
\begin{equation}
\label{holder0}
\sup_{\mu>0}\,E\,|\Pi_1\,\Gamma_\mu(z)|^p_{C^\rho([0,T]\times
[0,L])}\leq c, \end{equation} for any  $z
\in\,L^p(\Omega;C([0,T];\mathcal{H}_{0}))$.

 }
\end{Remark}

\subsection{The convergence result}
\label{subsec4.3}

If we set
$z:=\le(u,\partial u/\partial t\r)$  equation \eqref{ieq4} can be
written in the following abstract form
\begin{equation}
\label{astratta} dz(t)=\le[A_\mu
z(t)+B_\mu(z(t))\r]\,dt+G_\mu(z(t))\,dw(t),\ \ \ \
z(0)=z_0=(u_0,v_0).
\end{equation}

As in Definition \ref{def4.1}, a process $z_\mu=(u_\mu,v_\mu) \in\,
L^p(\Omega;C([0,T];\mathcal{H}_\d))$ is a mild solution of \eqref{astratta}, if
\[u^\mu \in\,L^p(\Omega;C([0,T];H^\delta(0,L))),\ \ \ \
v^\mu
\in\,L^p(\Omega;C([0,T];H^{\delta-1}(0,L))),\]   and
\[z^\mu(t)=S_\mu(t)(u_0,v_0)+ \int_0^t
S_\mu(t-s)B_\mu(z^\mu(s))\,ds+\int_0^t
S_\mu(t-s)G_\mu(z^\mu(s))\,dw(s).\]

The existence and uniqueness of a mild solution of equation
\eqref{astratta} for any fixed $\mu>0$ is a well known fact in the
literature, both under Hypothesis \ref{H3} (see \cite{carmona})
and under Hypothesis \ref{H4} (see \cite{millet2} in the more
delicate case of space dimension $d=2$, under more restrictive
conditions on the noise and on the initial data $u_0$ and $v_0$,
due to the higher dimension). Namely, we have

\begin{Theorem}\cite[Theorem 4.2]{smolu1}
\label{prop32}
 Fix   $\mu>0$ and  $\delta \in\,[0,1/2)$ and assume Hypothesis \ref{H1}. Then,
both under Hypothesis \ref{H3} and under Hypothesis \ref{H4}, for
any $T>0$ and $p\geq 1$ and for any initial datum $z_0=(u_0,v_0)
\in\,\mathcal{H}_{1}$ there exists a unique mild solution $z^\mu$
to problem \eqref{astratta} in
$L^p(\Omega;C([0,T];\mathcal{H}_\delta))$.

\end{Theorem}

The proof of the theorem above is obtained by considering an analogous of equation \eqref{ieq4}, obtained by replacing the coefficient $b$ with the truncated coefficients $b_n$. Uniform estimates are obtained for the solutions of the approximating problem and a global solution is obtained by introducing stopping times.The uniform bounds are obtained by using a splitting method and the uniform bounds for the stochastic convolution given in Proposition \ref{holder}.

\begin{Remark}
{\em \begin{enumerate} \item By looking at the proof of the
previous theorem, one sees that, in the case of Lipschitz
continuous $b$, in order to have solutions in
$L^p(\Omega;C([0,T];\mathcal{H}_\d))$ it is not necessary to take
the initial data $z_0=(u_0,v_0)$ in $\mathcal{H}_{1}$, but  it is
sufficient to take them in $\mathcal{H}_\d$.

\item From the proof
we see  that the solution is also unique in
$L^p(\Omega;L^p(0,T;\mathcal{H}_\d))$. \end{enumerate}}
\end{Remark}

\medskip

Now, the key point in the proof of our convergence result is
showing that  the family of probabilities
$\{\mathcal{L}(u^\mu)\}_{\mu \in\,(0,1]}$ is tight in
$C([0,T];H)$.

\begin{Proposition}\cite[Proposition 4.4]{smolu1}
\label{was30} Let $z_0 \in\,\mathcal{H}_1$ and assume that
Hypothesis  \ref{H1} holds. Moreover, assume that \ref{H3} or Hypothesis \ref{H4} hold. Then the family
of probability measures $\{\mathcal{L}(u^\mu)\}_{\mu \in\,(0,1]}$
is tight in $C([0,T];H)$.
\end{Proposition}

The proposition above is obtained by considering separately the case $b$ is Lipschitz and $b$ is locally Lipschitz.

\begin{Remark}
{\em In the proof of Proposition \ref{was30} it is shown that if $b$ is Lipschitz continuous, then
\[\sup_{\mu \in\,(0,1]}\,\E\sup_{s
\in\,[0,T]}|u^\mu(s)|^p_{H^\delta}\leq
c_p(T)\le(\,1+|v_0|_{H^{\delta-1}}^p+|u_0|_{H^{\delta}}^p\r).
\]
This implies that the family $\{u^\mu\}_{\mu \in\,(0,1]}$ is bounded
in $L^p(\Omega;C([0,T];H^\d))$, for any $p\geq 1$ and $\d<1/2$. In
the case that $b$ is not Lipschitz continuous we cannot prove
that. Nevertheless,  from the proof of the proposition above we
have that
 for any $\theta<1/4$ the family $\{u^\mu\}_{\mu \in\,(0,1]}$ is uniformly integrable in
 $C([0,T];C^\theta([0,L]))$, that is
\begin{equation}
\label{was22} \lim_{R\to \infty}\,\sup_{\mu
\in\,(0,1]}\,\mathbb{P}\,\le(\,\sup_{t
\in\,[0,T]}\,|u^\mu(t)|_{C^\theta([0,L])}>R\,\r)=0.
\end{equation}
Actually, if we denote by $f_T$  the inverse of the function $x
\mapsto c(1+x^{2\la})\exp(c x^2 T)$, by following the proof of the
second step in the proposition above we have
\[\begin{array}{l}
\ds{\mathbb{P}\,\le(\,\sup_{t
\in\,[0,T]}\,|u^\mu(t)|_{C^\theta([0,L])}>R\,\r)\leq
\mathbb{P}\,\le(\sup_{t
\in\,[0,T]}\,|\Pi_1\Gamma_\mu(z^\mu)(t)|_{C([0,L])}>f_T(R/2)\,\r)}\\
\vs \ds{+\mathbb{P}\,\le(\,\sup_{t
\in\,[0,T]}\,|\Pi_1\Gamma_\mu(z^\mu)(t)|_{C^\theta([0,L])}>R/2\,\r)}\\
\vs \ds{\leq \le(\frac 1 {f_T(R/2)}+\frac 2 R\r)\,\E\, \sup_{t
\in\,[0,T]}\,|\Pi_1\Gamma_\mu(z^\mu)(t)|_{C^\theta([0,L])}.}
\end{array}\]
Hence, as $f_T(R)$ diverges to $+\infty$ as $R\to \infty$, due to
\eqref{holder0} we obtain \eqref{was22}. }
\end{Remark}

Now, as in the case of additive noise, it is possible to prove the following integration by parts formula. For any
 $\varphi \in\,C^{2}([0,T]\times [0,L])$, such
that $\varphi(t,0)=\varphi(t,L)= 0$,
\begin{equation}
\label{fr64bis}
\begin{array}{l}
\ds{\int_{0}^L u^\mu(t,x)\varphi(t,x)\,dx=\int_{0}^L
u_0(x)\varphi(0,x)\,dx+\int_0^t\! \int_{0}^L
u^\mu(s,x)\le[\frac{\partial \varphi}{\partial
t}(s,x)+\Delta \varphi(s,x)\r] ds\,dx}\\
\vs \ds{+\int_0^t \int_{0}^L
b(x,u^\mu(s,x))\varphi(s,x)\,ds\,dx+\int_0^t
\int_{0}^L\varphi(s,x)g(x,u^\mu(s,x))\,w(ds,dx)+R_\mu(t),}
\end{array}
\end{equation}
where
\[ \begin{array}{l} \ds{R_\mu(t):=\mu\le(1-e^{-\frac
t\mu}\r)\int_{0}^Lv_0(x)\varphi(0,x)\,dx-\int_0^t
 e^{-\frac{t-s}\mu}\,M_\mu(s)\,ds}\\
 \vs
 \ds{-\int_0^t
 e^{-\frac{t-s}\mu}\le[\int_{0}^L\le(u_0(x)\frac{\partial \varphi}{\partial t}
 (0,x)-u^\mu(s,x)\frac{\partial \varphi}{\partial t}(s,x)+
 \int_0^s u^\mu(r,x)\frac{\partial
 ^2\varphi}{\partial
 t^2}(r,x)\,dr\r)\,dx\r]\,ds}\\
 \vs
 \ds{
 -\int_0^t \int_{0}^Le^{-\frac{t-s}\mu}\varphi(s,x)g(x,u^\mu(s,x))w(ds,dx),}
 \end{array}
 \]
 and
 \[
M_\mu(t):=
\int_0^L\le(u^\mu(t,x)\Delta\varphi(t,x)+b(x,u^\mu(t,x))\varphi(t,x)\r)\,dx.
\]

\bigskip

Then, also in this case, we have to show that the remainder term
$R_\mu(t)$ converges to zero, as the parameter $\mu$ goes to zero.
In the proof of Lemma \ref{lemma4}, $b$  was Lipschitz continuous and the noise was additive, so that
 we could prove the mean-square
convergence of  $R_\mu(t)$  to zero, for any fixed $t \geq 0$.
Here, in the case of non-Lipschitz $b$, it is  only possible to prove
convergence in probability, but, as we will show later on, this is
sufficient  in order to establish the validity of the
Smoluchowski-Kramers approximation.

\begin{Lemma}\cite[Lemma 4.7]{smolu1}
\label{rmu}
Assume Hypothesis \ref{H1}. Then, under either Hypothesis \ref{H3} or Hypothesis
\ref{H4}, for any $\e>0$
\[\lim_{\mu \to 0} \mathbb{P} \,\le(\,|R_\mu(t)|>\e\r)=0.\]
\end{Lemma}

\begin{Remark}
{\em In fact, under Hypothesis \ref{H3} we have for any fixed
$t>0$ mean-square convergence of $R_\mu(t)$ to zero, as $\mu \to
0$.}
\end{Remark}

As in Theorem \ref{teo4.6}, the tightness of the family of measures
$\{\mathcal{L}(u^\mu)\}_{\mu \in\,(0,1]}$ in $C([0,T];H)$,
together with  the integration by parts formula \eqref{fr64bis}
and the limit of Lemma \ref{rmu}, imply the following result.

\begin{Theorem}\cite[Theorem 4.9]{smolu1}
For any $\mu>0$, let $u^\mu=\Pi_1 z^\mu$, where $z^\mu$ is the
mild solution of equation \eqref{astratta}. Then, under either
Hypothesis \ref{H3} or Hypothesis \ref{H4}, for any $z_0=(u_0,v_0)
\in\,\mathcal{H}_1$, $T>0$ and $\e>0$ we have
\[\lim_{\mu\to
0}\mathbb{P}\,\le(\,|u^\mu-u|_{C([0,T];H)}>\e\r)=0,\] where $u$ is
the solution of the semi-linear stochastic heat equation
\eqref{ieq5}.

\end{Theorem}

\section{The Smoluchowski-Kramers approximation in presence of a magnetic field}
\label{sec5}

We consider here the following two dimensional system of stochastic PDEs

\begin{equation} \label{intro-eq}
\left\{
  \begin{array}{l}
    \ds{\mu \frac{\partial^2 u_{\mu}}{\partial t^2}(t,x) = \Delta u_{\mu}(t,x) + B(u_\mu(\cdot,t),t)  +\vec{m}\times \frac{\partial u_{\mu}}{\partial t}(t,x) + G(u_\mu(\cdot,t),t)\,\frac{\partial  w^{\,Q}}{\partial t}(t,x),}\\
    \vs
 \ds{u_\mu(0,x) = u_0(x), \ \ \ \frac{\partial u_\mu}{\partial t}(0,x) = v_0(x), \ \ \ x \in \mathcal{O},\ \ \ \ \ \ \ u_\mu(t,x) = 0, \ \ x \in \partial \mathcal{O}, }\\
  \end{array} \right.
\end{equation}
where $\mathcal{O}$ is a bounded regular domain in $\reals^d$, with $d\geq 1$, $B$ and $G$ are suitable nonlinearities, $\vec{m}=(0,0,m)$ is a constant vector and $w^Q(t,x)$ is a cylindrical Wiener process, white in time and colored in space, in the case of space dimension $d>1$.

By Newton's law, the vector field $u_\mu: [0,+\infty)\times \mathcal{O}\to \reals^2$ models the displacement of a continuum of electrically charged particles with constant mass density $\mu>0$ in the region $\mathcal{O} \subset \mathbb{R}^d$, in the presence of a noisy perturbation and a constant magnetic field $\vec{m}=(0,0,m)$, which is orthogonal to the plane where the motion occurs (in what follows we shall assume just for simplicity of notations $m=1$). For example, if $d=1$ and $\mathcal{O}=[0,1]$, this could model the displacement of a charged one-dimensional string, with fixed endpoints, that can move through two other spacial dimensions, where the Laplacian $\Delta$ models the forces neighboring particles exert on each other, the uniform magnetic field points in the direction of $\mathcal{O}$, $B$ is some nonlinear forcing, and $\partial w^Q/\partial t$ is a Gaussian random forcing field, whose intensity $G$ may depend on the state $u_\mu$.

In Section \ref{sec3} and Section \ref{sec4}, we have studied   the validity of the  so-called  Smoluchowski-Kramers approximation, in the case the magnetic field is  replaced by a constant friction. Namely, it has been shown that, as $\mu$ tends to $0$, the solutions of the second order system converge to the solution of the first order system which is obtained simply by taking $\mu=0$.

\medskip

One might hope that a similar result would be true in the case treated in the present section. Namely, one would expect that for any $T>0$, $\delta>0$,
\[
\lim_{\mu\to 0}\,\Pro\left(\sup_{t \in\,[0,T]}|u_\mu(t)-u(t)|_{ L^2(\mathcal{O};\reals^2)}>\delta\right)=0,
\]
where $u(t)$ is the solution of the following system of stochastic PDEs
\begin{equation}
\label{intro-eq2}
\left\{
  \begin{array}{l}
    \ds{\frac{\partial u}{\partial t}(t,x) = J_{0}^{-1}\Big[\Delta u(t,x) +  B(u(\cdot,t),t) + G(u(\cdot,t),t) \frac{\partial w^{Q}}{\partial t}(t,x)\Big] }\\
    \vs
   \ds{u(t,x) = 0, \ \ x \in \partial \mathcal{O},\ \ \ \ \ \ u(0,x) = u_0(x), \ \ \  x \in \mathcal{O}, }
  \end{array} \right.
\end{equation}
where
\[J_0^{-1}=-\,J_0=\begin{pmatrix} 0 & -1 \\ 1 & 0 \end{pmatrix}.\]

Unfortunately, as shown in \cite{smolu3} such a limit is not valid, even for finite dimensional analogues of this problem.
Actually, one can prove that if the stochastic term in \eqref{intro-eq} is replaced by a continuous function, then $u_\mu$ would converge uniformly in $[0,T]$ to the solution of \eqref{intro-eq2}. But if we have the white noise term, this is not true anymore. An explanation of this lies in the fact that, while for any continuous function $\varphi(s)$ it holds
\[\lim_{\mu\to 0}\int_0^t\sin (s/\mu)\,\varphi(s)\,ds=0,\]
if we consider a stochastic integral and replace $\varphi(s)ds$ with $dB(s)$, we have
\[\lim_{\mu\to 0}\int_0^t\sin (s/\mu)\,dB(s)\neq 0,\]
since
\[\text{Var}\le(\int_0^t\sin (s/\mu)\,dB(s)\r)=\int_0^t\sin^2(s/\mu)\,ds\to \frac t2,\ \ \ \ \text{as}\ \ \mu\downarrow 0.\]

\medskip

Nevertheless, the problem under consideration can be regularized in such a way that a counterpart of the Smoluchowski-Kramers approximation is still valid.
To this purpose, there are various ways to regularize the problem. One possible way consists in regularizing the noise (to this purpose, see \cite{smolu3}  and \cite{jjl} for  the analysis of finite dimensional systems, both in the case of constant and in the case of state dependent magnetic field). Another possible way, which is the one we are using here, consists  in introducing a small friction proportional to the velocity in equation \eqref{intro-eq}  and considering the regularized problem
\begin{equation} \label{regular-intro}
\left\{
  \begin{array}{l}
    \ds{\mu \frac{\partial^2 u^\e_{\mu}}{\partial t^2}(t) = \Delta u_{\mu}^\e(t) + B(u^\e_\mu(\cdot,t),t)  +\vec{m}\times \frac{\partial u^\e_{\mu}}{\partial t}(t) -\e\,\frac{\partial u^\e_{\mu}}{\partial t}(t)+ G(u^\e_\mu(\cdot,t),t)\frac{\partial w^{Q}}{\partial t}(t),}\\
    \vs
     \ds{u^\e_\mu(0) = u_0, \ \ \ \frac{\partial u^\e_\mu}{\partial t}(0) = v_0, \ \ \ \ \ \ \ u^\e_\mu(t,x) = 0, \ \ x \in \partial \mathcal{O}, }\\
  \end{array} \right.
\end{equation}
which now depends on two small positive parameters $\e$ and $\mu$.
  Our purpose here is showing that, for any fixed $\e>0$, we can take the limit as $\mu$ goes to $0$. Namely, we want to prove that for any $T>0$ and $p\geq1$
\[
\lim_{\mu\to 0}\E\sup_{t \in\,[0,T]}|u^\e_\mu(t)-u_\e(t)|_{L^2(\mathcal{O};\reals^2)}^p=0,\]
where $u_\e(t)$ is the unique mild solution of the problem
\[
\left\{
  \begin{array}{l}
    \ds{\frac{\partial u_{\e}}{\partial t}(t,x) = \le(J_0+\e\,I\r)^{-1}\Big[ \Delta u_{\e}(t,x) + B(u_\e(\cdot,t),t) + G(u_\e(\cdot,t),t) \frac{\partial w^{Q}}{\partial t}(t,x) \Big],}\\
    \vs
    \ds{u_\e(t,x) = 0, \ \ x \in \partial \mathcal{O},\ \ \ \ \ \ u_\e(0,x) = u_0(x),\ \ \ \ \   x \in \mathcal{O}, }
  \end{array} \right.
\]
which is precisely what we get from \eqref{regular-intro} when we formally set $\mu=0$.

\subsection{Assumptions and notations}
\label{subsec5.1}

In the present section, unlike in the rest of the paper, we shall denote by $H$ the Hilbert space $L^2(\mathcal{O},\reals^2)$, endowed with the scalar product
\[\le<(h_1,k_1),(h_2,k_2)\r>_H=\int_\mathcal{O} h_1(x)h_2(x)\,dx+\int_\mathcal{O} k_1(x)k_2(x)\,dx,\]
and the corresponding norm $|\cdot|_H$.

Now, if $\Delta$ denotes the realization of the Laplace operator in $L^2(\mathcal{O})$, endowed with Dirichlet boundary conditions,   there exists an orthonormal basis $\{\hat{e}_k\}$ for $L^2(\mathcal{O})$ and a positive sequence $\{\hat{\alpha}_k\}$ such that $\Delta \hat{e}_k = -\hat{\alpha}_k \hat{e}_k$, with $0<\hat{\alpha}_1\leq  \hat{\alpha}_k \leq \hat{\alpha}_{k+1}$. Thus, if we define for any $k \in \nat$,
 \[e_{2k-1} = (\hat{e}_k ,0), \ \alpha_{2k} = \hat{\alpha_k},\]
 \[e_{2k} = (0, \hat{e}_k), \ \alpha_{2k+1} = \hat{\alpha_k},\]
 we have that $\{e_k\}_{k=1}^\infty$ is a complete orthonormal basis of $H$. Moreover, if we define
 \[D(A)=D(\Delta)\times D(\Delta),\ \ \ A(h,k)=(\Delta h,\Delta k),\ \ \ (h,k) \in\,D(A),\]
 we have that
 \[A e_k=-\alpha_k e_k,\ \ \ \ k \in\,\nat.\]

 Next, in the same way as in section \ref{sec3} and \ref{sec4},
for any  $\delta \in \mathbb{R}$, we define $H^\delta$ to be the completion of $C^{\infty}_0(\mathcal{O};\mathbb{R}^2)$ with respect to  the norm
\[|u|_{H^\delta}^2 = \sum_{k=1}^\infty \alpha_k^\delta \left<u, e_k \right>_H^2.\]
Moreover, we define
$\H_\delta := H^\delta \times H^{\delta -1}$, and in the case $\delta=0$ we  simply set   $\H := \H_0$. Finally, for any $(h,k) \in\,\H_\delta$, we  denote
\[\Pi_1(h,k)=h,\ \ \ \ \Pi_2(h,k)=k.\]

The cylindrical Wiener process $w^Q(t,x)$ is defined as the formal sum
\[w^Q(t,x) = \sum_{k=1}^\infty Qe_k(x) \beta_k(t),\]
where $Q=(Q_1,Q_2) \in\,{\cal L}(H)$,
 $\{\beta_k\}_{k \in\,\mathbb{N}}$ is a sequence of identical, independently distributed one-dimensional, Brownian motions defined on some probability space $(\Omega,\cal{F},\mathbb{P})$ and $\{e_k\}_{k \in\,\mathbb{N}}$ is the orthonormal basis of $\H$ introduced above.

Concerning the non-linearity $B$ we assume the following conditions
\begin{Hypothesis} \label{H11}
The mapping  $B:H\times [0,+\infty)\to H$ is measurable. Moreover, for any $T>0$ there exists $\kappa_B(T)>0$ such that
  \[\left|B(u_1,t) - B(u_2,t) \right|_H \leq \kappa_B(T)|u_1-u_2|_H,\ \ \ \ u_1, u_2 \in\,H,\ \ \ t \in\,[0,T],\]
  and
  \[\sup_{t \in\,[0,T]}|B(0,t)|_H\leq \kappa_B(T).\]
\end{Hypothesis}
In the case there exists some measurable $b:\reals\times \mathcal{O}\times [0,+\infty)\to \reals$ such that for any $h \in\,H$ and $t\geq 0$
\[B(h,t)(x)=b(h(x),x,t),\ \ \ x \in\,\mathcal{O},\]
then Hypothesis \ref{H11} is satisfied if $b(\cdot,x,t):\reals\to\reals$ is Lipschitz continuous and has linear growth, uniformly with respect to $x \in\,\mathcal{O}$ and $t \in\,[0,T]$, for any $T>0$.

Concerning the diffusion coefficient $G$, we assume the following
\begin{Hypothesis}
\label{H12}
The mapping $G:H\times [0,+\infty)\to{\cal L}(L^\infty(\mathcal{O});H)$ is measurable and for any $T>0$ there exists $\kappa_G(T)>0$ such that
 \[\left|\le[G(h_1,t) - G(h_2,t)\r]z \right|_H \leq \kappa_G(T)|h_1-h_2|_H|z|_\infty,\ \ \ \ h_1, h_2 \in\,H,\ \ z \in\,L^\infty(\mathcal{O}),\ \ \ t \in\,[0,T],\]
  and
  \[\sup_{t \in\,[0,T]}|G(0,t)z|_H\leq \kappa_G(T)|z|_\infty,\ \ \ \ z \in\,L^\infty(\mathcal{O}),\ \ \  t \in\,[0,T].\]
\end{Hypothesis}
In particular, this implies that for any $h_1, h_2, ,z \in\,H$
\[
|[G^\star(h_1,t)-G^\star(h_2,t)]z|_{(L^\infty(\mathcal{O}))^\prime}\leq \kappa_G(T)|h_1-h_2|_H\,|z|_H, \ \ \ \ t \in\,[0,T].
\]
If for any $h \in\,L^2(\mathcal{O})$ and $z \in\,L^\infty(\mathcal{O})$ we define
\[[G(h,t)z](x)=g(h(x),x,t)z(x),\ \ \ \ \ x \in\,\mathcal{O},\]
for some measurable
$g:\reals^2 \times \mathcal{O} \times [0,+\infty] \to {\cal L}(\reals^2)$, then Hypothesis \ref{H12} is satisfied if
  \[\sup_{x \in \mathcal{O}}\sup_{t \in [0,T]} | g(h_1,x,t) - g(h_2,x,t)|_{{\cal L}(\reals^2)} \leq \kappa_T |h_1-h_2|_{\reals^2}\]
  and it has linear growth
  \[\sup_{x \in \mathcal{O}} \sup_{t \in [0,T]} |g(h,x,t)|_{{\cal L}(\reals^2)} \leq \kappa_T(1 + |h|_{\reals^2}).\]
Actually, in this case
\[\begin{array}{l}
\ds{|(G(h_1,t) - G(h_2,t))y |_H^2 = \int_\mathcal{O} |(g(h_1(x),x,t) - g(h_2(x),x,t))y(x)|_{\reals^2}^2 dx}\\
\vs
\ds{\leq \kappa_T \int_\mathcal{O} |h_1(x)-h_2(x)|_{\reals^2}^2 |y(x)|_{\reals^2}^2 dx \leq |h_1 -h_2|_H^2 |y|_\infty^2,}
\end{array}
\]
and by the same reasoning
\[
  |G(h,t)y|_{H} \leq \kappa_T(1+ |h|_H) |y|_\infty.
\]

\medskip

Now, for any $\mu>0$ and $\d \in\,\reals$, as in Section \ref{sec2}we define on $\H_\d$ the unbounded linear operator
\[A_\mu (u,v)=\frac 1\mu(\mu v,Au-J_0v),\ \ \ \ (u,v) \in\,D(A_\mu)=\H_{\d+1},\]
where
$J_0$ is the skew symmetric $2 \times 2$ matrix
\[J_0 = \begin{pmatrix} 0 & 1 \\ -1 & 0 \end{pmatrix}.\]
It can be proven that $A_\mu$ is the generator of a strongly continuous group of bounded linear operators $\{S_\mu(t)\}_{t\geq 0}$ on each $\H_\d$ (for a proof see \cite[Section 7.4]{pazy}).

Moreover, for any $\mu>0$ we define
\[B_\mu:\H\times [0,+\infty)\to\H,\ \ \ \ \ (z,t) \in\,\H\times [0,+\infty)\mapsto \frac 1\mu(0,B(\Pi_1 z,t)),\]
and
\[G_\mu:\H\times [0,+\infty)\to{\cal L}(L^\infty(\mathcal{O}),\H),\ \ \ \ \ (z,t) \in\,\H\times [0,+\infty)\mapsto \frac 1\mu (0,G(\Pi_1 z,t)).\]
With these notations, if we set
\[z_\mu(t)=\le(u_\mu(t),\frac{\partial u_\mu}{\partial t}(t)\r),\]
system \eqref{intro-eq} can be rewritten as the following stochastic equation in the Hilbert space $\H$
\[
dz_\mu(t)=\le[A_\mu z_\mu(t)+B_\mu(z_\mu(t),t)\r]\,dt+G_\mu(z_\mu(t),t)dw^Q(t),\ \ \ \ \ z_\mu(0)=(u_0,v_0).\]

\subsection{The approximating semigroup}
\label{subsec5.2}

For any $\mu,\e>0$ and $\d \in\,\reals$, we define
\[A_{\mu}^{\e}(u,v)=\frac 1\mu(\mu\,v,Au-J_\e v),\ \ \ \ \ (u,v) \in\,D(A_{\mu}^{\e})=\H_{\d+1},\]
where
\[J_\e = J_0 +\e I =  \begin{pmatrix} \e & 1 \\ -1 & \e \end{pmatrix},\ \ \ \ \e>0. \]
As we have seen for $A_\mu$, it is possible to prove that for any $\mu, \e>0$ the operator $A_{\mu}^{\e}$ generates a strongly continuous group of bounded linear operators $S_{\mu}^\e(t)$, $t\geq 0$,  on $\H_\d$.

\begin{Lemma}\cite[Lemma 3.1]{csal}
\label{energy-est-lem}
  For any $(u,v) \in\,\H_\theta$, with $\theta \in\,\reals$, and for any $\mu, \e>0$ let us define
  \[u_{\mu}^\e(t): = \Pi_1 S_\mu^\e(t)(u,v),\ \ \ \ \ \ v_\mu^\e(t):=\Pi_2 S_\mu^\e(t)(u,v).\]
Then
\begin{equation} \label{energy-est-1}
  \mu \left|v_\mu^\e(t) \right|_{H^{\theta-1}}^2 + |u_\mu^\e(t)|_{H^\theta}^2+2\e\int_0^t|v_\mu^\e(s)|_{H^{\theta-1}}^2\,ds
 = \mu|v|_{H^{\theta-1}}^2 + |u|_{H^\theta},
\end{equation}
and
\[
  \mu |u_\mu^\e(t)|_{H^\theta}^2 + \left| \mu v_\mu^\e(t) + J_\e u_\mu^\e(t) \right|_{H^{\theta-1}}^2+2\e\int_0^t|u_\mu^\e(s)|_{H^\theta}^2\,ds
=\mu |u|_{H^\theta}^2 + |\mu v + J_\e u|_{H^{\theta -1}}^2.
\]
\end{Lemma}

Notice that in particular this implies that for any $\mu,\e>0$ there exists $c_{\mu,\e}>0$ such that for any $(u,v) \in\,\H_\theta$
\[\int_0^\infty |S^\e_\mu(t)(u,v)|_{\H_\theta}^2\,dt\leq \frac {c_{\mu,\e}}{2\e}|(u,v)|_{\H_\theta}^2.\]
As a consequence of the Datko theorem, this allows to conclude that there exist $M_{\mu,\e},$ and  $\omega_{\mu,\e}>0$ such that
\[\|S_\mu^\e(t)\|_{{\cal L}(\H_\theta)}\leq M_{\mu,\e}\,e^{-\omega_{\mu,\e}t},\ \ \ \ t\geq 0.\]

\begin{Lemma}\cite[Lemma 3.2]{csal}
 \label{Smu-bound-lem}
  For any $\mu, \e>0$,  and for any $\theta \in \reals$ and $\gamma \in [0,1]$ it holds
\[
\left|\Pi_1 S_\mu^\e(t)(0,v) \right|_{H^\theta} \leq 2^{\gamma} \mu^{\frac{1+\gamma}{2}}  |v|_{H^{\theta+\gamma-1}},\ \ \ t\geq 0,\ \ \ \ v \in\,H^{\theta+\gamma-1}.\]
\end{Lemma}

Now, analogously as in Sections \ref{sec3} and \ref{sec4} for any $\mu>0$ we define the bounded linear operator
\[Q_\mu:H\to \H,\ \ \ \ u \in\,H\mapsto \frac1\mu(0,Qu) \in\,\H.\]

\begin{Lemma}\cite[Lemma 3.3]{csal}
\label{l1}
Assume that there exists a non-negative sequence $\{\la_k\}_{k \in\,\nat}$ such that
\[Qe_k=\la_k e_k,\ \ \ \ k \in\,\nat.\]
Then,
for any $0<\delta<1$ and $\e>0$ there exists a constant $c= c(\e, \d)>0$ such that for any $k \in \nat$ and $\theta>0$
  \[
    \sup_{\mu>0} \int_0^\infty s^{-\d} \left| \Pi_1 S_\mu^\e(s) Q_\mu e_k \right|_H^2 ds \leq c\, \frac{\la_k^2}{\alpha_k^{1-\d}},
  \]
  and
  \[
    \sup_{\mu>0} \mu^{1+ \d} \int_0^\infty s^{-\d} \left| \Pi_2 S_\mu^\e(s) Q_\mu e_k \right|_{H^{\theta-1}}^2 ds \leq c\, \frac{\la_k^2}{\alpha_k^{1-\theta}}.
  \]
\end{Lemma}

Now, for any $\e>0$ we define
\[A_\e:=J_\e^{-1} A=\frac{1}{1+\e^2}\begin{pmatrix} \e  & -1 \\ 1 & \e \end{pmatrix}\Delta,\]
and we denote by $T_\e(t)$, $t\geq 0$, the strongly continuous semigroup generated by $A_\e$ in $H^\theta$, for any $\theta \in\,\reals$.
Moreover, we denote
\[Q_\e=J_\e^{-1}Q.\]

\begin{Lemma}\cite[Lemma 3.4]{csal}
\label{l7}
We have
  \[
    \|T_\e(t) \|_{{\cal L}(H^\theta)}  \leq e^{-\frac{\e\a_1}{1+\e^2}t},\ \ \ \ \ t\geq 0.
  \]

 Moreover, if there exists a non-negative sequence $\{\la_k\}_{k \in\,\nat}$ such that
\[Qe_k=\la_k e_k,\ \ \ \ k \in\,\nat,\]
then, for any $0<\d<1$ and $\e>0$ there exists a constant $c=c(\d,\e)$ such that for any $k \in\,\nat$
  \begin{equation} \label{bound-2}
    \int_0^\infty s^{-\d} \left| T_\e(s) Q_\e e_k \right|_H^2 ds \leq c \frac{\lambda_k^2}{\alpha_k^{1-\d}}.
  \end{equation}
  Finally, for any  $k \in\,\nat$
  \[
    \int_0^T s^{-\d} \left| T_\e(s) Q_\e e_k \right|_H^2 ds \leq \frac{1}{1-\d}\,T^{1-\d}\,\lambda_k^2.
    \]

\end{Lemma}

In view of the previous estimates for $S_\mu^\e(t)$ and $T_\e(t)$,  the following convergence result holds.
\begin{Theorem}\cite[Theorem 3.5]{csal}
  For any $\e >0$, $0<t_0<T$, and   $n \in \nat$,
  \[
    \lim_{\mu \to 0} \sup_{t \leq T} \sup_{|u|_H\leq 1} \left| \Pi_1 S_\mu^\e(t)(P_n u,0) - T_\e(t) P_n u \right|_{H} =0,
  \]
  and
  \[
    \lim_{\mu \to 0} \sup_{0<t_0 \leq t\leq T} \sup_{|v|_H \leq 1} \left|\frac{1}{\mu} \Pi_1 S_\mu^\e(t)(0,P_n v) - T_\e(t)J_\e^{-1} P_n v \right|_H = 0,
  \]
  where $P_n$ is the projection of $H$ onto the $n$-dimensional subspace $H_n:=\text{\em{span}}\{e_1,\ldots,e_{2n}\}$.
\end{Theorem}

The two limits above imply that for any $\e>0$ and $T>0$ and for any $(u,v) \in \H$,
  \begin{equation}
  \label{Smu-to-Se-at-a-point}
    \lim_{\mu \to 0} \sup_{t \leq T} |\Pi_1 S_\mu^\e(t)(u,v) - T_\e(t) u|_H=0,
  \end{equation}
and, for any $v \in H$ and $0< t_0 \leq T$,
  \begin{equation} \label{Smu-to-Se-3-at-a-point}
    \lim_{\mu \to 0} \sup_{t_0 \leq t \leq T} \left|\frac 1\mu\Pi_1 S_\mu^\e(t) \left(0, v \right) - T_\e(t) J_\e^{-1} v \right|_H =0.
  \end{equation}
Moreover,  for any $\e>0$, $T>0$ and $p\geq 1$ and for any $\psi \in L^p(\Omega; L^p([0,T];H))$,
  \begin{equation} \label{Smu-to-Se-2}
    \lim_{\mu \to 0} \,\E \sup_{t \in\,[0,T]} \left|\frac{1}{\mu} \int_0^t \Pi_1 S_\mu^\e(t-s) (0,\psi(s)) ds - \int_0^t T_\e(t-s) J_\e^{-1} \psi(s) ds \right|_H^p = 0.
  \end{equation}

\subsection{Approximation by small friction for additive noise }
\label{subsec5.3}
We assume here that the noisy perturbation in system \eqref{intro-eq} is of additive type, that is $G(u,t) = I$, for any $u \in\,H$ and $t\geq 0$. Moreover, we assume that the covariance operator $Q$ satisfies the following condition.
\begin{Hypothesis}
\label{H13}
There exists a non-negative sequence $\{\la_k\}_{k \in\,\nat}$ such that $Q e_k=\la_k e_k$, for any $k \in\,\nat$. Moreover, there exists $\d>0$ such that
\[\sum_{k=1}^\infty \frac{\la_k^2}{\a_k^{1-\d}}<\infty.\]
\end{Hypothesis}

With the notations we have introduced above and \ref{sec3}, if we denote
\[z_\mu^\e(t)=(u_\mu^\e(t),\frac{\partial u_\mu^\e}{\partial t}(t)),\ \ \ t \geq 0,\]
the regularized system \eqref{regular-intro} can be rewritten as the abstract evolution equation
\begin{equation}
\label{reg-abst}
dz_{\mu}^\e(t)=\le[A^\e_\mu z^\e_\mu(t)+B_\mu(z_\mu^\e(t),t)\r]\,dt+Q_\mu dw(t),\ \ \ \ \ z_\mu^\e(0)=(u,v)
\end{equation}
in the Hilbert space $\H$.

Our purpose here is to show that for any fixed $\e>0$ the process $u_\mu^\e(t)$ converges to the solution $u_\e(t)$ of the following system of stochastic PDEs
\begin{equation} \label{first-order-eq-w-frict}
\left\{
  \begin{array}{l}
    \ds{\frac{\partial u_\e}{\partial t}(t) = J_\e^{-1} \Delta u_\e(t) + B_\e(u_\e(t),t) + \frac{\partial w^{Q_\e}}{\partial t}(t)}\\
    \vs
    \ds{u_\e(0) = u_0,\ \ \ \ u_\e(t,x)=0,\ \ \ x \in\,\partial \mathcal{O},}
  \end{array}
  \right.
\end{equation}
where for any $\e>0$ we have defined $Q_\e=J_\e^{-1} Q$ and
\[B_\e(u,t)=J_\e^{-1} B(u,t),\ \ \ u \in\,H,\ \ \ t\geq 0.\]
Notice that with these notations, system \eqref{first-order-eq-w-frict} can be rewritten as the abstract evolution equation
\begin{equation}
\label{first-abst}
du_\e(t)=\le[A_\e u_\e(t)+B_\e(u_\e(t),t)\r]\,dt+Q_\e dw(t),\ \ \ \ u_\e(0)=u_0,\end{equation}
in the Hilbert space $H$.

\medskip

According to Lemma \ref{l1}, due to Hypothesis \ref{H13} for any $t\geq 0$ we have
\[\int_0^t s^{-\d}\sum_{k=1}^\infty |S^\e_\mu(t-s)Q_\mu e_k|^2_{\H}\,ds\leq c\le(1+\mu^{-(1+\d)}\r)\sum_{k=1}^\infty\frac{\la_k^2}{\a_k^{1-\d}}.\]
This implies that the stochastic convolution
\[\Gamma^\e_\mu(t):=\int_0^t S^\e_\mu(s) Q_\mu\,dw(s),\ \ \ \ t\geq 0,\] takes values in $L^p(\Omega;C([0,T];\H))$, for any $T>0$ and $p\geq 1$ (for a proof see \cite{dpz1}).
Therefore, as the mapping $B_\mu(\cdot,t):\H\to\H$ is Lipschitz-continuous, uniformly with respect to $t \in\,[0,T]$, we have that there exists a unique process $z^\e_\mu \in\,L^p(\Omega;C([0,T];\H))$  which solves equation \eqref{reg-abst} in the mild sense, that is
\[z^\e_\mu(t)=S^\e_\mu(t)(u_0,v_0)+\int_0^t S^\e_\mu(t-s)B_\mu(z^\e_\mu(s),s)\,ds+\Gamma^\e_\mu(t).\]

In the same way, due to \eqref{bound-2} we have that the stochastic convolution
\[\Gamma_\e(t):=\int_0^t T_\e(s) Q_\e\,dw(s),\ \ \ \ t\geq 0,\] takes values in $L^p(\Omega;C([0,T];H))$, for any $T>0$ and $p\geq 1$, so that, as the mapping $B_\e(\cdot,t):H\to H$ is Lipschitz-continuous, uniformly with respect to $t \in\,[0,T]$, we can conclude that there exists a unique process $u_\e \in\,L^p(\Omega;C([0,T];H))$ solving equation \eqref{first-abst} in mild sense, that is
\[
  u_\e(t) = T_\e(t)u_0 + \int_0^t T_\e(t-s) B_\e(u_\e(s),s) ds + \Gamma_\e(t).
\]

\begin{Theorem}\cite[Theorem 4.1]{csal}
 \label{conv-in-mu-thm}
Under Hypotheses \ref{H11} and \ref{H13}, for any $\e>0$, $T>0$ and $p\geq 1$ and for any initial conditions $z_0=(u_0,v_0) \in\,\H$, we have
\[    \lim_{\mu \to 0} \E \sup_{t\leq T} \left|u_\mu^\e(t) - u_\e(t) \right|_H^p  = 0.
  \]
\end{Theorem}

\subsection{Approximation by small friction for multiplicative noise } \label{subsec5.4}

In this section we assume that the space dimension $d=1$ and $\mathcal{O}$ is a bounded interval, the diffusion coefficient $G$ satisfies Hypothesis \ref{H12} and the covariance operator $Q$ satisfies the following condition.
\begin{Hypothesis}
\label{H14}
There exists a bounded non-negative sequence $\{\la_k\}_{k \in\,\nat}$ such that
\[Q e_k=\la_k e_k,\ \ \ \ k \in\,\nat.\]
\end{Hypothesis}

We begin by studying the stochastic convolutions
\[
\Gamma_\mu^\e(z)(t) := \int_0^t S_\mu^\e(t-s) G_\mu(z(s),s)dw^Q(s),\ \ \ \ z \in\,L^p(\Omega,C([0,T];\H)),\]
and
\[\Gamma_\e(u)(t) = \int_0^t T_\e(t-s) G_\e(u(s),s) dw^Q(s), u \in\,L^p(\Omega,C([0,T];H)).\]
With the notations introduced above, the regularized system \eqref{regular-intro} can be rewritten as
\begin{equation}
\label{s9}
dz^\e_\mu(t)=\le[A^\e_\mu z^\e_\mu(t)+B_\mu(z^\e_\mu(t),t)\r]\,dt+G_\mu(z^\e_\mu(t),t)\,dw^Q(t),\ \ \ \ z^\e_\mu(0)=(u_0,v_0),
\end{equation}
and the limiting problem \eqref{first-order-eq-w-frict} can be rewritten as
\begin{equation}
\label{s10}
du_\e(t)=\le[A_\e u_\e(t)+B_\e (u_\e(t),t)\r]\,dt+G_\e(u_\e(t),t)\,dw^Q(t),\ \ \ \ \ u_\e(0)=u_0,\end{equation}
where
\[G_\e(u,t)=J_\e^{-1}G(u,t).\]

\begin{Lemma}\cite[Lemma 5.1]{csal}
 \label{lem:stoch-conv-Lip}
Under Hypotheses \ref{H12} and \ref{H14}, for any $\mu, \e>0$, $T \geq 0$ and $p>4$ we have
  \[z \in\,L^p(\Omega;C([0,T];\H))\Longrightarrow \Gamma_\mu^\e(z) \in L^p(\Omega;C([0,T];\H)). \]
Moreover, there exists a constant $c:=c(\e,\mu,p,T)$ such that
  \[
  \E | \Gamma_\mu^\e(z_1) - \Gamma_\mu^\e(z_2) |^p_{C([0,T];\H)} \leq c_p\, \int_0^T\E |\Pi_1 z_1-\Pi_1 z_2|^p_{C([0,\si];H)}\,d\si.\]
\end{Lemma}

\begin{Remark}
{\em We can show that the first component $\Pi_1\Gamma^\e_\mu$ satisfies
  \begin{equation}
  \label{s30}
  \sup_{\mu>0} \E | \Pi_1\Gamma_\mu^\e(z_1) - \Pi_1\Gamma_\mu^\e(z_2) |^p_{C([0,T];\H)} \leq c_p\, \int_0^T\E |\Pi_1 z_1-\Pi_1 z_2|^p_{C([0,\si];H)}\,d\si
,\end{equation}
for a constant $c=c(\e,p,T)>0$ that is independent of $\mu$.}
\end{Remark}

Lemma \ref{lem:stoch-conv-Lip} states that the mapping
\[z \in\,L^p(\Omega;C([0,T];\H))\mapsto \Gamma_\mu^\e(z) \in\,L^p(\Omega;C([0,T];\H)),\]
is Lipschitz continuous. Therefore, as the mapping $B_\mu(\cdot,t):\H\to \H$, is Lipschitz continuous, uniformly for $t \in\,[0,T]$, we have that for any initial condition $z_0=(u_0,v_0) \in\,\H$, system \eqref{s9} admits a unique adapted mild solution $z_\mu^\e \in\,L^p(\Omega;C([0,T];\H))$.

In the same way, we have the following result for the stochastic convolution associated with the parabolic problem.

\begin{Lemma}\cite[Lemma 5.3]{csal}
\label{lem:stoch-conv-epsilon-Lip}
 Under Hypotheses \ref{H12} and \ref{H14}, for any $\e, T \geq 0$ and any $p>4$
   \[u \in\,L^p(\Omega;C([0,T];H)\Longrightarrow \Gamma_\e(u) \in L^p(\Omega;C([0,T];H)). \]
   Moreover,
there exists a constant $c:=c(\e,p,T)$ such that for any $u,v \in\,L^p(\Omega;C([0,T];H))$
 \begin{equation}
 \label{s46}
 \E \left| \Gamma_\e(u) - \Gamma_\e(v) \right|_{C([0,T];H)}^{p} \leq c\, \int_0^T\E|u-v|_{C([0,\si];H)}^{p}\,d\si.\end{equation}
 If we assume that
 \[\sum_{k=1}^\infty \la_k^2<\infty,\]
 then the constant $c$ in \eqref{s46} is independent of $\e>0$.
\end{Lemma}

As a consequence of this lemma, since the mapping $B_\e(\cdot,t):H\to H$ is Lipschitz continuous, uniformly for $t \in\,[0,T]$, we have that for any initial condition $u_0 \in\,\H$, system \eqref{s9} admits a unique adapted mild solution $u_\e \in\,L^p(\Omega;C([0,T];H))$.

\begin{Theorem}\cite[Theorem 5.4]{csal}
 \label{thm:Gamma-mu-to-Gamma-e}
  For any  fixed $\e>0$, $T>0$ and $p\geq 1$,
    \[\lim_{\mu \to 0} \E \left|\Pi_1\Gamma_\mu^\e((u,0)) - \Gamma_\e(u) \right|_{C([0,T];H)}^{p}=0.\]
\end{Theorem}

This follows from a stochastic factorization argument, and limits \eqref{Smu-to-Se-at-a-point}, \eqref{Smu-to-Se-3-at-a-point} and \eqref{Smu-to-Se-2}, combined with arguments analogous to those used in the proof of Lemma \ref{lem:stoch-conv-Lip} (see \cite[Lemma 5.1]{csal}).

\begin{Theorem}
\label{ts31}
Let $z^\e_\mu=(u^\e_\mu,v^\e_\mu)$ and $u_\e$ be the mild solutions of problems \eqref{s9} and \eqref{s10}, with initial conditions $z_0 \in\,\H$ and $u_0=\Pi_1 z_0 \in\,H$, respectively.
 Then, under Hypotheses \ref{H2}, \ref{H3} and \ref{H7}, for any $T>0$,  $\e>0$ and $p\geq 1$ we have
  \[\lim_{\mu \to 0}\E |u^\e_\mu - u_\e|_{C([0,T];H)}^{p}=0.\]
\end{Theorem}

Actually, we have
  \[u^\e_\mu(t) = \Pi_1 S_\mu^\e(t)(u_0,v_0) + \Pi_1 \int_0^t S_\mu^\e(t-s) B_\mu(z_\mu^\e(s),s) ds + \Pi_1 \Gamma_\mu^\e(z^\e_\mu)(t), \]
and
  \[u_\e(t) = T_\e(t)u_0 + \int_0^t T_\e(t-s) B_\e(u_\e(s),s) + \Gamma_\e(u_\e)(t).\]
  Then
  \[\begin{array}{l}
  \ds{\left|u_\mu^\e(t) - u_\e(t) \right|_H \leq \left| \Pi_1 S_\mu^\e(t) (u_0,v_0) - T_\e(t)u_0\right|_H}\\
  \vs
  \ds{ + \left| \int_0^t \Pi_1 S_\mu^\e(t-s)[B_\mu(z^\e_\mu(s),s) - B_\mu((u_\e(s),0),s)] ds \right|_H }\\
  \vs
  \ds{ + \left| \frac 1\mu\int_0^t \Pi_1 S_\mu^\e(t-s) (0,B(u_\e(s),s))ds - \int_0^t T_\e(t-s) J_\e^{-1} B(u_\e(s),s) ds \right|_H}\\
  \vs
   \ds{+ \left|\Pi_1\le[\Gamma_\mu^\e(z^\e_\mu)(t) - \Gamma_\mu^\e((u_\e(t),0))\r] \right|_{H} + \left|\Pi_1 \Gamma_\mu^\e(u_\e(t),0) - \Gamma_\e(u_\e)(t) \right|_H.}
  \end{array}\]
 By Lemma \ref{Smu-bound-lem}, and Hypothesis \ref{H11}, there is a constant independent of $\mu$ and of  $0<s<t$, such that
  \[\left|\Pi_1 S_\mu^\e(t-s)[B_\mu(z^\e_\mu(s),s) -B_\mu((u_\e(s),0),s)] \right|_H \leq c\, |u_\mu^\e(s) - u_\e(s)|_H, \]
  so that for any $p\geq 2$
  \[\left|\int_0^t \Pi_1 S_\mu^\e(t-s)[B_\mu(z^\e_\mu(s),s) - B_\mu((u_\e(s),0),s)] ds\right|_H^{p}
  \leq c_p\, t^{p-1} \int_0^t |u^\e_\mu-u_\e|_{C([0,s];H)}^{p} ds. \]
Thanks to \eqref{s30}, this implies
  \[\begin{array}{l}
    \ds{\E\left|u^\e_\mu - u_\e \right|_{C([0,t];H)}^{p} \leq c_p\, T^{p-1} \int_0^t \E \left| u^\e_\mu - u_\e \right|_{C([0,s];H)}^{p} ds}\\
    \vs
     \ds{+c_p\,\sup_{s\leq t} \left| \Pi_1 S_\mu^\e(s) (u_0,v_0) - T_\e(t)u_0\right|_H^{p} + c_p\E\left| \Pi_1\Gamma_\mu^\e((u_\e,0)) -
     \Gamma_\e(u_\e) \right|_{C([0,t];H)}^{p} }\\
     \vs
  \ds{ +c_p\,\E\sup_{s \leq t} \left|\frac 1\mu\int_0^s \Pi_1 S_\mu^\e(s-r) (0,B(u_\e(r),r))dr - \int_0^s T_\e(s-r) B_\e(u_\e(r),r) dr \right|_H^{p},}
    \end{array}\]
and the Gr\"onwall's inequality yields
  \[\begin{array}{l}
    \ds{\E \left| u^\e_\mu - u_\e \right|_{C([0,T];H)}^{p} }\\
    \vs
 \ds{   \leq c_p(T)\le( \sup_{s\leq T} \left| \Pi_1 S_\mu^\e(s) (u_0,v_0) - T_\e(t)u_0\right|_H^{p} + \E\left| \Pi_1\Gamma_\mu^\e((u_\e,0)) -
     \Gamma_\e(u_\e) \right|_{C([0,T];H)}^{p} \r)}\\
     \vs
  \ds{ +c_p(T)\,\E\sup_{s \leq T} \left|\frac 1\mu\int_0^s \Pi_1 S_\mu^\e(s-r) (0,B(u_\e(r),r))dr - \int_0^s T_\e(s-r) B_\e(u_\e(r),r) dr \right|_H^{p}.}
  \end{array}\]
  Finally, the result follows because of \eqref{Smu-to-Se-at-a-point}, \eqref{Smu-to-Se-2}, and Theorem \ref{thm:Gamma-mu-to-Gamma-e}.

\subsection{The convergence for $\e \downarrow 0$ } \label{subsec5.6}

In the previous sections, we have shown that under suitable conditions on the coefficients and the noise, for any fixed $\e>0$, $T>0$ and $p\geq 1$
\[\lim_{\mu \to 0} \E \left|u^\e_\mu-u_\e  \right|_{C([0,T];H)}^p =0. \]
This limit is not uniform in $\e>0$, and the limit is not true for $\e=0$. In this section we want  to show that
\begin{equation}
\label{limite}
\lim_{\e \to 0} \E \left|u_\e - u \right|_{C([0,T];H)}^p = 0,\end{equation}
where $u$ is the mild solution of the problem
\begin{equation}
\label{s33}
du(t)=\le[A_0 u(t)+B_0(u(t),t)\r]\,dt+G_0(u(t),t)\,dw^Q(t),\ \ \ \ \ u(0)=u_0,\end{equation}
with
\[A_0:=J_0^{-1} A,\ \ \ B_0=J_0^{-1} B,\ \ \ \ G_0=J_0^{-1} G.\]

This statement is not true unless we strengthen Hypothesis \ref{H14}. Actually, Hypothesis \ref{H14} is the weakest assumption on the regularity of the noise that implies Theorem \ref{conv-in-mu-thm} and Theorem \ref{ts31}, for $\e>0$. But in order to prove \eqref{limite} we need to assume the following stronger condition on the covariance $Q$.

\begin{Hypothesis} \label{H15}
There exists a non-negative sequence $\{\la_k\}_{k \in\,\nat}$ such that   $Q e_k = \lambda_k e_k$, for any $k \in\,\nat$, and
   \[\sum_{k=1}^\infty \lambda_k^2 < +\infty.\]
\end{Hypothesis}

In what follows, we shall denote by $T_0(t)$, $t\geq 0$, the semigroup generated by the differential operator $A_0$ in $H$, with $D(A_0)=D(A)$.
The semigroup $T_0(t)$ is strongly continuous in $H$. Moreover, if we define $u(t)=T_0(t)x$, for $x \in\,D(A_0)$, we have
\[\le\{\begin{array}{l}
\ds{\frac{\partial u_1}{\partial t}(t)=-\Delta u_2(t),\ \ \ \ u_1(0)=x_1}\\
\vs
\ds{\frac{\partial u_2}{\partial t}(t)=\Delta u_1(t),\ \ \ \ u_2(0)=x_2}
\end{array}\r.\]
This means that if we take the scalar product in $H^\theta$ of the first equation by $u_1$ and of the second equation by $u_2$, we get
\[\frac d{dt}|u(t)|^2_{H^\theta}=0,\]
so that
\begin{equation}
\label{s32}
|T_0(t)x|_{H^\theta}=|x|_{H^\theta},\ \ \ \ t\geq 0,\end{equation}
for any $\theta \in\,\reals$ and $x \in\,H$.

Now, let us consider the stochastic convolution associated with problem \eqref{s33}, in the simple case $G=I$
\[\Gamma(t)=\int_0^t T_0(t-s)Qdw(s),\ \ \ \ t \geq0.\]
As a consequence of \eqref{s32}, we have
\[\E\,|\Gamma(t)|_H^2=\int_0^t\sum_{k=1}^\infty |T_0(s)Qe_k|_H^2\,ds=\int_0^t\sum_{k=1}^\infty |Qe_k|_H^2\,ds=t \sum_{k=1}^\infty \la_k^2,\]
and this implies that Hypothesis \ref{H15} is necessary in order to have a solution in $H$ for the limiting equation \eqref{s33}.

As proven in \cite[Lemma 6.1 and Lemma 6.2]{csal},
  The matrix $J_\e^{-1}$ converges to $J_0^{-1}$ in ${\cal L}(\mathbb{R}^2)$. Furthermore,  for any $T\geq0$,
  \[
  \lim_{\e \to 0} \sup_{t \in\,[-T,T]} \left\|e^{tJ_\e^{-1}} - e^{t J_0^{-1}} \right\|_{{\cal L}(\mathbb{R}^2)} =0
  \]
Moreover, for any $n \in \nat$, and $T\geq0$,
   \begin{equation} \label{Se-to-S0-1}
     \lim_{\e \to 0}\,\sup_{t  \in\,[0,T]} \sup_{|x|_H \leq 1} \left|T_\e(t)P_{n} x - T_0(t)P_n x \right|_H = 0.
   \end{equation}
Notice that, as in the proof of Theorem \ref{conv-in-mu-thm}, this implies that for any $x \in\,H$
\[\lim_{\e\to 0}\,\sup_{t \in\,[0,T]}|T_\e(t)x-T_0(t)x|_H=0.
\]

Now, as a consequence of \eqref{Se-to-S0-1}, 
and by arguments similar to those used to prove \eqref{Smu-to-Se-2}
\begin{Lemma}\cite[Lemma 6.3]{csal}
\label{s34}
  For any $\psi \in L^1(\Omega; L^1([0,T];H))$
  \[
  \lim_{\e \to 0}\E \sup_{t  \in\,[0,T]} \left| \int_0^t \left( T_\e(t-s)J_\e^{-1}\psi(s) - T_0(t-s)J_0^{-1} \psi(s) \right)ds \right|_H = 0.
  \]
\end{Lemma}

Now, $u_\e$ is the unique mild solution in $L^p(\Omega;C([0,T];H)$ of problem \eqref{first-abst} (in the case of additive noise) or problem \eqref{s10} (in the case of multiplicative noise), so that
\[u_\e(t)=T_\e(t)u_0+\int_0^t T_{\e}(t-s)B_\e(u_\e(s),s)\,ds+\Gamma_\e(u_\e)(t).\]
Moreover, $u(t)$ is the unique mild solution in $L^p(\Omega;C([0,T];H))$ of the problem
\[du(t)=\le[A_0 u(t)+B_0(u(t),t)\r]\,dt+G_0(u(t),t)\,dw^Q(t),\ \ \ \ u(0)=u_0,\]
with
$G_0=J_{0}^{-1}I$ or $G_0=J_0^{-1} G$, so that
\[u(t)=T_0(t)u_0+\int_0^t T_{0}(t-s)B_0(u_\e(s),s)\,ds+\Gamma_0(u_\e)(t).\]

Then, in view the previous two lemmas, we have that the arguments used in the proof of Theorem \ref{conv-in-mu-thm} and Theorems \ref{thm:Gamma-mu-to-Gamma-e} and \ref{ts31} can be repeated and we have the following result.

\begin{Theorem}\cite[Theorem 6.4]{csal}
\label{ts56}
Assume either $G$ satisfies Hypothesis \ref{H3} or $G(x,t)=I$. Then, under Hypotheses \ref{H2} and \ref{H6}, we have that for any $T>0$ and $p\geq 1$
\[
\lim_{\e\to 0}\,\E\,|u_\e-u|^p_{C([0,T];H)}=0.
\]
\end{Theorem}

In fact, it is possible to show  that the convergence result proved above for $\e\downarrow 0$ is also valid for the second order system, that is for every $\mu>0$ fixed.

\begin{Theorem}\cite[Theorem 6.5]{csal}
\label{ts55}
Assume either $G$ satisfies Hypothesis \ref{H12} or $G(x,t)=I$. Then, under Hypotheses \ref{H11} and \ref{H15}, we have that
  for any initial conditions $(u_0,v_0)$ and $\mu>0$
  \[\lim_{\e \to 0} \E\,|z_\mu^\e - z_\e|_{C([0,T];\H)}^{2m},\]
  for any $T>0$ and $p\geq 1$.
\end{Theorem}
As long as we can show that $S^\e_\mu(t) P_n z\to S^0_\mu(t) P_n z$ for any fixed $n$, we can prove Theorem \ref{ts55} by following the arguments of Theorems  \ref{conv-in-mu-thm} and \ref{ts31}. Fortunately, we can prove something stronger. Actually, it is possible to show  that $ \sup_{t \geq0}\|S^\e_\mu(t) - S^0_\mu(t)\|_{\mathcal{L}(\H)}=0$, for any fixed $\mu>0$.

To this purpose, it is useful to introduce an equivalent norm  on $H \times H^{-1}$, depending on $\mu>0$,  by setting
\[|(u,v)|_{\H(\mu)}^2 = |x|_H^2 + \mu|y|_{H^{-1}}^2.\]
Because of \eqref{energy-est-1}, for any $\e\geq0$,
\[
  \sup_{t \geq 0}\| S_\mu^\e(t)\|_{{\cal L}(\H)} \leq 1.
\]
Note that if $\e=0$, then, by \eqref{energy-est-1}, for any $z \in \H$ and $t \geq0$,
\[|S_\mu^0(t) z|_{\H(\mu)} = |z|_{\H(\mu)}.\]
\begin{Lemma}\cite[Lemma 6.6]{csal}
  For fixed $\mu>0$ and  $T>0$,
  \[
    \lim_{\e \to 0} \sup_{t \in\,[0,T]} \|S_\mu^\e(t) - S_\mu^0(t) \|_{{\cal L}(\H(\mu))} =0.
  \]
\end{Lemma}

\section{The long time behavior}
\label{sec6}

In this section we study the relation between the stationary
distributions of the processes $u^\mu(t)$ and $u(t)$, defined
respectively as the solution of the semi-linear stochastic damped
wave equation \eqref{quasilinear} and as the solution of the
semi-linear stochastic heat equation \eqref{heat}.

If we set
\[z^\mu(t):=(u^\mu(t),v^\mu(t)),\ \ \ \ t\geq 0,\ \ \mu>0,\]
with the notations introduced in Section
\ref{sec3} we can write equation \eqref{quasilinear} as the
abstract evolution equation on the Hilbert space
$\mathcal{H}=H\times H^{-1}=L^2(\mathcal{O})\times H^{-1}(\mathcal{O})$
\begin{equation}
\label{cage} dz^\mu(t)= \le[ A_\mu z^\mu(t)+B_\mu
(z^\mu(t))\r]\,dt+dw^{\mathcal{Q}_\mu}(t),\ \ \ \ \
z^\mu(0)=(u_0,v_0),
\end{equation}
where  $B_\mu$ and $\mathcal{Q}_\mu$ are the operators
defined in \eqref{pp1} and \eqref{pp2}, respectively.

Note that the adjoint of the operator
$\mathcal{Q}_\mu:H\to \mathcal{H}$ is the
operator $\mathcal{Q}^\star_\mu:\mathcal{H}\to H$
defined by
\[\mathcal{Q}^\star_\mu(u,v)=\frac 1\mu (-\Delta)^{-1} Q v.\]
In particular we have that $\mathcal{Q}_\mu
\mathcal{Q}^\star_\mu:\mathcal{H}\to \mathcal{H}$ is given by
\[
\mathcal{Q}_\mu \mathcal{Q}^\star_\mu\,(u,v)=\frac
1{\mu^2}\,(0,(-\Delta)^{-1} Q^2 v),\ \ \ \ \ (u,v)
\in\,\mathcal{H}.\]

Next, for any $\mu>0$, we introduce the operator $C_{\mu}
\in\,\L^+(\mathcal{H})$ by setting
\[C_{\mu}:=\int_0^\infty S_\mu(s)\, \mathcal{Q}_\mu \mathcal{Q}_\mu ^\star\,
S_\mu^\star(s)\,ds,\] where $\{S_\mu^\star(t)\}_{t\geq 0}$ is the
semigroup generated by $A_\mu^\star$, the adjoint   to the
operator $A_\mu$.

Next proposition  provides an explicit expression for the
operator $C_\mu$.

\begin{Proposition}\cite[Proposition 5.1]{smolu2}
\label{cmu} For every $\mu>0$,  we have
\begin{equation}
\label{cmu20} C_\mu(u,v)=\frac 12 \le((-\Delta)^{-1} Q^2 u,\frac
1\mu(-\Delta)^{-1} Q^2 v\r),\ \ \ \ \ \ (u,v) \in\,\mathcal{H}.
\end{equation}
In particular, if we assume that $Q e_k=\la_k$, for every $k \in\,\mathbb{N}$,  where $\{\la_k\}_{k \in\,\mathcal{N}}$ is a non-negative sequence such that
\[\sum_{k=1}^\infty\frac{\la_k^2}{\a_k}<\infty,\]
then $C_\mu$ is a trace-class operator with
\[\mbox{{\em Tr}}\,
C_{\mu}=\frac 12\le(1+\frac 1\mu\r)\sum_{k=1}^\infty
\frac{\la_k^2}{\a_k}.\]
\end{Proposition}

In particular, 
we have that $C_\mu$ is a non-negative, symmetric trace-class operator in $\mathcal{H}$, and we conclude that the centered Gaussian measure
\[\nu_\mu:=\mathcal{N}(0,C_\mu),\]
of  covariance $C_\mu$ is well defined in $\mathcal{H}$. Moreover, $\nu_\mu$ can be written as a product of two centered gaussian measures on $H$ and $H^{-1}$, respectively,  that is
\[\nu_\nu=\mathcal{N}(0,(-\Delta)^{-1}Q^2/2)\times \mathcal{N}(0,(-\Delta)^{-1}Q^2/2\mu).\]

\subsection{The linear case}

Our aim here is  studying the  invariant measure of the linear system
\begin{equation} \label{linsist}
dz(t)=A_\mu z(t)\,dt+d\,w^{\mathcal{Q}_\mu}(t),\ \ \ \ \
z(0)=(u_0,v_0) \in\,\mathcal{H},
\end{equation}
and showing that the stationary distribution for the solution of
the linear damped wave equation
\begin{equation}
\label{sta1} \mu \frac{\partial^2 u}{\partial t^2}(t,x)=\Delta
u(t,x)-\frac{\partial u}{\partial t}(t,x)+\frac{\partial
w^Q}{\partial t}(t,x),\ \ \ \ \ \ u(t,x)=0,\ \ \ x \in\,\partial
\mathcal{O},
\end{equation}
coincides for all $\mu>0$ with the unique invariant measure of the
stochastic heat equation
\begin{equation}
\label{sta2} \frac{\partial u}{\partial t}(t,x)=\Delta u
(t,x)+\frac{\partial w^Q}{\partial t}(t,x),\ \ \ \ \ \ u(0,x)=0,\
\ \ x \in\,\partial \mathcal{O}.
\end{equation}

\medskip

 In what follows, we shall assume the following conditions on $Q$.
 \begin{Hypothesis}
\label{H6}
The linear operator $Q$ is bounded in $H$ and diagonal with respect to the basis $\{e_k\}_{k \in\,\nat}$ which diagonalizes $A$. Moreover, if $\{\la_k\}_{k \in\,\nat}$ is the corresponding sequence of eigenvalues, we have
\[
  \sum_{k=1}^\infty \frac{\lambda_k^2}{\alpha_k} < + \infty.
\]
\end{Hypothesis}
In particular, if $d=1$ we can take $Q=I$, but if $d>1$ the noise has to colored in space.

\begin{Theorem}\cite[Theorem 5.2]{smolu2}
\label{stationary} Under Hypothesis \ref{H6}, the
Gaussian measure $\nu_\mu=\mathcal{N}(0,C_\mu)$ is the unique invariant
measure  of system \eqref{linsist}, for each $\mu>0$, and for any
$\varphi \in\,C_b(\mathcal{H})$ and $z_0 \in\,\mathcal{H}$
\begin{equation}
\label{chav} \lim_{t\to
\infty}\E^{z_0}\,\varphi(z^\mu(t))=\int_{\mathcal{H}}\varphi(z)\,\mathcal{N}(0,C_\mu)(dz),
\end{equation}
so that $\mathcal{N}(0,C_\mu)$ is ergodic and strongly mixing.

Moreover the Gaussian measure
$\Pi_1 \nu_\mu=\mathcal{N}(0,(-\Delta)^{-1}Q^2/2)$ is the stationary
distribution of \eqref{sta1}. In particular, $\Pi_1 \nu_\mu$ does not
depend on $\mu>0$  and coincides with the unique invariant measure
$\nu$ of the stochastic heat equation \eqref{sta2}.

\end{Theorem}

As stated in Proposition \ref{cmu}, the operator $C_\mu$ is
non-negative, symmetric  and of trace-class on $\mathcal{H}$.
Thus problem \eqref{linsist} admits an invariant measure  of the
form
\[\la_\mu\star \mathcal{N}(0,C_{\mu}),\] where $\la_\mu$ is an invariant
measure for the semigroup $S_\mu(t)$ and $\mathcal{N}(0,C_{\mu})$
is the Gaussian measure, with zero mean and covariance  operator
$C_{\mu}$ (for a proof see e.g. \cite[Theorem 11.7]{dpz1}).
Moreover, as the semigroup $\{S_\mu(t)\}_{t\geq 0}$ is of negative
type (see Proposition \ref{asy}), due to \cite[Theorem
11.11]{dpz1} $\mathcal{N}(0,C_{\mu})$ is the unique  invariant
measure for \eqref{cage}  and \eqref{chav} holds. As well known
this implies that $\mathcal{N}(0,C_{\mu})$ is ergodic and strongly
mixing.

Next, due to \eqref{cmu20} the measure $\mathcal{N}(0,C_\mu)$
defined on $\mathcal{H}$ is the product of two Gaussian
measures, defined respectively on $L^2(\mathcal{O})$ and on
$H^{-1}(\mathcal{O})$. Namely
\[ \mathcal{N}(0,C_\mu)=\mathcal{N}\le(0,
(-\Delta)^{-1}Q^2/2\r)\times \mathcal{N}\le(0,
(-\Delta)^{-1}Q^2/2\mu\r).
\]
In particular the marginal measure  $\Pi_1
\mathcal{N}(0,C_\mu)$ equals  $\mathcal{N}(0,(-\Delta)^{-1}Q^2/2)$,
so that it does not depend on $\mu>0$ and coincides with the
unique invariant measure $\nu$ of the Ornstein-Uhlenbeck process
solving problem \eqref{sta2}.

This allows us to conclude, as the process
$\bar{u}^\mu(t)=\Pi_1 \bar{z}^\mu(t)$ with
\[\bar{z}^\mu(t)=(\bar{u}^\mu(t),\bar{v}^\mu(t)):=\int_{-\infty}^t
S_\mu(t-s)d\bar{w}^{\mathcal{Q}_\mu}(s)\] is the stationary solution to
problem \eqref{sta1} and its distribution does coincides with  the measure $\Pi_1
\mathcal{N}(0,C_\mu)$.

\subsection{The semi-linear case}

We show that an analogous result holds also in the non-linear
case, when \eqref{quasilinear} is a gradient system. To this purpose, we need to assume that the non-linearity $B$ has some special structure.
\begin{Hypothesis}
\label{H7}
There exists $F : H \to \mathbb{R}$ of class $C^1$, with $F(0)=0$, $F(h) \geq 0$ and $\left<D\!F(h),h \right> \geq 0$ for  all $h \in\,H$, such that
\[B(h) = -Q^2 D\!F(h),\ \ \ \ h \in\,H.\]
Moreover, there exists some $\kappa>0$ such that
\[  \left| D\!F(h) - D\!F(k) \right|_H \leq \kappa \left|h - k \right|_H,\ \ \ \ h, k \in\,H.
\]
\end{Hypothesis}

\begin{Example}
\em{ \begin{enumerate}
\item  Assume $d=1$ and take $Q=I$. Let $b: \mathbb{R} \to \mathbb{R}$ be a decreasing Lipschitz continuous function with $b(0)=0$.  Then the composition operator $B(h)(x) = b(h(x))$, $x \in\,\mathcal{O}$, is of gradient type.  Actually, if we set
  \begin{equation*}
    F(h) = -\int_\mathcal{O} \int_0^{h(x)} b(\eta) d \eta d x,\ \ \ \ h \in\,H,
  \end{equation*}
  we have
  \[B(h)=-D\!F(h),\ \ \ \ h \in\,H.\]
  Moreover,
 it is clear that $F(0)=0$, $F(h) \geq 0$ for all $h \in\,H$,  and
 \[\left<D\!F(h),h \right>=-\int_\mathcal{O} b(h(x))h(x)\,dx\geq 0,\ \ \  h \in\,H.\]
   \item
  Assume now $d\geq 1$, so that $Q$ is a general bounded operator in $H$, satisfying Hypothesis \ref{H1}.  Let $b: \mathbb{R} \to \mathbb{R}$ be a function of class $C^1$, with Lipschitz-continuous first derivative, such  that $b(0) = 0$ and $b(x) \geq 0$, for all $x \in\,\reals$. Moreover, the only local minimum of $b$ occurs at $0$.  Let
  \[F(h)= \int_\mathcal{O} b( h(x)) dx,\ \ \ h \in\,H.\]
It is immediate to check that $F(0)=0$ and $F(h)\geq 0$, for all $h \in\,H$.  Furthermore, for any $h \in\,H$
\[D\!F(h)(x)= b^\prime (h(x)),\ \ \ \  x \in\,\mathcal{O}.\]
 Therefore, the nonlinearity
 \[B(h) = -Q^2  b^\prime(h(\cdot)),\ \ \ h \in\,H,\]
satisfies Hypothesis \ref{H7}.
  \end{enumerate}}
\end{Example}

Our aim first is showing that, under the above conditions, system \eqref{cage} is of gradient
type and admits an invariant measure of the following type
\[\nu_\mu(dz)=c_\mu\,e^{-\,2 F(u)}\,\mathcal{N}(0,C_{\mu})(dz),\]
for some normalizing constant $c_\mu$.

To this purpose we introduce some notations.
 For any $n \in\,\nat$ and $\delta \in\,\reals$ we define
 \[T_n:\reals^n\to H^\delta(\mathcal{O}),\ \ \ \ \ (x_1,\ldots,x_n)\mapsto
 \sum_{k\leq n} x_k e_k,\]
 and
 \[\bar{T}_n:\reals^n\times \reals^n\to \mathcal{H}_\delta,\ \ \ \
 \ (x,\eta)\mapsto (T_n x, T_n \eta).\]
 Moreover, we define
\[R_n:H^\delta(\mathcal{O}) \to \reals^n,\ \ \ \ \ u\mapsto
\le(\le<u,e_1\r>_{L^ 2(\mathcal{O})},\ldots,\le<u,e_n\r>_{L^
2(\mathcal{O})}\r).\] Clearly we have $R_n
T_n=\mbox{Id}_{\reals^n}$. Furthermore, if we set
\[P_n:=T_n R_n,\]
we have that $P_n$ is the projection of $H^\delta(\mathcal{O})$
onto the finite dimensional space generated by
$\{e_1,\ldots,e_n\}$ and for any fixed $u
\in\,H^\delta(\mathcal{O})$ we have that $P_n u$ converges to $u$
in $H^\delta$, as $n$ goes to infinity. In particular, setting
\[\bar{P}_n(z):=(P_n u,P_n v), \ \ \ \ \ z=(u,v)
\in\,\mathcal{H}_\delta,\] we have
 \begin{equation} \label{sl1}
\lim_{n\to \infty} \bar{P}_n z=z,\ \ \ \ \ \text{in}\ \
\mathcal{H}_\d.
\end{equation}

In what follows, for any Banach space $X$ we denote by $B_b(X)$
the Banach space of Borel and bounded functions from $X$ into
$\reals$, endowed with the sup-norm, and we denote by $C_b(X)$ the
subspace of uniformly continuous functions.

We recall that the transition semigroup $\{P^\mu(t)\}_{t\geq 0}$
associated with system \eqref{cage} in $\mathcal{H}$ is defined
for any $t\geq 0$ and $\varphi \in\,B_b(\mathcal{H})$ by
\[P^\mu(t)\varphi(z)=\E\,\varphi(u^\mu(t),v^\mu(t)),\ \ \ \ \
z=(u,v) \in\,\mathcal{H},\] where $(u^\mu(t),v^\mu(t))$ is the
solution of \eqref{cage} with initial datum $z=(u,v)$.

Next, we denote by $z^\mu_n(t)$ the solution to the finite
dimensional  problem
\begin{equation}
\label{sl8} dz(t)=\le[A_\mu  z(t)+\bar{P}_n B_\mu(\bar{P}_n
z(t))\r]\,dt+dw^{\mathcal{Q}_{\mu,n}}(t),\ \ \ \ \ z(0)=\bar{P}_n
z, \end{equation}
 where $\mathcal{Q}_{\mu,n}=\bar{P}_n
\mathcal{Q}_{\mu}$. Due to \eqref{sl1},  to the fact that $B_\mu$
is Lipschitz continuous and to   estimate
\eqref{sl5}, it is possible to prove the following approximation
result
\[ \lim_{n\rightarrow
\infty}\E\,|z^\mu_n(t)-z^\mu(t)|^2_{\mathcal{H}}=0,\ \ \ \ \ \
t\geq 0.\] An important consequence of this fact is that the
semigroup $P^\mu(t)$ can be approximated by the semigroup
$P_\mu^n(t)$ associated with equation \eqref{sl8}. Namely,  for any
$\varphi \in\,C_b(\mathcal{H})$ and $t\geq 0$ it holds
\begin{equation}
\label{sl3} \lim_{n\to \infty}P^\mu_n(t)\varphi(z)=\lim_{n\to
\infty}\E\,\varphi(z^\mu_n(t))=P^\mu(t)\varphi(z),\ \ \ \ \ \ z
\in\,\mathcal{H}.
\end{equation}

 Since $DF:H \to
H$ is Lipschitz continuous,  it is not difficult
to check that
\[|F(u)|\leq c\le(1+|u|^2_{H}\r),\ \ \ \
|F(u)-F(v)|\leq
c\,|u-v|_{H}\le(1+|u|_{H}+|v|_{H}\r),\]
so that, in particular, $F:H\to \reals$ is  locally
Lipschitz continuous.

\begin{Remark}
{\em \begin{enumerate} \item The assumption that $F\geq 0$ implies that
\[Z:=\int_{H} e^{\,-2 F(u)}\,\mathcal{N}(0,(-\Delta)^{-1}Q^2/2)(du)<\infty,\]
and
\[Z_n:=\int_{H}
e^{\,-2 F(P_n u)}\,\mathcal{N}(0,(-\Delta)^{-1}Q^2/2)(du)<\infty,\ \ \
\ \ \ n \in\,\nat.\]

\item From the proof of Theorem \ref{gradiente} one sees that it is sufficient to assume a
weaker condition than $F\geq 0$. Namely, what is needed
is that $Z$ and $Z_n$ are  finite and
\[\lim_{n\to \infty}\int_{\mathcal{H}}
\varphi_n(z)e^{\,-2 F(P_n
u)}\,\mathcal{N}(0,C_\mu)\,dz=\int_{\mathcal{H}}
\varphi(z)e^{\,-2 F(u)}\,\mathcal{N}(0,C_\mu)\,dz,\] for any
sequence $\{\varphi_n\}\subset C_b(\mathcal{H})$ uniformly
bounded and pointwise convergent  to some $\varphi
\in\,C_b(\mathcal{H})$.
\end{enumerate}}
\end{Remark}
\begin{Theorem}\cite[Theorem 5.4]{smolu2}
\label{gradiente} Assume that Hypotheses \ref{H6} and \ref{H7} hold.  Then, the probability measure
\[\nu_\mu(dz):=\frac 1{Z}\,e^{\,-2 F(u)}\,\mathcal{N}(0,C_{\mu})(dz),\ \ \ \ \mu>0,\]
 is   invariant  for system \eqref{cage}.

Moreover the distribution
\[\nu(du):=\frac 1{Z}\,e^{\,-2
F(u)}\,\mathcal{N}(0,Q^2(-\Delta)^{-1}/2)(du)\] is stationary for
equation \eqref{quasilinear} and coincides with the unique
invariant measure for the stochastic semi-linear heat equation
\begin{equation}
\label{cagebis}
du(t)=\le[\Delta u(t)-Q^2DF(u(t))\r]\,dt+dw^Q(t),\ \ \ \ u(0)=u_0.
\end{equation}
\end{Theorem}

The proof of the theorem above is obtained by first considering the finite dimensional case and then by a limiting argument.

Actually, for any $\mu>0$, $n \in\,\nat$ and $(q,p) \in\,\reals^n\times
\reals^n$, we  denote by $\zeta^\mu_n(t):=(q^\mu_n(t),p^\mu_n(t))$
the solution of the system in $\reals^n$
\begin{equation}
\label{sl7} \le\{ \begin{array}{l} \ds{\dot{q}^{\,\mu}_n(t)=p^\mu_n(t),\ \ \ \ \ q^\mu_n(0)=q}\\
\vs \ds{\mu\,\dot{p}^{\,\mu}_n(t)=R_n \Delta T_n q^\mu_n(t)-Q_n^2R_n
DF(T_n q^\mu_n(t))-p^\mu_n(t)+ Q_n\dot{w}_n(t),\ \ \ \ \
p^\mu_n(0)=p,}
\end{array}\r.
\end{equation}
where \[w_n(t)=(\beta_1(t),\ldots,\beta_n(t)),\ \ \ \ \ t\geq 0.\]
Here
\[\beta_k(t) =\int_0^t\lambda_k^{-1}\left<e_k, dw^Q(s)\right>_H\]
are independent one dimensional Brownian motions
and \[Q_n(\xi_1,\ldots,\xi_n)=(\la_1\xi_1,\ldots,\la_n\xi_n),\ \ \ \ \xi=(\xi_1,\ldots,\xi_n) \in\,\mathbb{R}^n.\]
The transition semigroup associated with system \eqref{sl7} is
defined for any $\varphi \in\,C_b(\reals^{2n})$ by
\[\hat{P}_n^\mu(t)\varphi(q,p)=\E\,\varphi\le(\zeta^\mu_n(t)\r),\ \ \ \
t\geq 0.\] Note that if we define $F_n(q):=F(T_n q)$, we have
\[DF_n(q)=R_n\, DF(T_n q),\ \ \ \ q \in\,\reals^n,\]
and
\[ \frac 12 \,D \le<R_n \Delta T_n q, q\r>_{\reals^n}=R_n
\Delta T_n q\ \ \ \ q \in\,\reals^n.\] Moreover, since
\[\begin{array}{l}
\ds{\int_{\reals^n}\exp\le(\le<Q_n^{-2}R_n \Delta T_n q,q\r>_{\reals^n}-2\, F(T_n
q)\r)\,dq}\\
\vs
\ds{=c_n\,\int_{\reals^n} e^{-2\,F(T_n
q)}\,\mathcal{N}(0,Q_n^2R_n(-\Delta)^{-1} T_n/2)\,dq,}
\end{array}\] for the obvious
normalizing constant $c_n$,  by a change of variable from
Hypothesis \ref{H3} we have
\[\begin{array}{l}
\ds{\int_{\reals^n}\exp\le(\le<Q_n^{-2}R_n \Delta T_n q,q\r>_{\reals^n}-2\, F(T_n
q)\r)\,dq}\\
\vs
\ds{= c_n\,\int_{H}e^{-2\,F(P_n
u)}\,\mathcal{N}(0,Q^2(-\Delta)^{-1}/2)\,du<\infty.}
\end{array}\]

As a well-known fact,  the Boltzmann distribution
\[\begin{array}{l}
\ds{\hat{\nu}_{\mu,n}(dq,dp)=c_{\mu,n}\exp\le(\textcolor{blue}{-}\le<Q_n^{-2}R_n \Delta T_n
q, q\r>_{\reals^n}-2\, F_n(q)\r)\exp\le(-\mu
|Q_n^{-1}p|^2_{\reals^n}\r)(dq,dp)}\\
\vs \ds{=\frac 1{Z_n}\,
e^{-2\,F_n(q)}\mathcal{N}(0,Q_n^{2}R_n(-\Delta)^{-1}T_n/2)(dq)\times
\mathcal{N}(0,Q_n^2/2\mu)(dp)}
\end{array}\] is invariant for system \eqref{sl7}, so that  for any
$\hat{\varphi} \in\,C_b(\reals^n)$ and $t\geq 0$
\begin{equation}
\label{sl9} \int_{\reals^n\times
\rn}\hat{P}_n^\mu(t)\hat{\varphi}(q,p)\ \hat{\nu}_{\mu,n}(dq,dp)=
\int_{\reals^n\times \rn}\hat{\varphi}(q,p)\
\hat{\nu}_{\mu,n}(dq,dp). \end{equation}

 Now, it is immediate to
check that the $\mathcal{H}$-valued process $\bar{T}_n
\zeta^\mu_n(t)$ coincides with the solution $z^\mu_n(t)$ of the
approximating system \eqref{sl8} with initial datum $\bar{T}_n
(q,p)$. For any $\varphi \in\,C_b(\mathcal{H})$, this yields
\[P^\mu_n(t)\varphi(\bar{T}_n (q,p))=\hat{P}^\mu_n(t)(\varphi \circ
\bar{T}_n)(q,p),\ \ \ \ \ (q,p) \in\,\reals^n\times \reals^n,\]
and hence from \eqref{sl9}  for any $\varphi
\in\,C_b(\mathcal{H})$ we obtain
\begin{equation}
\label{sl10}
\begin{array}{l} \ds{\int_{\reals^{n}\times \reals^n}P_n^\mu(t)\varphi(\bar{T}_n (q,p))\
\hat{\nu}_{\mu,n}(dq,dp)
=\int_{\reals^{n}\times \reals^n} \varphi(\bar{T}_n(q,p))\
\hat{\nu}_{\mu,n}(dq,dp).} \end{array}
\end{equation}
If $T_n$ is considered as a mapping from $\reals^n$ into
$H$ by reasoning as above we have
\begin{equation}
\label{sl13}
\le[e^{-2\,F_n(q)}\mathcal{N}(0,Q_n^2R_n(-\Delta)^{-1}T_n/2)\r]\circ
T_n^{-1}(du)=e^{-2\,F(u)}\mathcal{N}(0,Q^2(-\Delta)^{-1}P_n/2)(du).
\end{equation}
Moreover, if $T_n$ is considered as a mapping from $\reals^n$ into
$H^{-1}(\mathcal{O})$ we have
\begin{equation}
\label{sl14} \mathcal{N}(0,Q_n^2/2\mu)\circ
T_n^{-1}=\mathcal{N}(0,Q^2(-\Delta)^{-1}P_n/2\mu).
\end{equation}
Actually, for any $\la \in\,H^{-1}(\mathcal{O})$ we have
\[\begin{array}{l}
\ds{\int_{H^{-1}(\mathcal{O})}\exp\le(i\,\le<\la,v\r>_{H^{-1}(\mathcal{O})}\r)\,
\le[\mathcal{N}(0,Q_n^2/2\mu)\circ
T_n^{-1}\r]\,dv}\\
\vs \ds{=\int_{\reals^n}\exp\le(i\,\le<(-\Delta)^{-1}\la,T_n
p\r>_{H}\r)\mathcal{N}(0,Q_n^2/2\mu)\,dp}\\
\vs \ds{=\exp\le(-\frac 1{4\mu}
\le<Q_n^2R_n(-\Delta)^{-1}\la,R_n(-\Delta)^{-1}
\la\r>_{\reals^n}\r)}\\
\vs
\ds{=\exp\le(-\frac 1{4\mu}
\le<Q^2(-\Delta)^{-1}P_n\la,P_n \la\r>_{H^{-1}(\mathcal{O})}\r)}\\
\vs
\ds{=\int_{H^{-1}(\mathcal{O})}\exp\le(i\,\le<\la,v\r>_{H^{-1}(\mathcal{O})}\r)
\mathcal{N}(0,Q^2(-\Delta)^{-1}P_n/2\mu)\,dv,}
\end{array}\]
and  by uniqueness of the Fourier transform we obtain
\eqref{sl14}.

Therefore, from \eqref{sl13} and \eqref{sl14} we have
\[\begin{array}{l}
\ds{\hat{\nu}_{\mu,n}\circ \bar{T}_n^{-1}(dz)= \frac
1{Z_n}\,e^{-2\,F(u)} \mathcal{N}(0,Q^2(-\Delta)^{-1}P_n/2)\times
\mathcal{N}(0,Q^2(-\Delta)^{-1}P_n/2\mu)\,(dz)}\\
\vs \ds{=\frac 1{Z_n}\,e^{-2\,F(u)} \le[\mathcal{N}(0,C_\mu)\circ
\bar{P}_n^{-1}\r](dz),}
\end{array}\] and hence, since
\[P^\mu_n(t)\varphi(z)=P^\mu_n(t)\varphi(\bar{P}_nz),\ \ \ \ \ z \in\,\mathcal{H},\]
from \eqref{sl10} it follows
\begin{equation}
\label{sl15}
\begin{array}{l}\ds{ \frac
1{Z_n}\,\int_{\mathcal{H}}P_n^\mu(t)\varphi(z) \,e^{-2\,F(P_n
u)}\mathcal{N}(0,C_\mu)\,dz
=\frac 1{Z_n}\,\int_{\mathcal{H}} \varphi(\bar{P}_n
z)\,e^{-2\,F(P_n u)}\mathcal{N}(0,C_\mu)\,dz.}
\end{array}
\end{equation}

Now, due to \eqref{sl1} and \eqref{sl3} we have
\[\lim_{n\to \infty}P_n^\mu(t)\varphi(z)
\,e^{-2\,F(P_n u)}=P^\mu(t)\varphi(z) \,e^{-2\,F(u)}.\]
 Then, thanks
to Hypothesis \ref{H3}, by the dominated convergence theorem we
can take the limit as $n$ goes to infinity in both sides of
\eqref{sl15} and we get
\[
\begin{array}{l}\ds{ \frac
1{Z}\,\int_{\mathcal{H}}P^\mu(t)\varphi(z) \,e^{-2\,F(
u)}\mathcal{N}(0,C_\mu)\,dz=\frac 1{Z}\,\int_{\mathcal{H}}
\varphi(z)\,e^{2\,U(u)}\mathcal{N}(0,C_\mu)\,dz,}
\end{array}\]
for any $\varphi \in\,C_b(\mathcal{H})$. By a monotone class
argument the same identity follows for arbitrary $\varphi
\in\,B_b(\mathcal{H})$. This in particular implies that the
measure
\[\nu_\mu=\frac 1Z\, e^{-2\,F(u)}\mathcal{N}(0,C_\mu)(dz)\] is
invariant for $P^\mu(t)$.

Finally, since
\[\Pi_1\,[\frac 1Z\, e^{-2\,F(u)}\mathcal{N}(0,C_\mu)(dz)]=\frac
1Z\, e^{-2\,F(u)}\mathcal{N}(0,Q^2(-\Delta)^{-1}/2)(dz),\] we obtain
the second part of the theorem.

\section{Large deviations in the gradient case}
\label{sec7}

We are here interested in the small noise behavior of the solution of the evolution equation \eqref{cage} in $\mathcal{H}$ and in comparing it with the small mass behavior of the solution of the evolution equation \eqref{cagebis}.

Namely, we introduce the equation
\begin{equation}
\label{wave-eq}
dz^\mu_\e(t)=\le[A_\mu z^\mu_\e(t)+B_\mu(z^\mu_\e(t))\r]\,dt+\sqrt{\e}\,dw^{\mathcal{Q}_\mu}(t),\ \ \ \ z^\mu_\e(0)=(u_0,v_0) \in\,\mathcal{H},
\end{equation}
for some parameters $0<\e,\mu\ll1$. We assume that Hypotheses \ref{H6} and \ref{H7} are verified so that the non-linearity $B$ has a gradient structure.

 We keep $\mu>0$ fixed and let $\epsilon$ tend to zero, to study some relevant quantities associated with the large deviation principle for this system, as the quasi-potential that describes also the asymptotic behavior of the expected exit time  from a domain and the corresponding exit places.  Due to the gradient structure of \eqref{wave-eq}, as in the finite dimensional case studied in \cite{freidlin}, we are here able to calculate explicitly the quasi-potentials $V^\mu(u,v)$ for system \eqref{wave-eq}
as
\begin{equation}
\label{m2-i1}
   V^\mu(u,v) = \left|(-\Delta)^{\frac{1}{2}}Q^{-1} u \right|_{H}^2 + 2F(u) + \mu \left|Q^{-1} v \right|_{H}^2.\end{equation}
   Actually, we can prove that for any $\mu>0$
   \begin{equation}
   \label{m2-i2}
   V^{\mu}(u,v) = \inf \left\{ I^{\mu}_{-\infty}(z)\, :\, z(0) = (u,v),\  \lim_{t \to - \infty} \left|C_\mu^{-1/2} z(t) \right|_H = 0 \right\},\end{equation}
   where
   \[I^\mu_{-\infty}(z)=\frac 12 \int_{-\infty}^0\le|Q^{-1}\le( \mu \frac{\partial^2 \varphi}{\partial t^2} (t) - \Delta \varphi(t) + \frac{\partial \varphi}{\partial t} (t) -B(\varphi(t))\r)\r|_{H}^2\,dt\]
   with $\varphi(t)=\Pi_1z(t)$, and $C_\mu$ is the operator defined in \eqref{cmu20}.
From \eqref{m2-i2}, we obtain that
\[   V^\mu(u,v) \geq \left|(-\Delta)^{\frac{1}{2}}Q^{-1} u \right|_{H}^2 + 2F(u) + \mu \left|Q^{-1} v \right|_{H}^2.\]
Thus, we obtain the equality \eqref{m2-i1} by constructing  a path which realizes the minimum.
An immediate consequence of \eqref{m2-i1} is that for each $\mu>0$
\[
\bar{V}_\mu(u)=\inf_{v \in\,H^{-1}(\mathcal{O})}V^\mu(u,v)=V^\mu(u,0)=V(u),\]
where $V(u)$ is the quasi-potential associated with the equation
\begin{equation}
\label{heat-eq}
du(t)=\le[\Delta u(t)-Q^2DF(u(t))\r]\,dt +\sqrt{\e}\,dw^Q(t),\ \ \ \ u(0)=u_0 \in\,H.
\end{equation}

\subsection{A characterization of the quasi-potential}

For any $\mu>0$ and $t_1<t_2$, and for any $z \in\,C([t_1,t_2];\H)$ and $z_0 \in\,\H$, we define
\[
  I^\mu_{t_1,t_2}(z) = \frac{1}{2} \inf \left\{ \left| \psi \right|_{L^2((t_1,t_2);H)}^2 : z = z^\mu_{z_0,\psi} \right\}
\]
where $z^\mu_{\psi}$ solves the  skeleton equation associated with \eqref{wave-eq}
\[
  \frac{dz^\mu_{z_0,\psi}}{dt}(t) = A_\mu z^\mu_{z_0,\psi}(t)+B_\mu(z^\mu_{z_0,\psi}(t))+\sqrt{\e}\mathcal{Q}_\mu\psi(t),\ \ \ \ z^\mu_{z_0,\psi}(0)=z_0 \in\,\mathcal{H}.  \]
Analogously,
for $t_1<t_2$, and for any $\varphi\in\,C([t_1,t_2];H)$ and $u_0 \in\,H$, we define
\[
  I_{t_1,t_2} ( \varphi) = \frac{1}{2} \inf \left\{ \left| \psi \right|_{L^2((t_1,t_2);H)}^2 : \varphi = \varphi_{\psi,u_0} \right\},
\]
where $\varphi_{u_0,\psi}$ solves the  problem
\[
 \frac{ d\varphi_{u_0,\psi}}{dy}(t) = \Delta \varphi_{u_0,\psi}-Q^2DF(\varphi_{u_0,\psi}(t))+Q\,\psi(t),\ \ \ \ \ \varphi_{u_0,\psi}(0)=u_0.
\]
In what follows, we shall also denote
\[  I^\mu_{-\infty}(z) = \sup_{t<0} I^\mu_{t,0}(z),\ \ \ \ \ I_{-\infty}(\varphi) = \sup_{t<0} I_{t,0}(\varphi).
\end{equation*}

Since \eqref{wave-eq} and \eqref{heat-eq} have additive noise, as a consequence of the contraction lemma we have that the family $\{z^{\mu}_{\epsilon}\}_{\epsilon>0}$ satisfies the large deviations principle in $C([0,T];\H)$, with respect to the rate function $I^\mu_{0,T}$ and the family $\{u^\epsilon\}_{\epsilon>0}$ satisfies the large deviations principle in $C([0,T];H)$, with respect to the rate function $I_{0,T}$.

In what follows, for any fixed $\mu>0$ we shall denote by $V^\mu$  the quasi-potential associated with system \eqref{wave-eq}, namely
\[   V^\mu(u,v) = \inf \left\{ I^\mu_{0,T}(z)\, :\, z(0)= 0,\  z(T) = (u,v),\  T>0 \right\}.
 \]
Analogously, we shall denote by $V$ the  quasi-potential associated with equation \eqref{heat-eq}, that is
  \[  V(u) = \inf \left\{ I_{0,T}(\varphi)\,:\, \varphi(0) = 0,\  \varphi(T) = u,\  T>0 \right\}.\]
Moreover, for any $\mu>0$ we shall define
\[\bar{V}_\mu(u)=\inf_{v \in\,H^{-1}(\mathcal{O})} V^\mu(u,v)=\inf \left\{ I^\mu_{0,T}(z)\,:\, z(0) = 0,\  \Pi_1 z(T) = u,\  T>0 \right\}.\]
In \cite[Proposition 5.4]{cerrok}, it has been proven   that $V(u)$ can be represented as
\begin{equation*}
  V(u) = \inf \left\{ I_{-\infty} (\varphi)\,:\, \varphi(0) = u,\  \lim_{t \to -\infty} \left|\varphi(t) \right|_{H} = 0 \right\}.
\end{equation*}
Here, we want   to show a similar representations for $V^\mu(u,v)$, for any fixed $\mu>0$.

\begin{Theorem}\cite[Theorem 3.1]{cs-asy}
  For any $\mu>0$ and $(u,v) \in \H$ we have
 \[   V^{\mu}(u,v) = \inf \left\{ I^{\mu}_{-\infty}(z)\, :\, z(0) = (u,v),\  \lim_{t \to - \infty} \left|C_\mu^{-1/2} z(t) \right|_H = 0 \right\}.
  \]
\end{Theorem}

 First we observe that by the definitions of $I^\mu_{t_1,t_2}$ and $V^\mu(x,y)$,
  \begin{equation*}
    V^\mu(u,v) = \inf \left\{I^\mu_{-T,0}(z): z(-T)=0, z(0) = (u,v), T>0 \right\}.
  \end{equation*}
Now, for any $\mu>0$ and $(u,v) \in\,\H$,  we define
  \begin{equation*}
    M^\mu(u,v):= \inf \left\{ I^{\mu}_{-\infty}(z)\, :\, z(0) = (u,v),\ \lim_{t \to - \infty} \left|C_\mu^{-1/2} z(t) \right|_\H = 0 \right\}.
  \end{equation*}
  Clearly, we want to prove that $M^\mu(u,v)=V^\mu(u,v)$, for all $(u,v) \in\,\H$.

  If $z$ is a continuous path with $z(-T)=0$ and $z(0) = (u,v)$, we can extend it  in $C((-\infty,0);\H)$, by defining $z(t) = 0$, for $t<-T$.  Then, since $D\!F(0)=0$,  we see that
  \[I^\mu_{-\infty}(z) = I^\mu_{-T,0}(z),\]
so that $M^\mu(u,v) \leq V^\mu(u,v)$.

Now, let us prove that the opposite inequality holds.  If $M^\mu(u,v)= +\infty$, there is nothing else to prove.  Thus, we assume that $M^\mu(u,v) < +\infty$.  This means that for any $\e>0$ there must be some  $z^\mu_\e \in C((-\infty,0);\H)$, with the properties that $z^\mu_\e(0) = (u,v)$ and
\[\lim_{t \to -\infty} |C_\mu^{-1/2} z^\mu_\e(t) |_\H = 0,\ \ \ I^\mu_{-\infty}(z^\mu_\e) \leq M^\mu(u,v) + \epsilon.\]

Now, as proven in \cite[Lemma 3.2]{cs-asy}, for any $\mu>0$, there exists  $T_\mu > 0$ such that for any $t_1<t_2-T_\mu$ we have $\text{Im}(L^\mu_{t_1,t_2}) = D(C_\mu^{-1/2})$ and
\[    \left|(L^\mu_{t_1,t_2})^{-1} z \right|_{L^2((t_1,t_2);H)} \leq c(\mu, t_2 - t_1) \left| C_\mu^{-1/2} z \right|_\H,\ \ \ z \in\,\text{Im}(L^\mu_{t_1,t_2}),
  \]
  where $C_\mu$ is the operator in \eqref{cmu20} and
 $L^\mu_{t_1,t_2} : L^2((t_1,t_2):H) \to \H$, is the operator defined as
\[
  L^\mu_{t_1,t_2} \psi = \int_{t_1}^{t_2} \Smu(t_2-s) \Qmu \psi(s) ds.
\]

Thus,  let $T_\mu$ be the constant from \cite[lemma 3.2]{cs-asy} and let $t_\epsilon<0$ be such that
 \[|C_\mu^{- 1/2} z^\mu_\e(t_\epsilon)|_{\H} < \epsilon.\]
Moreover, let
  $\psi^\mu_\e := (L^\mu_{t_\epsilon - T_\mu, t_\epsilon})^{-1} z^\mu_\e(t_\epsilon)$.  We have
  \begin{equation*}
    \left| \psi^\mu_\e \right|_{L^2(t_\epsilon - T_\mu, t_\epsilon;H)} \leq c_\mu \epsilon.
  \end{equation*}
Next, we define
  \begin{equation*}
    \hat{z}^\mu_\e(t) =  \int_{t_\epsilon - T_\mu}^t \Smu(t-s) Q_\mu \psi^\mu_\e(s) ds.
  \end{equation*}
We have
  \begin{equation*}
  \begin{array}{l}
    \displaystyle{ \int_{t_\epsilon - T_\mu}^{t_\epsilon} |\hat{z}^\mu_\e(t) |_\H^2 dt
    \leq \int_{t_\epsilon - T_\mu}^{t_\epsilon} \left( \int_{t_\epsilon - T_\mu}^t \frac{M_\mu}{\mu} e^{-\omega_\mu (t-s)}
    |Q \psi^\mu_\e(s)|_{H^{-1}} ds \right)^2 dt}\\
    \vs
    \displaystyle{ \leq \frac{M_\mu^2}{2 \mu^2 \omega_\mu}\int_{t_\epsilon - T_\mu}^{t_\epsilon}  |Q \psi^\mu_\e|^2_{L^2((t_\epsilon - T_\mu, t\epsilon);H^{-1})} dt \leq T_\mu \frac{M_\mu^2}{2\mu^2 \omega_\mu} |Q \psi^\mu_\e|^2_{L^2((t_\epsilon - T_\mu, t\epsilon);H^{-1})}\leq c_\mu \e^2.}
    \end{array}
  \end{equation*}
Furthermore, $\hat{z}^\mu_\e(t_\epsilon - T_\mu) = 0$ and $\hat{z}^\mu_\e(t_\epsilon) = z(t_\epsilon)$.  Finally, we notice that
\[\begin{array}{l}
\ds{    \hat{z}^\mu_\e(t) = -\int_{t_\epsilon - T_\mu}^t \Smu(t-s) \Qmu Q D\!{F}(\hat{z}^\mu_\e(s)) ds}\\
\vs
\ds{+ \int_{t_\epsilon - T_\mu}^t \Smu(t-s) \Qmu \left( \psi^\mu_\e(s) + Q D\!{F}( \hat{z}^\mu_\e(s)) \right) ds,}
\end{array}
  \]
so that
  \begin{equation*}
    I^\mu_{t_\epsilon - T_\mu, t_\epsilon}(\hat{z}^\mu_\e) = \frac{1}{2} \left| \psi^\mu_\e + Q D\!{F}(\hat{z}^\mu_\e(s)) \right|_{L^2((t_\epsilon - T_\mu, t_\epsilon);H)}^2
    \leq c_\mu \epsilon^2.
  \end{equation*}
  Now if we define
  \begin{equation*}
    \tilde{z}^\mu_\e(t) = \begin{cases}
                    \hat{z}^\mu_\e(t) & \text{ if } t_\epsilon - T_\mu \leq t < t_\epsilon \\
                    z^\mu_\e(t) & \text{ if } t_\epsilon \leq t \leq 0,
                   \end{cases}
  \end{equation*}
  we see that $\tilde{z}^\mu_\e \in C((t_\epsilon - T_\mu,0);\H)$ and
  \begin{equation*}
    V^\mu(u,v) \leq I^\mu_{t_\epsilon - T_\mu,0}(\tilde{z}^\mu_\e) = I^\mu_{t_\epsilon - T_\mu, t_\epsilon}(\hat{z}^\mu_\e) + I^\mu_{t_\epsilon,0}(z^\mu_\epsilon) \leq c_\mu \epsilon^2 + M^\mu(u,v) + \epsilon.
  \end{equation*}
Due to the arbitrariness  of $\epsilon>0$, we can conclude.

\subsection{The main result}

  If $z \in\,C((-\infty,0];\H)$ is such that $I^\mu_{-\infty,0}(z) < +\infty$,  then we have for $\varphi = \Pi_1 z$,
  \begin{equation} \label{I-mu-explicit-eq}
    I^\mu_{-\infty}(z) = \frac{1}{2} \int_{-\infty}^0 \left| Q^{-1} \left( \mu \frac{\partial^2 \varphi}{\partial t^2}(t) + \frac{\partial \varphi}{\partial t}(t) - \Delta \varphi(t) + Q^2 D\!F(\varphi(t)) \right) \right|_H^2 dt.
  \end{equation}
  Actually, if $I^\mu_{-\infty,0}(z) < +\infty$, then there exists $\psi \in\,L^2((-\infty,0);H)$ such that
  $\varphi= \Pi_1 z$ is a weak solution to
  \begin{equation*}
    \mu \frac{\partial^2 \varphi}{\partial t^2} (t) = \Delta \varphi(t) - \frac{\partial \varphi}{\partial t}(t) - Q^2 D\!F(\varphi(t)) + Q \psi.
  \end{equation*}
This means that
  \begin{equation*}
    \psi(t) = Q^{-1} \left( \mu \frac{\partial^2 \varphi}{\partial t^2}(t) + \frac{\partial \varphi}{\partial t}(t) - \Delta \varphi(t) + Q^2 D\!F(\varphi(t))  \right)
  \end{equation*}
  and \eqref{I-mu-explicit-eq} follows.

  By the same argument, if $I_{-\infty,0}(\varphi) < +\infty$, then it follows that
 \[
    I_{-\infty}(\varphi) = \int_{-\infty}^0 \left|Q^{-1} \left(\frac{\partial \varphi}{\partial t}(t) - \Delta \varphi(t) + Q^2 D\!F(\varphi(t)) \right) \right|_H^2 dt.
    \]

\begin{Theorem}\cite[Theorem 4.1]{cs-asy}
\label{m2-10}
 For any fixed $\mu>0$ and $(u,v) \in\,D((-\Delta)^{1/2}Q^{-1})\times D(Q^{-1})$ it holds
\[    V^\mu(u,v) = \left|(-\Delta)^{\frac{1}{2}}Q^{-1} u \right|_H^2 + 2F(u) + \mu \left|Q^{-1} v \right|_H^2.
  \]
  Moreover
  \[
    V(u) = \left|(-\Delta)^{\frac{1}{2}} Q^{-1} u \right|_H^2 + 2F(u).
  \]
In particular,  for any $\mu>0$,
  \begin{equation*}
    \bar{V}_\mu(u) = \inf_{v \in H^{-1}} V^\mu(u,v) = V^\mu(u,0) = V(u).
  \end{equation*}
\end{Theorem}

First, we observe that if $\varphi(t) = \Pi_1 z(t)$, then
  \begin{equation} \label{I-mu-rewrite-eq}
\begin{array}{l}
\ds{I^\mu_{-\infty}(z) = \frac{1}{2} \int_{-\infty}^0 \left|Q^{-1} \left( \mu \frac{\partial^2 \varphi}{\partial t^2}(t) - \frac{\partial \varphi}{\partial t} (t) - \Delta \varphi(t) + Q^2 D\!F(\varphi(t))  \right) \right|_H^2 dt }\\
\vs
\ds{+ 2 \int_{-\infty}^0 \left< Q^{-1}\frac{\partial \varphi}{\partial t} (t), Q^{-1} \left(\mu \frac{\partial^2 \varphi}{\partial t^2} (t) - \Delta \varphi(t) \right) + Q D\!F(\varphi(t)) \right>_H dt.}
  \end{array}
  \end{equation}
Now, if
\[\lim_{t \to -\infty} |C_\mu^{-1/2} z(t)|_\H =0,\]
  then
  \begin{equation*}
    \lim_{t \to -\infty} \left|(-\Delta)^{\frac{1}{2}} Q^{-1} \varphi(t) \right|_H + \left|Q^{-1} \frac{\partial \varphi}{\partial t} (t) \right|_H = 0,
  \end{equation*}
so that
\[  \begin{array}{l}
 \ds{2 \int_{-\infty}^0 \left< Q^{-1}\frac{\partial \varphi}{\partial t} (t), Q^{-1} \left(\mu \frac{\partial^2 \varphi}{\partial t^2} (t) - \Delta \varphi(t) \right) + Q D\!F(\varphi(t)) \right>_H dt}\\
\vs
\ds{= \left| (-\Delta)^{\frac{1}{2}} Q^{-1} \varphi(0) \right|_H^2 + 2 F(\varphi(0)) + \mu \left|Q^{-1} \frac{\partial \varphi}{\partial t}(0) \right|_H^2.}
  \end{array}\]
This yields
\[V^\mu(u,v) \ge \left| (-\Delta)^{\frac{1}{2}} Q^{-1} u \right|_H^2 + 2 F(u) + \mu \left|Q^{-1}  v \right|_H^2.\]

Now, let $\tilde{z}(t)$ be a mild solution of the  problem
  \begin{equation*}
   \tilde{z}(t) =  \Smu(t) (u,-v) - \int_0^t \Smu(t-s) \Qmu Q D\!{F}(\tilde{z}(s)) ds,
  \end{equation*}
and let  $(u,v) \in D(C_\mu^{-1/2})$.
  Then $\tilde{\varphi} (t)= \Pi_1 \tilde{z}(t)$ is a weak solution of the problem
  \begin{equation*}
 \mu \frac{\partial^2  \tilde{\varphi}}{\partial t^2}(t) = \Delta \tilde{\varphi}(t) - \frac{\partial  \tilde{\varphi}}{\partial t} (t) - Q^2 D\!F(\tilde{\varphi}(t)),\ \ \ \  \tilde{\varphi}(0)=u,\ \ \ \frac{ \tilde{\varphi}}{\partial t}(0)=-v.
   \end{equation*}
Assume that
\begin{equation}
\label{rev1}
\lim_{t \to -\infty} \left| C_\mu^{-1/2} \tilde{z}(t) \right|_\H =0.\end{equation}
  Then, if we define $\hat{\varphi}(t) = \tilde{\varphi}(-t)$ for $t \leq 0$,  we see that $\hat{\varphi}(t)$ solves
  \begin{equation*}
    \mu \frac{\partial^2 \hat{\varphi}}{\partial t^2} (t) = \Delta \hat{\varphi}(t) + \frac{\partial \hat{\varphi}}{\partial t} (t) -Q^2 D\!F(\hat{\varphi}(t)),\ \ \ \hat{\varphi}(0)=u,\ \ \ \frac{\partial \hat{\varphi}}{\partial t}(0) = v.
  \end{equation*}
Thanks to \eqref{I-mu-rewrite-eq} this yields
  \begin{equation*}
    I^\mu_{-\infty}(\hat{\varphi}) = \left| (-\Delta)^{\frac{1}{2}} Q^{-1} u \right|_H^2 + 2 F(u) + \mu \left|Q^{-1} v \right|_H^2.
  \end{equation*}
  and then
  \begin{equation*}
    V^\mu(u,v) = \left| (-\Delta)^{\frac{1}{2}} Q^{-1} u \right|_H^2 + 2 F(u) + \mu \left|Q^{-1} v \right|_H^2.
  \end{equation*}

As known, an analogous result holds for $V(u)$.  In what follows, for completeness, we give a proof. We have
\[ \begin{array}{l}
\ds{    I_{-\infty}(\varphi) = \frac{1}{2} \int_{-\infty}^0 \left|Q^{-1} \left( \frac{\partial \varphi}{\partial t} (t) + \Delta \varphi(t) -Q^2 D\!F(\varphi(t))  \right) \right|_H^2 dt } \\
\vs
\ds{+ 2 \int_{-\infty}^0 \left<Q^{-1} \frac{\partial \varphi}{\partial t} (t), Q^{-1} \left( -\Delta \varphi(t) + Q^2 D\!F(\varphi(t)) \right)   \right>_H dt.}
\end{array}
  \]
From this we see that
  \begin{equation*}
    V(u) \ge \left|(-\Delta)^{\frac{1}{2}} Q^{-1} u \right|_H^2 + 2 F(u).
  \end{equation*}
Just as for the wave equation, for $u \in D((-\Delta)^{\frac{1}{2}} Q^{-1})$, we define $\tilde{\varphi}$ to be the solution of
  \begin{equation*}
    \tilde{\varphi}(t) = e^{t\Delta} u - \int_0^t e^{(t-s)\Delta} Q^2 D\!F(\tilde{\varphi}(s)) ds.
  \end{equation*}
We have
  \begin{equation*}
    \lim_{t \to +\infty}|(-\Delta)^{\frac{1}{2}} Q^{-1}\tilde{\varphi}(t)|_H =0.
  \end{equation*}
Then, if we define $\hat{\varphi}(t) = \tilde{\varphi}(-t)$ we get
  \begin{equation*}
    \frac{\partial \hat{\varphi}}{\partial t} (t) = -\Delta \hat{\varphi}(t) + Q^2 D\!F( \hat{\varphi}(t)),
  \end{equation*}
   so that
   \begin{equation*}
     I_{-\infty}(\hat{\varphi}) = \left|(-\Delta)^{\frac{1}{2}} Q^{-1} u \right|_H^2 + 2 F(u)
   \end{equation*}
   and
   \begin{equation*}
     V(u) = \left|(-\Delta)^{\frac{1}{2}} Q^{-1} u \right|_H^2 + 2 F(u).
   \end{equation*}

Therefore, in order to conclude the proof of  Theorem \ref{m2-10}, one has to show that \eqref{rev1} holds.

\begin{Lemma}\cite[Lemma 4.2]{cs-asy}
 \label{go-to-zero-lemma}
  Let $(u,v) \in\, D((-\Delta)^{\frac{1}{2}}Q^{-1})\times D(Q^{-1})$ and let $\varphi$ solve the problem
  \[
 \mu \frac{\partial^2 \varphi}{\partial t^2} (t) = \Delta \varphi(t) - \frac{\partial \varphi}{\partial t} (t) - Q^2 D\!F(\varphi(t)),\ \ \ \ \  \varphi(0) = u,\ \
      \frac{\partial \varphi}{\partial t} (0) = v.
        \]
  Then \[z(t) = \left( \varphi(t), \frac{\partial \varphi}{\partial t}(t) \right) \in\,D((-\Delta)^{\frac{1}{2}})\times D(Q^{-1}),\ \ \ \ t \geq 0,\]
  and
  \[
    \lim_{t \to +\infty} \left|C_1^{-1/2} z(t)  \right|_\H =0.
  \]
\end{Lemma}

\section{Large deviations in the non-gradient case}
\label{sec8}

Here, we assume that the smooth bounded domain $\mathcal{O} \subset \reals^d$ is regular enough so that
\begin{equation}
\label{m18}
\a_k\sim k^{2/d},\ \ \ \ k \in\,\nat,
\end{equation}
where $\{\a_k\}_{k \in\,\mathbb{N}}$ are the eigenvalues of the Laplacian $\Delta$ in $H$.
This happens for example in the case of a  {\em strongly regular} open set $\mathcal{O}$, both with Dirichlet and with Neumann boundary conditions, see \cite[Theorem 1.9.6]{davies}.

So far, we have assumed that the noise $w^Q(t,x)$ is a cylindrical Wiener process that can be written as
\[w^Q(t,x)=\sum_{k=1}^\infty Q e_k(x) \beta_k(t),\]
where $\{e_k\}_{k \in\,\mathbb{N}}$ is the complete orthonormal basis in $H$ that diagonalizes $\Delta$ and $\{\beta_k(t)\}_{k \in\,\mathbb{N}}$ is a sequence of mutually independent Brownian motions, all defined on the same stochastic basis. Moreover, we have assumed that $Q$ is diagonal with respect to the basis $\{e_k\}_{k \in\,\mathbb{N}}$, with $Q e_k=\la_k\,e_k$. In what follows we shall assume the following for the sequence $\{\la_k\}_{k \in\,\mathbb{N}}$.

\begin{Hypothesis}
\label{H8}
If $\{\la_k\}_{k \in\,\nat}$ is the  sequence of eigenvalues of $Q$, corresponding to the eigenbasis $\{e_k\}_{k \in\,\mathbb{N}}$, we have
\begin{equation}
\label{m1}
\frac 1c \a_k^{-\beta}\leq \la_k\leq c\,\a_k^{-\beta},\ \ \ \ k \in\,\nat,
\end{equation}
for some $c>0$ and $\beta>(d-2)/4$.
\end{Hypothesis}
\begin{Remark}
{\em \begin{enumerate}
\item If $d=1$, according to Hypothesis \ref{H8} we can consider space-time white noise, thai is $Q=I$.
\item  Thanks to \eqref{m18}, condition \eqref{m1} implies that if $d\geq 2$, then there exists $\gamma<2d/(d-2)$ such that
\[\sum_{k=1}^\infty \la_k^\gamma<\infty.\]
Moreover
\[\sum_{k=1}^\infty \frac{\la_k^2}{\a_k}<\infty.\]
\item As a consequence of \eqref{m1}, for any $\d \in\,\reals$
\[\textnormal{Dom}((-A)^{\d/2}Q^{-1})=H^{\d+2\beta}\]
and there exists $c_\d>0$ such that for any $x \in\,H^{\d+2\beta}$
\[\frac 1{c_\d}\,|(-A)^{\d/2}Q^{-1}x|_H\leq |x|_{\d+2\beta}\leq c_\d\,|(-A)^{\d/2}Q^{-1}x|_H\]
\end{enumerate}}
\end{Remark}

Concerning the nonlinearity $B$, we shall assume the following conditions.

\begin{Hypothesis}
\label{H9}
For any $\d \in\,[0,1+2\beta]$, the mapping $B:H^\d\to H^\d$ is Lipschitz continuous, with
\[ [B]_{{\tiny{\text{{\em Lip}}}(H^{\delta})}}=:\gamma_\d<\a_1.\]
Moreover $B(0)=0$.
We also assume that $B$ is differentiable in the space $H^{2\beta}$, with
\[\sup_{u \in H^{2\beta}(D)}\|DB(u)\|_{L(H^{2\beta})} = \gamma_{2\beta}.\]
\end{Hypothesis}

\begin{Remark}
{\em \begin{enumerate}
\item   The assumption that $B$ is differentiable is made for convenience to simplify the proof of lower bounds in Theorem \ref{t.82}. We believe that by approximating the Lipschitz continuous $B$ with a sequence of differentiable functions whose $C^1$ semi-norm is controlled by the Lipschitz semi-norm of $B$, the results proved in Theorem \ref{t.82} should remain true.
\item If for any $h \in\,H$
\[B(h)(x)=b(x,h(x)),\ \ \ x \in\,\mathcal{O},\]
as in Section \ref{sec3} and after,
and if we assume that $b(x,\cdot) \in\,C^{2k}(\reals)$, for $k \in\,[\beta+\d/2-5/4,\beta+\d/2-1/4]$, and
\[\frac{\partial^j b}{\partial \si^j}(x,\si)_{|_{\si=0}}=0,\ \ \ x \in\,\overline{\mathcal{O}},\]
then $B$ maps $H^{\d}$ into itself, for any $\d \in\,[0,1+2\beta]$.
The Lipschitz continuity of $B$ in $H^{\d}$ and the bound on the Lipschitz norm, are satisfied if the derivatives of $b(x,\cdot)$ are small enough.
\end{enumerate}
}
\end{Remark}

In Section \ref{sec6} we have compared the small noise asymptotic behavior of the solution of equation \eqref{wave-eq} with the small noise asymptotic behavior of the solution of equation \eqref{heat-eq}. We have seen that in the gradient case (that is, when $B$ satisfies Hypothesis \ref{H7}), the quasi-potential $V^\mu(u,v)$ associated with equation \eqref{wave-eq} is given by
\begin{equation}
\label{cs8-1}
V^\mu(u,v)=|(-\Delta)^{\frac 12}Q^{-1}u|_H^2+2F(u)+\mu\,|Q^{-1}v|_H^2,
\end{equation}
so that
\begin{equation}
\label{cs8-2}
\bar{V}_\mu(u):=\inf_{v \in\,H^{-1}}V^\mu(u,v)=|(-\Delta)^{\frac 12}Q^{-1}u|_H^2+2F(u)=:V(u).
\end{equation}
And $V(u)$ is, as well known, the quasi-potential for equation \eqref{heat-eq}.

In \cite{cs-annals}, we have studied the same problem in the more delicate situation the system under consideration is not of gradient type and hence there is no explicit expression for $V^\mu(u,v)$ and $V(u)$, as in \eqref{cs8-1} and \eqref{cs8-2}.

\subsection{The unperturbed equation} \label{sec:3-unperturbed}

We consider here equation \eqref{wave-eq} with $\e=0$. Namely,
\begin{equation}
\label{m-fine202}
\frac{dz}{dt}(t)=A_\mu z(t)+B_\mu(z(t)),\ \ \ \ z(0)=z_0=(u_0,v_0).
\end{equation}
The solution to \eqref{m-fine202} will be denoted by $z^\mu_{z_0}(t)$. We recall here that  $\gamma_0$ denotes the Lipschitz constant of $B$ in $H$ (see  Hypothesis \ref{H9}).

In what follows, we will need the following bounds and limit for the solution $z^\mu_{z_0}$ to equation \eqref{m-fine202}.

\begin{Lemma}\cite[Lemma 3.1]{cs-annals} \label{X-sup-bounded-lemma}
  If $\mu < (\alpha_1 - \gamma_0)\gamma_0^{-2}$, there exists a constant $c_1(\mu)>0$ such that
  \[
    \sup_{t \ge 0}\left|z^{\mu}_{z_0}(t) \right|_\H +  \left| z^{\mu}_{z_0} \right|_{L^2((0,+\infty);\H)} \leq c_1(\mu) |z_0|_\H,\ \ \ \ z_0 \in\,\H.
  \]
Moreover,  for any $R>0$,
  \[
    \lim_{t \to +\infty} \sup_{|z_0|_\H \leq R} \left| z^{\mu}_{z_0}(t) \right|_\H =0.
  \]
\end{Lemma}

In particular, this means that
 the unperturbed system is uniformly attracted to $0$ from any bounded set in $\H$. Next, we show that if the initial velocity is large enough, $\Pi_1 z^{\mu}_{z_0}$ will leave any bounded set.
Actually, as shown in \cite[Lemma 3.3]{cs-annals},
for any $\mu>0$ and $t>0$, there exists $c_2(\mu,t)>0$ such that
  \[
   \sup_{s \leq t} \left| \Pi_1 \Smu(s) \left(0,v_0 \right) \right|_H \geq c_2(\mu,t) \left|v_0 \right|_{H^{-1}},\ \ \ \  v_0 \in H^{-1}.
  \]

Therefore, by using the mild formulation of the unperturbed equation \eqref{m-fine202},
 we can conclude that the following lower bound estimate holds for the solution of \eqref{m-fine202}.

\begin{Lemma}\cite[Lemma 3.4]{cs-annals} \label{large-init-veloc-causes-exit-lem}
  For any $\mu>0$ and $t> 0$ there exists $c(\mu,t)>0$ such that
  \[
    \sup_{s \leq t} \left|\Pi_1 z^{\mu}_{z_0}(s) \right|_H \ge c(\mu,t) |\Pi_2 z_0|_{H^{-1}},\ \ \ \ z_0 \in\,\H.
  \]
\end{Lemma}

\subsection{The skeleton equation} \label{sec:4-skeleton}

For any  $\mu>0$ and $s< t$ and for any $\psi \in\,L^2((s,t);H)$ we define
\[L^\mu_{s,t} \psi = \int_{s}^{t} \Smu(t -r) \Qmu \psi(r) dr.\]
Clearly $L^\mu_{s,t}$ is a continuous bounded linear operator from $L^2([s,t];H)$ into $\mathcal{H}$. If we define the pseudo-inverse of $L^\mu_{s,t}$ as
\[(L^\mu_{s,t})^{-1}(x)=\arg\!\min\,\le\{|(L^\mu_{s,t})^{-1}(\{x\})|_{L^2([s,t];H)}\r\},\ \ \ \ x \in\,\Im(L^\mu_{s,t}),\]
  we have the following bounds.

\begin{Theorem} \cite[Theorem 4.1]{cs-annals}          \label{L-inverse-thm}
For any $\mu>0$ and $s<t$, it holds
\begin{equation}
\label{m11}
\left | \left(L^\mu_{s,t} \right)^{-1} z \right |_{L^2((s,t);H)} =
    \sqrt{2}\left |   (C_\mu - \Smu(t-s) C_\mu S^{\star}_\mu(t-s))^{-1/2} z \right |_{\mathcal{H}},\ \ \ \ z \in\,\Im(L^\mu_{s,t}),
    \end{equation}
      where
  \[
    C_\mu (u,v) = \left(Q^2(-A)^{-1} u, \frac{1}{\mu} Q^2 (-A)^{-1} v \right),\ \ \ \ (u,v) \in\,\mathcal{H}.
  \]
    Moreover,  for every $\mu>0$ there exists $T_\mu>0$ such that
  \[    \Im(L^\mu_{s,t}) = \Im((C_\mu)^{1/2})=\H_{1+2\beta},\ \ \ \ t-s \geq T_\mu,
  \]
  and
  \[
    \left | (L^\mu_{s,t})^{-1} z  \right |_{L^2((s,t);H)} \leq c(\mu, t-s) \left |z \right |_{\H_{1+2\beta}},\ \ \ z \in\,\H_{1+2\beta},
  \]
  for some constant $c(\mu,r)>0$, with $r\geq T_\mu$.
\end{Theorem}

\begin{Remark}
{\em \begin{enumerate}
\item In fact,  it is possible to show that $\Im(L^\mu_{s,t}) = \Im((C_\mu)^{1/2})$, for all $t-s>0$, by using the explicit representation of $S_\mu^\star(t)$.
\item As $S_\mu(t)$ is of negative type, from \eqref{m11}, it easily follows that
\[|(L^\mu_{-\infty,t})^{-1}z|_{L^2((-\infty,t);H)}=\sqrt{2}\,|C_\mu^{-1/2}z|_\H,\ \ \ z \in\,\text{Im}(L^\mu_{-\infty,t}).
\]
\end{enumerate}}
\end{Remark}

For any $\psi \in L^2((-\infty,0);H^{2\a})$, with $\a \in\,[0,1/2]$, and for any $\mu>0$, we denote by $z^\mu_\psi \in\,C((-\infty,0);\H)$ the solution of the skeleton equation
\begin{equation}
\label{m16}
z^\mu_\psi(t)=\int_{-\infty}^t S_\mu(t-s) B_\mu(z^\mu_\psi(s))\,ds+\int_{-\infty}^tS_\mu(t-s) Q_\mu\psi(s)\,ds,\ \ \ t\in\,\reals.
\end{equation}
\begin{Lemma} \cite[Lemma 4.3]{cs-annals} \label{wave-high-reg-conv-to-zero-lem}
If
\begin{equation}
\label{m31}\lim_{t\to-\infty}|z^\mu_\psi(t)|_\H=0,
\end{equation}
 we have $z^\mu_\psi \in\,C((-\infty,0);\H_{1+2(\a+\beta)})$ and
   \begin{equation}
   \label{m30}
      \lim_{t \to -\infty} \left |z^\mu_\psi(t) \right|_{\H_{1+2(\a+\beta)}} =0.
    \end{equation}
\end{Lemma}

\begin{Remark}
{\em \begin{enumerate}
\item From the previous lemma, we have that if $z^\mu_\psi \in\,C((-\infty,0);\H)$ solves equation \eqref{m16} and limit \eqref{m31} holds, then $z^\mu_\psi(t) \in\,\H_{1+2\beta}$, for any $t\leq 0$. In particular $z^\mu_\psi(0) \in\,\H_{1+2\beta}$.
\item In \cite[Lemma 3.5]{cerrok}, it has been proven that the same holds for
equation \eqref{heat-eq}. Actually, if $\varphi_\psi \in\,C((-\infty,0);H)$ is the solution to
\[    \varphi_\psi(t) = \int_{-\infty}^t e^{(t-s)A} B(\varphi_\psi(s)) ds + \int_{-\infty}^t e^{(t-s)A} Q \psi(s) ds,
  \]
  for $\psi \in\,L^2((-\infty,0);H)$, and
  \[\lim_{t\to-\infty}|\varphi_\psi(t)|_H=0,\]
then $\varphi_\psi \in\,C((-\infty,0);H^{1+2\beta})$ and there exists a constant such that for all $t\leq 0$,
  \[
    \left|\varphi_\psi(t) \right|_{H^{1+2\beta}} \leq c\, |\psi|_{L^2((-\infty,0;H)}.
  \]
  Moreover,
  \[
  \lim_{t\to-\infty}|\varphi_\psi(t)|_{H^{1+2\beta}}=0.
  \]
  \end{enumerate}
 }
\end{Remark}

The next lemma shows that the solution of the skeleton equation \eqref{m16}, decay as $t \to -\infty$ and they depend continuously on the control $\psi$.

\begin{Lemma}\cite[Lemma 4.3]{cs-annals} \label{energy-est-lemma}
Let $\a \in\,[0,1/2]$ and let $\psi_1, \psi_2 \in\,L^2((-\infty,0);H^{2\a})$. In correspondence of each $\psi_i$, let $z^\mu_{\psi_i} \in\,C((-\infty,0);\H_{1+2(\a+\beta)})$ be a solution of equation \eqref{m16}, verifying
 \eqref{m31}.
Then,
$z^\mu_{\psi_i}\in\,L^2((-\infty,0);\H_{1+2(\a+\beta)})$, for $i=1,2$, and
there exist $\mu_0>0$  and $c>0$ such that for any $\mu\leq \mu_0$ and $\tau\leq 0$
\begin{equation}
\label{m20}
\begin{array}{l}
\ds{|z^\mu_{\psi_1}-z^\mu_{\psi_2}|_{L^2((-\infty,\tau);\H_{1+2(\a+\beta)})}^2 + \sup_{t\leq \tau} \left|\mathcal{I}_\mu(z^\mu_{\psi_1}(t) - z^\mu_{\psi_2}(t)) \right|_{\H_{1+ 2(\alpha +\beta)}}^2}\\
\vs
\ds{ \leq c\,|\psi_1-\psi_2|_{L^2((-\infty,\tau);H^{2\a})}^2,}
\end{array}
\end{equation}
where $\mathcal{I}_\mu(u,v)=(u,\sqrt{\mu} v)$, for any $(u,v) \in\,\mathcal{H}$.

\end{Lemma}

\begin{Remark}
{\em
\begin{enumerate}
\item
Notice that, since $B(0)=0$, we have $z^\mu_0=0$, so that from \eqref{m20} we get
\begin{equation}
\label{m43}
\begin{array}{l}
\ds{|z^\mu_\psi|^2_{L^2((-\infty,\tau);\H_{1+2(\a+\beta)})} + \sup_{t\leq \tau}  \left|\mathcal{I}_\mu z^\mu_{\psi}(t) \right|_{\H_{1+ 2(\alpha +\beta)}}^2 \leq c\,|\psi|^2_{L^2((-\infty,\tau);H^{2\a})}},
\end{array}
\end{equation}
for any $\mu\leq \mu_0$ and $\tau\leq 0$.
\item By proceeding as in the proof of Lemma \ref{energy-est-lemma}, it is possible to prove that
\[
|z^\mu_{\psi_1}-z^\mu_{\psi_2}|_{L^2((-\infty,\tau);\H_{2\beta})}^2 + \sup_{t \leq \tau} \left|\mathcal{I}_\mu (z^\mu_{\psi_1}(t) - z^\mu_{\psi_2}(t)) \right|_{\H_{2\beta}}^2\leq c\,|\psi_1-\psi_2|_{L^2((-\infty,\tau);H^{-1})}^2.
\]
and
\[|z^\mu_{\psi}|_{L^2((-\infty,\tau);\H_{2\beta})}^2 + \sup_{t \leq \tau} \left|\mathcal{I}_\mu z^\mu_{\psi}(t)  \right|_{\H_{2\beta}}^2\leq c\,|\psi|_{L^2((-\infty,\tau);H^{-1})}^2.\]
\end{enumerate}
}
\end{Remark}

\subsection{A characterization of the quasi-potential} \label{sec:5-quasipotential}
For any $t_1<t_2$, $\mu>0$ and $z \in\,C((t_1,t_2);\H)$, we define
\begin{equation}
\label{m5}
I^{\mu}_{t_1,t_2}(z) =\frac{1}{2}  \inf \left\{ | \psi|_{L^2((t_1,t_2);H)}^2\,:\, z=z^\mu_{\psi,z_0}\right\},
\end{equation}
 where $z^\mu_{\psi,z_0}$ is a mild solution of the skeleton equation associated with equation  \eqref{wave-eq}, with deterministic control $\psi \in\,L^2((t_1,t_2);H)$ and initial conditions $z_0$, namely
\[
\frac{dz^\mu_{\psi,z_0}}{dt}(t)=A_\mu z^\mu_{\psi,z_0}(t)+B_\mu(z^\mu_{\psi,z_0}(t))+Q_\mu\psi(t),\ \ \ \ t_1\leq t\leq t_2.
\]
For $\e,\mu>0$ and $z_0 \in\,\H$ we denote by $z^\mu_{\e,z_0} \in\,L^2(\Omega;C([0,T];\H))$ the mild solution of equation \eqref{wave-eq}.
Since the mapping $B_\mu:\mathcal{H}\to\mathcal{H}$ is Lipschitz-continuous and the noisy perturbation in \eqref{wave-eq} is of additive type, as an immediate consequence of  the contraction lemma, for any fixed $\mu>0$ the family $\{\mathcal{L}(z^\mu_{\e,z_0})\}_{\e>0}$  satisfies a large deviation principle in $C([t_1,t_2];\H)$, with action functional  $I^{\mu}_{t_1,t_2}$. In particular,
for any $\delta>0$ and $T>0$,
  \[
    \liminf_{\epsilon \to 0}  \epsilon \log \left( \inf_{z_0 \in \H} \P \left( \left|z^\mu_{\e,z_0}- z^{\mu}_{\psi,z_0} \right|_{C([0,T];\H)} < \delta \right) \right) \ge - \frac{1}{2} |\psi|_{L^2((0,T);H)}^2
  \]
  and, if $K^\mu_{0,T}(r) = \{z \in C([0,T];\H): I^\mu_{0,T}(z) \leq r \}$,
  \[
    \limsup_{\epsilon \to 0} \epsilon \log \left( \sup_{z_0 \in \H} \P \left( \textnormal{dist}_\H(z^{\mu}_{\epsilon,z_0}, K^\mu_{0,T}(r)) > \delta  \right) \right) \leq - r.
  \]

Analogously, if for any $\e>0$ $u_\e$ denotes the mild solution of equation \eqref{heat-eq}, the family $\{\mathcal{L}(u_\e)\}_{\e>0}$  satisfies a large deviation principle in $C([t_1,t_2];H)$ with action functional
\[
I_{t_1,t_2}(\varphi) = \inf \left\{ \frac{1}{2} \left| \psi \right|_{L^2([t_1,t_2];H)}^2\, :\,
                               \varphi=\varphi_{\psi}\r\},\]
where     $\varphi_{\psi}$ is a mild solution of the skeleton equation associated with equation \eqref{heat-eq}
\[\frac{du}{dt}(t)=A u(t)+B(u(t))+Q\psi(t),\ \ \ \ t_1\leq t\leq t_2.\]
In particular, the functionals $I^{\mu}_{t_1,t_2}$ and $I_{t_1,t_2}$ are lower semi-continuous and have compact level sets. Moreover, it is not difficult to show that for any compact sets $E\subset H$ and $\mathcal{E}\subset \mathcal {H}$, the level sets
\[ K_{E,t_1,t_2}(r) = \left\{\varphi \in C([t_1,t_2];H)\ ;\  I_{t_1,t_2}(\varphi) \leq r,\ \varphi(t_1) \in\, E \right\}
\]
and
\[
K^\mu_{\mathcal{E},t_1,t_2}(r) = \left\{ z \in C([t_1,t_2];\H)\ ;\  I^{\mu}_{t_1,t_2}(z) \leq r,\ z(t_1) \in \,\mathcal{E} \right\}
\]
are compact.

In what follows, for the sake of brevity, for any $\mu>0$  and $t \in\,(0,+\infty]$  we shall define $I^{\mu}_t:=I^{\mu}_{0,t}$ and $I^{\mu}_{-t}:=I^{\mu}_{-t,0}$ and, analogously, for any $t \in\,(0,+\infty]$ we shall define $I_t:=I_{0,t}$ and $I_{-t}:=I_{-t,0}$. In particular, we shall set
\[I^\mu_{-\infty}(z)=\sup_{t>0}I^\mu_{-t}(z),\ \ \ \ \ I_{-\infty}(\varphi)=\sup_{t>0}I_{-t}(\varphi).\]
Moreover, for any $r>0$ we shall set
\[K^\mu_{-\infty}(r)=  \left\{ z \in C((-\infty,0];\H) \ ;\  \lim_{t \to -\infty} |z(t)|_\H =0,\   I^\mu_{-\infty}(z) \leq r \right\}\]
and
\[  K_{-\infty}(r) = \left\{\varphi \in C((-\infty,0];H)\ ;\  \lim_{t \to -\infty} |\varphi(t)|_H =0,\  I_{-\infty}(\varphi) \leq r \right\}.\]

Once we have introduced the action functionals $I^{\mu}_{t_1,t_2}$ and $I_{t_1,t_2}$, as we have already seen in Section \ref{sec6} we can introduce the corresponding {\em quasi-potentials}, by setting for any $\mu>0$ and $(u,v) \in\,\H$
\[V^\mu(u,v) = \inf \left\{ I^\mu_{0,T}(z) \ ;\  z(0) = 0,\  z(T) =(u,v), T > 0 \right\},\]
and for any $u \in\,H$
\[  V(u) = \inf \left\{ I_{0,T}(\varphi)\ ;\  \varphi(0)=0,\  \varphi(T) = u, \ T\ge 0 \right\}.\]
Moreover, for any $\mu>0$ and $u \in\,H$, we shall define
  \[
    \bar{V}_\mu(u) = \inf_{v \in H^{-1}} V^\mu(u,v).
  \]
  In \cite[Proposition 5.1]{cerrok} it has been proved that the level set $K_{-\infty}(r)$ is compact in the space $C((-\infty,0];H)$, endowed with the uniform convergence on bounded sets, and in \cite[Proposition 5.4]{cerrok} it has been proven that
  \[V(u)=\min \le\{\,I_{-\infty}(\varphi)\ ;\ \varphi \in\,C((-\infty,0];H),\ \lim_{t \to -\infty} |\varphi(t)|_H =0, \ \varphi(0)=u\,\r\}.\]
In what follows we want to prove an analogous result for $K^\mu_{-\infty}$, $V^\mu(u,v)$ and $\bar{V}_\mu(u)$.

\begin{Theorem} \cite[Theorem 5.1]{cs-annals}\label{compact-level-sets-half-line-thm}
  For small enough $\mu>0$, the level sets $K^\mu_{-\infty}(r)$ are compact in the topology of uniform convergence on bounded intervals.
\end{Theorem}

As a consequence of previous theorem, we have that there exists $\mu_0>0$ such that  for any $\psi \in L^2((-\infty,0);H)$ and $\mu\leq \mu_0$ there exists $z^\mu_\psi \in C((-\infty,0];\H)$ such that
 \begin{equation}
 \label{m56}
   z^\mu_\psi(t)= \int_{-\infty}^t \Smu(t-s)  B_\mu(z^\mu_\psi(s))  ds+ \int_{-\infty}^t \Smu(t-s) \Qmu \psi(s)  ds,\ \ \ \ t\leq 0.
   \end{equation}
Moreover,
 \begin{equation}
 \label{rev2}
   \lim_{t \to -\infty} |z_\psi^\mu(t)|_\H =0.
 \end{equation}

Actually, a standard fixed point argument shows that for any $\mu>0$ and $ N \in\,\nat$ there exists $z^{\mu}_{N} \in C([-N,0];\H)$ satisfying
   \[
     z^{\mu}_{N}(t) = \int_{-N}^t \Smu(t-s) B_\mu(z^{\mu}_{N}(s))  ds +\int_{-N}^t \Smu(t-s) \Qmu \psi(s)  ds.
   \]
Each $z^\mu_N$ can be seen as an element of $C((-\infty,0];H)$, just  by extending it to  $z^{\mu}_{N}(t) = 0$, for all $t< -N$.
   According to  Theorem \ref{compact-level-sets-half-line-thm}, there exists a subsequence $\{z^\mu_{N_k}\}$ converging to  some   $z^\mu \in K^\mu_{-\infty}\left(\frac{1}{2} |\psi|_{L^2((-\infty,0);H)}^2\right)$, uniformly on compact sets.  We notice that for any fixed $N_0 \in\,\nat$ and $t \ge -N_0$
   \[
     z^{\mu}_{N}(t) = \Smu(t+ N_0)z^\mu_N(-N_0)+\int_{-N_0}^t \Smu(t-s) B_\mu(z^{\mu}_{N}(s))  ds +\int_{-N_0}^t \Smu(t-s) \Qmu \psi(s)  ds.
   \]
Therefore, by taking the limit as $N \to +\infty$, we obtain
   \[
     z^{\mu}(t) = \Smu(t+ N_0)z^{\mu}(-N_0)+ \int_{-N_0}^t \Smu(t-s) B_\mu(z^{\mu}(s))  ds +\int_{-N_0}^t \Smu(t-s) \Qmu \psi(s)  ds.
   \]
Finally, if we let $N_0 \to +\infty$, we see that $z^\mu$ solves equation \eqref{m56}.

\medskip

As $K_{-\infty}(r)$ is compact in $C((-\infty,0];H)$ with respect to the uniform convergence on bounded intervals, we have analogously that
for any $\varphi \in\,L^2((-\infty,0)$ there exists $\varphi_\psi \in\,C((-\infty,0];H)$ such that
\[    \varphi_\psi(t) = \int_{-\infty}^t e^{(t-s)A} B(\varphi(s)) ds + \int_{-\infty}^t e^{(t-s)A} Q \psi(s) ds,
 \]
 and
 \[\lim_{t\to-\infty}|\varphi_\psi(t)|_H=0.\]

In \cite{cerrok}, it has been proved that the $V(u)$ can be characterized as
  \[
    V(u) = \inf \left\{ I_{-\infty}(\varphi) : \lim_{t \to -\infty} \varphi(t) = 0,  \varphi(0) = u \right\}.
  \]
  The next crucial result shows that an analogous result holds for $V^\mu(u,v)$ and $\bar{V}_\mu(u)$, at least for $\mu$ sufficiently small.

\begin{Theorem} \cite[Theorem 5.3]{cs-annals}
\label{quasipotential-representation-thm}
For small enough $\mu>0$, we have the following representation for the quasi-potentials $V^\mu(u,v)$
  \[
    V^\mu(u,v) = \min \left\{ I^\mu_{-\infty}(z): \lim_{t \to -\infty} |z(t)|_\H =0, z(0)=(u,v) \right\},
  \]
and for $\bar{V}_\mu(u)$
  \begin{equation}
  \label{m68}
    \bar{V}_\mu(u) = \min \left\{ I^\mu_{-\infty}(z): \lim_{t \to -\infty} |z(t)|_\H =0, \Pi_1z(0) = u \right\},
  \end{equation}
  whenever these quantities are finite.
\end{Theorem}

  From the definitions of $I^\mu_{t_1,t_2}$, it is clear that
 \[
    V^\mu(u,v) = \inf \left\{ I^\mu_{t_1,0}(z): z(t_1) = 0, z(0) =(u,v), t_1 \leq 0 \right\}.
  \]
Now, if we define
  \[
    M^\mu(u,v) = \inf \left\{ I^\mu_{-\infty}(\varphi): \lim_{t \to -\infty} |z(t)|_\H =0, z(0) = (u,v)  \right\},
  \]
it is immediate to check that $M^\mu(u,v) \leq V^\mu(u,v)$, for any $(u,v) \in\,\H$.  To see this, we observe that if $z \in C([t_1,0];\H)$, with $z(t_1) = 0$ and $z(0)=(u,v)$, then
  \[
    \hat{z} (t) = \begin{cases}
                            0, &t \leq t_1 \\
                            z(t), & t_1 < t \leq 0
                          \end{cases}
  \]  has the property that $\hat{z}(0)=(u,v)$, and $|\hat{z}(t)|_{\H}\to 0$, as $t\to -\infty$. Moreover,
  \[I^\mu_{-\infty}(\hat{z}) = I^\mu_{t_1,0}(z).\]
Therefore, we need to show that $V^\mu(u,v) \leq M^\mu(u,v)$, for all $(u,v) \in\,\H$.

\medskip

The characterization of $V^\mu(u,v)$ and $\bar{V}_\mu(u)$ given in Theorem \ref{quasipotential-representation-thm}, implies that $V^\mu$ and $\bar{V}_\mu$ have compact level sets.

More precisely,
for any $\mu>0$ and $r\geq 0$ the level sets
\[K^\mu(r)=\le\{ (u,v) \in\,\H\,:\ V^\mu(u,v)\leq r\r\}\]
and
\[K_\mu(r)=\le\{ x \in\,H\,:\ \bar{V}_\mu(u)\leq r \r\}\]
are compact, in $\H$ and $H$, respectively.

Actually, let $\{(u_n,v_n)\}_{n \in\,\nat}\subset K^\mu(r).$ In view of Theorem \ref{quasipotential-representation-thm}, for each $n \in\,\nat$ there exists $z^n \in\,C((-\infty,0];\H)$, with $z^n(0)=(u_n,v_n)$, and $|z^n(t)|_H\to 0$, as $t\downarrow -\infty$, such that $V^\mu(u_n,v_n)=I^\mu_{-\infty}(z^n).$ As $I^\mu_{-\infty}(z^n)\leq r$ and the level sets of $I^\mu_{-\infty}$ are compact in $C((-\infty,0];\H)$, as shown in Theorem \ref{compact-level-sets-half-line-thm}, there exists a subsequence $\{z^{n_k}\}\subseteq \{z^n\}$ converging to some $\hat{z} \in\,
C((-\infty,0];\H)$, with $I_{-\infty}^\mu(\hat{z})\leq r$. Since
\[\lim_{k\to\infty}(u_{n_k},v_{n_k})=\lim_{k\to\infty}z^{n_k}(0)=\hat{z}(0)=:(\hat{u},\hat{v}),\ \ \ \text{in}\,\H,\]
due to Theorem \ref{quasipotential-representation-thm} we have
\[V^\mu(\hat{u},\hat{v})\leq I^\mu_{-\infty}(\hat{z})\leq r,\]

\subsection{Continuity of $V^\mu$ and  $\bar{V}_\mu$} \label{sec:6-continuity-of-V}
As a consequence of the compactness of the level sets of $V^\mu$ and $\bar{V}_\mu$,  the mappings $V^\mu:\H\to[0,+\infty]$ and  $\bar{V}_\mu:H\to [0,+\infty]$ are lower semicontinuous.
In fact, it is possible to show that this implies that  the mappings
\[V^\mu:\H_{1+2\beta}\to [0,+\infty),\ \ \ \bar{V}_\mu:H^{1+2\beta}\to[0,+\infty)\]
are well defined and  continuous, uniformly in $0<\mu<1$.

Actually, there exists $c>0$ such that for any $\mu>0$ and $(u,v) \in\,\H_{1+2\beta},$ we have
  \begin{equation}
  \label{m73}
    V^\mu(u,v) \leq  c(1+\mu+\mu^2) |(u,v)|_{\H_{1+2\beta}}^2 
    \end{equation}
  and
  \begin{equation}
  \label{m74}
    \bar{V}_\mu(u) \leq c(1+\mu) \left|u \right|_{H^{1+2\beta}}^2.
  \end{equation}

First of all, we notice that \eqref{m74} is a consequence of \eqref{m73} and of the way $\bar{V}_\mu(x)$ has been defined. The proof of \eqref{m73} is based on the fact that
  \begin{equation*}
    V^\mu(u,v) \leq I^\mu_{-\infty} ( \Pi_1 \Smu(- \cdot) (u, -v))
  \end{equation*}
  and
  \begin{equation*}
    \bar{V}_\mu(u) \leq I^\mu_{-\infty} ( \Pi_1 \Smu(- \cdot)(u,0)).
  \end{equation*}
Now,  if we set $z(t) =\Smu(-t)(x,-y)$ and $\varphi(t)=\Pi_1z(t)$, due to Hypothesis \ref{H2} we have
  \[\begin{array}{l}
\ds{I^\mu_{-\infty}(z) = \frac{1}{2} \int_{-\infty}^0 \left|Q^{-1} \left( \mu \frac{\partial^2\varphi}{\partial t^2}(t) + \frac{\partial\varphi}{\partial t}(t) -A \varphi(t) - B(\varphi(t)) \right) \right|_H^2 dt }\\
\vs
\ds{\leq \int_{-\infty}^0 \left| Q^{-1} \left( \frac{\partial^2\varphi}{\partial t^2}(t) + \frac{\partial\varphi}{\partial t}(t) -A \varphi(t) \right) \right|_H^2 dt + c\,\gamma_{2\beta}^2 \int_{-\infty}^0 \left|\varphi(t) \right|_{H^{2\beta}}^2 dt.}
\end{array}
  \]
Therefore, we get  \eqref{m73} from the following lemma.

\begin{Lemma}\cite[Lemma 6.1]{cs-annals}
  Let  us fix $(u,v) \in\,\H_{1+2\beta}$ and   $\mu>0$ and let
$z(t) = \Smu(-t)(u,-v),$ $t \leq 0.$
  Then, if we denote $\varphi(t)=\Pi_1z(t),$ we have that   $\varphi$ is a weak solution to
    \[
 \le\{     \begin{array}{l}
\ds{        \mu \frac{\partial^2\varphi}{\partial t^2}(t) = A \varphi(t) + \frac{\partial \varphi}{\partial t}(t),\ \ \ \  t\leq 0}\\
\vs
\ds{\varphi(0)  =u, \quad \frac{\partial \varphi}{\partial t}(0) = v.}
      \end{array}\r.
    \]
and
\[\frac{1}{2} \int_{-\infty}^0 \left| Q^{-1} \left( \mu \frac{\partial^2\varphi}{\partial t^2}(t) + \frac{\partial\varphi}{\partial t}(t) - A \varphi(t) \right) \right|_H^2 dt = \left|(-A)^{\frac{1}{2}}Q^{-1} u \right|_H^2 + \mu \left| Q^{-1} v \right|_H^2.
    \]
Moreover, $\varphi \in L^2((-\infty,0);H^{1+2\beta})$ and
    \[
          \int_{-\infty}^0 \left|\varphi(t) \right|_{H^{1+2\beta}}^2 dt \leq c\,(1+\mu+\mu^2)|(u,v)|_{\H_{1+2\beta}}^2.
    \]
\end{Lemma}

Moreover, we have  the continuity of  $V^\mu$ and $\bar{V}_\mu$.

\begin{Theorem}\cite[Theorem 6.3]{cs-annals}
\label{continuity-of-V-mu-tilde-thm}
For each $\mu>0$ the mappings $V^\mu:\H_{1+2\beta}\to [0,+\infty)$ and $\bar{V}_\mu:H^{1+2\beta}\to [0,+\infty)$ are well defined and continuous. Moreover,
  \[
      \lim_{n \to \infty} \left|(u,v) - (u_n,v_n)\right|_{\H_{1+2\beta}}=0\Longrightarrow    \lim_{n \to \infty} \sup_{0 < \mu < 1} \left| V^\mu(u,v) - V^\mu(u_n,v_n) \right|   =0.
  \]
and
  \[\lim_{n \to \infty} \left|u - u_n\right|_{H^{1+2\beta}}=0\Longrightarrow   \lim_{n \to \infty} \sup_{0 < \mu < 1} \left| \bar{V}_\mu(u) - \bar{V}_\mu(u_n) \right|   =0.
  \]
\end{Theorem}

\subsection{Upper bound} \label{sec:7-upper}

In this section we show that for any closed set $N \subset H$
\[
  \limsup_{\mu \downarrow 0} \inf_{u \in N} \bar{V}_\mu(u) \leq \inf_{u \in N} V(u).
\]
First of all, we notice that if $I_{-\infty}(\varphi)<\infty$, then
\begin{equation}
\label{m48}
\varphi \in\,L^2((-\infty,0);H^{2(1+\beta)}),\ \ \ \frac{\partial \varphi}{\partial t} \in\,L^2((-\infty,0);H^{2\beta}),
\end{equation}
and
  \begin{equation}
  \label{m49}
    I_{-\infty}(\varphi) =\frac{1}{2} \int_{-\infty}^0 \left|Q^{-1} \left( \frac{\partial\varphi}{\partial t}(t) - A \varphi(t) -
    B(\varphi(t)) \right) \right|_H^2 dt.
  \end{equation}
Actually, if $\varphi$ solves
\begin{equation*}
    \varphi(t) = \int_{-\infty}^t e^{(t-s)A} B(\varphi(s)) ds + \int_{-\infty}^t e^{(t-s)A} Q \psi(s) ds
  \end{equation*}
  then we can check that \eqref{m48} holds and
    \begin{equation*}
    \psi(t) = Q^{-1} \left( \frac{\partial\varphi}{\partial t} (t) - A \varphi(t) - B(\varphi(t)) \right),
  \end{equation*}
so that \eqref{m49} follows.
Moreover, if
\[  \varphi \in\,L^2((-\infty,0);H^{2(1+\beta)}),\ \ \ \ \ \ \ \frac{\partial \varphi}{\partial t}, \frac{\partial^2 \varphi}{\partial t^2}\in\,L^2((-\infty,0);H^{2\beta}),\]
then
\[
    I^\mu_{-\infty} (z) = \frac{1}{2} \int_{-\infty}^0 \left|Q^{-1} \left( \mu \frac{\partial^2\varphi}{\partial t^2}(t) + \frac{\partial \varphi}{\partial t}(t) - A \varphi(t) - B(\varphi(t)) \right) \right|_H^2 dt,
  \]
  where
  \[z(t)=(\varphi(t),\frac{\partial \varphi}{\partial t}(t)).\]
Actually, if $I^\mu_{-\infty}(z)<\infty$, then $z$ solves
  \begin{equation*}
    z(t) = \int_{-\infty}^t \Smu(t-s) B_\mu(z(s))  ds
    +\int_{-\infty}^t \Smu(t-s) \Qmu \psi(s)ds
  \end{equation*}
  so that
  \begin{equation*}
    \psi(t) = Q^{-1} \left(\mu \frac{\partial^2\varphi}{\partial t^2} (t) + \frac{\partial\varphi }{\partial t} (t) - A \varphi(t)  - B(\varphi(t)) \right)
  \end{equation*}
  weakly.

In particular, as in \cite{cf}, where the finite dimensional case is studied, this means
\begin{equation} \label{I-mu=I+remainder}
\begin{array}{l}
\ds{  I^\mu_{-\infty} (z) = I_{-\infty}(\varphi) + \frac{\mu^2}{2} \int_{-\infty}^0 \left| Q^{-1} \frac{\partial^2\varphi}{\partial t^2}(t) \right|_H^2 dt}\\
\vs
\ds{  + \mu \int_{-\infty}^0 \left< Q^{-1} \frac{\partial^2\varphi}{\partial t^2}(t), Q^{-1} \left( \frac{\partial\varphi}{\partial t}(t) - A \varphi(t) - B(\varphi(t)) \right) \right>_H dt,}
\end{array}
\end{equation}
where $\varphi(t)=\Pi_1 z(t)$,
as long as all of these terms are finite.

Now, for any $\mu>0$ let us define
\begin{equation}
\label{m83}
  \rho_\mu(t) = \frac{1}{\mu^\alpha} \,\rho \left(\frac{t}{\mu^{\alpha}} \right),\ \ \ \ t \in\,\mathbb{R},
\end{equation}
for some $\alpha>0$ to be chosen later, where
 $\rho \in C^\infty(\mathbb{R})$ is the usual mollifier  function such that
\[\text{supp}(\rho) \subset \subset [0,2] ,\ \ \ \int_\mathbb{R} \rho(s) ds = 1,\ \ \ 0 \leq \rho \leq 1.\]
This scaling ensures that
\[
  \int_\mathbb{R} \rho_\mu(s) ds = 1.
\]
Next, we define $\varphi_\mu$ as the convolution for $t<0$
\begin{equation} \label{convolution-def}
  \varphi_\mu(t) = \int_{-\infty}^0 \rho_\mu(t-s) \varphi(s) ds.
\end{equation}

\begin{Lemma}\cite[Lemma 7.1]{cs-annals} \label{x_mu-conv-to-x-lem}
Assume that
\[\varphi \in L^2((-\infty,0);H^{2(1+\beta)}) \cap C((-\infty,0];H^{1+2\beta}),\ \ \   \frac{\partial\varphi}{\partial t} \in L^2((-\infty,0);H^{2\beta})\]
with
\[\varphi(0) =x \in H^{1+2\beta},\ \ \
\lim_{t\to -\infty}\left| \varphi(t) \right|_{H^{1+2\beta}}= 0.\]
Then,
\[
\varphi_\mu \in L^2((-\infty,0);H^{2(1+\beta)}) \cap C((-\infty,0];H^{1+2\beta}),\ \ \   \frac{\partial\varphi_\mu}{\partial t} \in L^2((-\infty,0);H^{2\beta}),
\]
and
\[   \lim_{t \to -\infty} \sup_{\mu >0} \left|\varphi_\mu(t) \right|_{H^{1+2\beta}} =0.
  \]
Moreover,
\[\frac{\partial^2\varphi_\mu}{\partial t^2} \in L^2((-\infty,0);H^{2\beta})\]
 and  for all $\mu >0$,
  \begin{equation}
  \label{m81}
    \left|\frac{\partial^2\varphi_\mu}{\partial t^2} \right|_{L^2((-\infty,0);H^{2\beta})} \leq \frac{c}{\mu^\alpha} \left|\frac{\partial\varphi_\mu}{\partial t}\right|_{L^2((-\infty,0);H^{2\beta})}.
  \end{equation}
  \end{Lemma}

Moreover, the following approximation results hold.
\begin{Lemma}\cite[Lemma 7.2]{cs-annals}
Under the same assumptions of Lemma \ref{x_mu-conv-to-x-lem},
we have
  \begin{equation}
  \label{m82}
    \lim_{\mu \to 0} \left|x- \varphi_\mu(0) \right|_{H^{1+2\beta}} =0,
  \end{equation}
  and
   \[    \lim_{\mu \to 0} \sup_{t \leq 0} \left|\varphi_\mu(t) - \varphi(t) \right|_{H^{1+2\beta}} =0.
  \]
  Moreover,
 \begin{equation} \label{convolution-L2-convergence-eq}
    \lim_{\mu \to 0} |\varphi_\mu - \varphi|_{L^2((-\infty,0);H^{2(1+\beta)})}=0,
  \end{equation}
  and
  \begin{equation} \label{convolution-deriv-L2-convergence-eq}
    \lim_{\mu \to 0} \left| \frac{\partial {\varphi}_\mu}{\partial t} -\frac{\partial {\varphi}}{\partial t}\right|_{L^2((-\infty,0);H^{2\beta})}=0.
  \end{equation}
\end{Lemma}

Using these estimates we can prove the main result of this section.

\begin{Theorem}\cite[Theorem 7.3]{cs-annals}
 \label{upper-bound-thm}
  For any $u \in H^{1+2\beta}$ we have
  \begin{equation}
  \label{m84}
    \limsup_{\mu \downarrow 0} \bar{V}_\mu(u) \leq V(u).
  \end{equation}
\end{Theorem}

  Let $\varphi$ be the minimizer of $V(u)$.  This means $\varphi(0)=u$, \eqref{m49} holds  and $I_{-\infty}(\varphi) = V(u).$  For each $\mu>0$, let $\varphi_\mu$ be the convolution given by \eqref{convolution-def} and let $u_\mu = \varphi_\mu(0)$.

  It is clear that
  \[
    \bar{V}_\mu(u_\mu) \leq  I^\mu_{-\infty}(z_\mu),
  \]
  where
  \[z_\mu(t)=(\varphi_\mu(t),\frac{\partial \varphi_\mu}{\partial t}(t)),\ \ \ t\leq 0.\]
  According to Lemma \ref{x_mu-conv-to-x-lem}, we can apply
 \eqref{I-mu=I+remainder} and we have
\[\begin{array}{l}
\ds{I^\mu_{-\infty}(z_\mu) \leq \frac{c\,\mu^2}{2} \int_{-\infty}^0 |\frac{\partial^2{\varphi}_\mu}{\partial t^2}(t)|_{H^{2\beta}}^2 dt + I_{-\infty}(\varphi_\mu)}\\
\vs
\ds{
+ \mu \int_{-\infty}^0 \left<Q^{-1} \frac{\partial^2{\varphi}_\mu}{\partial t^2}(t),Q^{-1} \left( \frac{\partial{\varphi}_\mu}{\partial t}(t) - A \varphi_\mu(t) - B(\varphi_\mu(t)) \right)\right>_H dt}\\
\vs
\ds{
\leq \frac{\mu^2}{2} \int_{-\infty}^0 |\frac{\partial^2{\varphi}_\mu}{\partial t^2}(t)|_{H^{2\beta}}^2 + I_{-\infty}(\varphi_\mu)+  \mu \left( \int_{-\infty}^0 |\frac{\partial^2{\varphi}_\mu}{\partial t^2}t)|_{H^{2\beta}}^2 dt \right)^{1/2} \left( I_{-\infty}(\varphi_\mu) \right)^{1/2}.}
  \end{array}
\]
  By \eqref{m81}, this implies
  \[
I^\mu_{-\infty}(z_\mu) \leq I_{-\infty}(\varphi_\mu) + c \mu^{2-2\alpha} \left|\frac{\partial{\varphi}}{\partial t} \right|_{L^2((-\infty,0);H^{2\beta})}^2+c \mu^{1 - \alpha} \left| \frac{\partial{\varphi}}{\partial t}  \right|_{L^2((-\infty,0);H^{2\beta})} \left(I_{-\infty}(\varphi_\mu) \right)^{1/2},\]
and by  \eqref{convolution-L2-convergence-eq} and \eqref{convolution-deriv-L2-convergence-eq}
\[    \lim_{\mu \downarrow 0} I_{-\infty}(\varphi_\mu) = I_{-\infty}(\varphi) = V(u).
  \]
  Therefore, if we pick $\alpha<1$ in  \eqref{m83}, we get
  \[
    \limsup_{\mu \downarrow 0} \bar{V}_\mu(u_\mu) \leq \limsup_{\mu \downarrow 0} I^\mu_{-\infty}(z_\mu) \leq V(u).
  \]
Since, in view of \eqref{m82}  and Theorem \ref{continuity-of-V-mu-tilde-thm},
  \begin{equation*}
    \limsup_{\mu \downarrow 0} \bar{V}_\mu(u_\mu) = \limsup_{\mu \downarrow 0} \bar{V}_\mu(u)
  \end{equation*}
we can conclude that \eqref{m84} holds.

\subsection{Lower bound}
\label{sec:8-lower}


Let $N \subset H$ be a closed set with $N \cap H^{1+ 2\beta} \not = \emptyset$.  In particular, by Theorem \ref{continuity-of-V-mu-tilde-thm} we have $\inf_{u \in N} \bar{V}_\mu(u) < +\infty$.
Due to \eqref{m68} and Theorem \ref{compact-level-sets-half-line-thm}, there exists   $z^\mu \in\,C((-\infty,0];\H)$ such that
\[u^\mu:=\Pi_1 z^\mu(0) \in N,\ \ \ \ I^\mu_{-\infty}(z^\mu) = \bar{V}_\mu(u^\mu) = \inf_{u \in N} \bar{V}_\mu(u).\] Now, let $\psi^\mu \in L^2((-\infty,0);H)$ be the minimal control such that
\[z^\mu(t) = \int_{-\infty}^t \Smu(t-s)  B_\mu(z^\mu(s))  ds+ \int_{-\infty}^t \Smu(t-s) \Qmu \psi^\mu(s)  ds,\]
and
\begin{equation}
\label{l2}
  \inf_{u \in N} \bar{V}_\mu(u) = \bar{V}_\mu(u^\mu) = \frac{1}{2} \left| \psi^\mu \right|_{L^2((-\infty,0);H)}^2.
\end{equation}

In what follows, we shall denote $v^\mu=\Pi_2 z^\mu(0).$
For any $\delta>0$, we  define the approximate control
\[
  \psi^{\mu,\delta}(t) = (I - \delta A)^{-\frac{1}{2}} \psi^{\mu}(t),\ \ \ t\leq 0,
\]
and in view of \eqref{m56} and \eqref{rev2} we can define $z^{\mu,\delta}$ to be the solution to the corresponding control problem
\[
  z^{\mu,\delta} (t) = \int_{-\infty}^t  \Smu(t-s)  B_\mu(z^{\mu,\delta}(s))  ds+ \int_{-\infty}^t  \Smu(t-s) \Qmu \psi^{\mu,\delta}(s)  ds.\]
Notice that, according to \eqref{m30},
\[\lim_{t\to-\infty}|z^{\mu,\delta}|_{\H_{1+2\beta}}=0.\]
Moreover, as $  \psi^{\mu,\delta} \in\,L^2((-\infty,0);H^1)$, thanks to \eqref{m30} we have
\[\lim_{t\to-\infty}|z^{\mu,\delta}|_{\H_{2(1+\beta)}}=0.\]
In what follows, we shall denote $(u^{\mu,\d},v^{\mu,\d})=z^{\mu,\d}(0).$

\begin{Lemma} \cite[Lemma 8.1]{cs-annals}
\label{endpoint-difference-theorem}
There exists $\mu_0>0$ such that,
 \begin{equation}
 \label{m89}
    \lim_{\d\to 0}\,\sup_{\mu\leq \mu_0}\,\left| u^\mu - u^{\mu,\delta}  \right|_{H^{2\beta}}^2  =0.
  \end{equation}
\end{Lemma}

From the previous lemma, we get  the main result of this section.
\begin{Theorem}\cite[Theorem 8.2]{cs-annals}
\label{t.82}
  For any closed $N \subset H$, we have
  \begin{equation}
  \label{m-fine101}
    \inf_{u \in N} V(u) \leq \liminf_{\mu \downarrow 0} \inf_{u \in N} \bar{V}_\mu(u).
  \end{equation}

\end{Theorem}

 If the right hand side in \eqref{m-fine101} is infinite,  the theorem is trivially true.  Therefore, in what follows we can assume that
 \begin{equation}
 \label{m-fine100}
\liminf_{\mu \to 0} \inf_{u \in N} \bar{V}_\mu(u) < +\infty.
 \end{equation}

  We first observe that, if we define
  \[\varphi^{\mu,\d}(t)=\Pi_1 z^{\mu,\d}(t),\ \ \ t\leq 0,\]
  in view of \eqref{I-mu=I+remainder}
\[
 \begin{array}{l}
 \ds{V(x^{\mu,\delta}) \leq I_{-\infty}(\varphi^{\mu,\delta}) =I^\mu_{-\infty}(z^{\mu,\delta}) - \frac{\mu^2}{2} \int_{-\infty}^0 \left|Q^{-1}\frac{\partial^2{\varphi}^{\mu,\delta}}{\partial t^2}(t) \right|_H^2 dt }\\
 \vs
 \ds{- \mu \int_{-\infty}^0 \left<Q^{-1} \frac{\partial^2{\varphi}^{\mu,\delta}}{\partial t^2}(t),  Q^{-1} \frac{\partial{\varphi}^{\mu,\delta}}{\partial t}(t) - Q^{-1} A \varphi^{\mu,\delta}(t) - Q^{-1}B(\varphi^{\mu,\delta}(t)) \right>_H dt.}
    \end{array}
    \]
Since
\begin{equation}
\label{m91}
|\psi^{\mu,\d}(t)|_H=|(I-\d A)^{-1/2}\psi^\mu(t)|_H\leq |\psi^\mu(t)|_H,\ \ \ \ t\leq 0,
\end{equation}
we have
\[I^\mu_{-\infty}(z^{\mu,\delta})\leq I^\mu_{-\infty}(z^{\mu})=\inf_{x \in N}\bar{V}_\mu(u),\]
so that
\[\begin{array}{l}
\ds{ V(u^{\mu,\delta}) \leq  \inf_{u \in N} \bar{V}_\mu(u)} \\
\vs
\ds{- \mu \int_{-\infty}^0 \left<Q^{-1} \frac{\partial^2{\varphi}^{\mu,\delta}}{\partial t^2}(t),  Q^{-1} \frac{\partial{\varphi}^{\mu,\delta}}{\partial t}(t) - Q^{-1} A \varphi^{\mu,\delta}(t) - Q^{-1}B(\varphi^{\mu,\delta}(t)) \right>_H dt.}
\end{array}\]
Thanks to \eqref{m30} and Hypothesis \ref{H9}, by integrating by parts
\begin{equation}
\label{remainder}
\begin{array}{l}
\ds{V(u^{\mu,\d})\leq  \inf_{u \in N} \bar{V}_\mu(u) }\\
\vs
\ds{-\frac{\mu}{2} |Q^{-1} v^{\mu,\delta} |_H^2
    - \mu \left<(-A) Q^{-1} u^{\mu,\delta}, Q^{-1} v^{\mu,\delta} \right>_H
    + \mu \left< Q^{-1}B(u^{\mu,\delta}), Q^{-1} v^{\mu,\delta} \right>_H}\\
    \vs
    \ds{+ c\,\mu \int_{-\infty}^0 \left| \frac{\partial{\varphi}^{\mu,\delta}}{\partial t}(t) \right|_{H^{1+2\beta}}^2 dt
    + c\,\gamma_{2\beta} \mu \int_{-\infty}^0 \left|\frac{\partial{\varphi}^{\mu,\delta}}{\partial t}(t) \right|_{H^{2\beta}}^2 dt= \inf_{u \in N}\bar{V}_\mu(u)+\sum_{i=1}^5I^{\mu,\d}_i.}
    \end{array}
    \end{equation}

First, we note that
\begin{equation}
  \label{I_1-est}
  I^{\mu,\delta}_1 \leq 0.
\end{equation}
Next, by \eqref{m43} see that
\begin{equation*}
\begin{array}{l}
  \ds{I^{\mu,\delta}_2 + I^{\mu,\delta}_4 \leq c\sqrt{\mu} \left( |u^{\mu,\delta}|_{H^{2\beta +2}}^2 + \mu |v^{\mu,\delta}|_{H^{2\beta +1 }}^2  \right) + c\mu \int_{-\infty}^0 |z^{\mu,\delta}(t)|_{\H^{2 + 2\beta}}^2 dt} \\
  \ds{\leq c(\mu + \sqrt{\mu}) \int_{-\infty}^0 |\psi^{\mu,\delta}(t)|_{H^1}^2 dt}.
\end{array}
\end{equation*}
Since for any $h \in\,H$
we have $(I - \delta A)^{-\frac{1}{2}} h \in \textnormal{Dom}(-A)^{\frac{1}{2}}$ and
\[\left| (-A)^{\frac{1}{2}} (I - \delta A)^{-\frac{1}{2}} h \right|_H \leq \delta^{-1/2} \left|h \right|_H,\]
we have
\[\left|\psi^{\mu,\delta}(t) \right|_{H^1}\leq \d^{-1/2}\,\left|\psi^{\mu}(t) \right|_{H},\ \ \ t\leq 0.\]
Therefore, by \eqref{l2},
\begin{equation}
  \label{I_2-I_4-est}
  I^{\mu,\delta}_2 + I^{\mu,\delta}_4 \leq c\, \d^{-1/2}(\mu + \sqrt{\mu}) \int_0^t |\psi^\mu(t)|_H^2 =2\, c \delta^{-1/2}(\mu + \sqrt{\mu}) \inf_{u \in N} \bar{V}_\mu(u).
\end{equation}
By the same arguments, \eqref{m43}, and \eqref{m91} give
\begin{equation}
  \label{I_3-I_5-est}
  I^{\mu,\delta}_3 + I^{\mu,\delta}_5 \leq c( \mu + \sqrt{\mu}) \inf_{u \in N} \bar{V}_\mu(u).
\end{equation}
Combining together \eqref{I_1-est}, \eqref{I_2-I_4-est}, and \eqref{I_3-I_5-est} with \eqref{remainder},  we obtain,
\[
  V(u^{\mu,\delta}) \leq \inf_{u \in N} \bar{V}_\mu(u)+ c(\mu + \sqrt{\mu})(1 + \delta^{-1/2})\, \inf_{u \in N} \bar{V}_\mu(u).
\]
From this, choosing $\delta=\sqrt{\mu}$, and due to \eqref{m-fine100}, we see that
\[\liminf_{\mu\to0}V(u^{\mu,\sqrt{\mu}})\leq \liminf_{\mu\to0} \inf_{u \in N}\bar{V}_\mu(u).\]
Since we are assuming \eqref{m-fine100}, and, by \cite[Proposition 5.1]{cerrok}, the level sets of $V$ are compact, there is a sequence $\mu_n \to 0$ and $u^0 \in\,H$ such that
\[\lim_{n\to \infty}|u^{\mu_n,\sqrt{\mu_n}}-u^0|_H=0,\ \ \ V(u^0) \leq \liminf_{\mu \to 0} V(u^{\mu,\sqrt{\mu}}).\]
By \eqref{m89}, we have that $u^{\mu_n}$ converges to $u^0$ in $H$, so that $u^0 \in N$.  This means that  we can conclude, as
\[\inf_{u \in N} V(u) \leq V(u^0) \leq \liminf_{\mu \to 0} V(u^{\mu,\sqrt{\mu}}) \leq \liminf_{\mu \to 0} \inf_{u \in N} \bar{V}_\mu(u).   \]

\subsection{Application to the exit problem} \label{sec:9-exit}
We are here interested in the problem of the exit of the solution $u^\mu_\e$ of equation \eqref{wave-eq} from a domain $D\subset H$, for any $\mu>0$ fixed. Then we  apply the limiting results proved in Theorems \ref{upper-bound-thm} and \ref{t.82} to show that, when $\mu$ is small, the relevant quantities in the exit problem from $D$ for the solution $u^\mu_\e$ of equation \eqref{wave-eq}  can be approximated by the corresponding ones arising  for equation \eqref{heat-eq}.

First, let us give some assumptions on the set $D$.

\begin{Hypothesis}
\label{H10}
The domain $D \subset H$ is an open, bounded, connected set, such that  $0 \in D$.
Moreover, for any $h \in \partial D \cap H^{1 + 2\beta}$ there exists a sequence $\{h_n\}_{n \in \nat} \subset \bar{D}^c \cap H^{1+2\beta}$ such that
\begin{equation} \label{boundary-reg-assum}
  \lim_{n \to +\infty} \left| h_n - h \right|_{\H_{1 + 2\beta}} =0.
\end{equation}
\end{Hypothesis}
Assume now that $D$ is an open, bounded and connected set such that, for any $h \in \partial D \cap H^{1 + 2\beta}$, there exists a $k \in \bar{D}^c \cap H^{1 + 2\beta}$ such that
\begin{equation}
\label{r.1}
\{tk+ (1-t)h: 0<t\leq1\} \subset \bar{D}^c.\end{equation}
Then it is immediate to check that \eqref{boundary-reg-assum} is satisfied.
Condition \eqref{r.1}  is true, for example, if $D$ is convex, because of the Hahn-Banach separation theorem and the density of $H^{1 + 2\beta}$ in $H$.

\begin{Lemma}\cite[Lemma 9.1]{cs-annals}
Under Hypothesis \ref{H3}
\[\bar{V}_\mu(\partial D):=  \inf_{u \in \partial D} \bar{V}_\mu(u)=\bar{V}_\mu(u_{D,\mu})<\infty,
\]
for some $u_{D,\mu} \in\,\partial D\cap \H_{1+2\beta}$.
\end{Lemma}

Now, if we denote by $z^\mu_{\e,z_0}=(u^\mu_{\e,z_0},v^\mu_{\e,z_0})$ the mild solution of \eqref{wave-eq}, with initial position and velocity $z_0=(u_0,v_0) \in\,\H$, we define the exit time
\[
  \tau^{\mu,\epsilon}_{z_0} = \inf \left\{t>0: u^\mu_{\epsilon,z_0} (t) \not \in D \right\}.
\]
Here is the main result of this section
\begin{Theorem}\cite[Theorem 9.2]{cs-annals}
\label{m-t81}
  There exists $\mu_0>0$ such that for $\mu<\mu_0$ the following conditions are verified.
  For any $z_0=(u_0,v_0) \in\,\H$ such that $u_0 \in\,D$ and $u^\mu_{0,z_0}(t) \in\,D$, for $t\geq 0$,
  \begin{enumerate}
    \item The exit time has the following asymptotic growth
    \[
      \lim_{\epsilon \to 0} \epsilon \log  \E \left( \tau^{\mu,\epsilon}_{z_0} \right)= \inf_{u \in \partial D} \bar{V}_\mu(u),
    \]
  and for any $\eta>0$,
    \[
      \lim_{\epsilon \to 0} \Pro \left( \exp\le(\frac{1}{\epsilon}(\bar{V}_\mu(\partial D) -\eta)\r) \leq \tau^{\mu,\epsilon}_{z_0} \leq \exp\le(\frac{1}{\epsilon}(\bar{V}_\mu(\partial D) + \eta)\r) \right) =1.
    \]

    \item  For any closed $N \subset \partial D$ such that $\ds{\inf_{u\in N}\bar{V}_\mu(u) > \inf_{u \in \partial D} \bar{V}_\mu (u)}$, it holds
    \[
      \lim_{\epsilon \to 0} \Pro\left( u^\mu_{\epsilon,z_0}(\tau^{\mu,\epsilon}_{z_0}) \in N \right) = 0.
    \]
  \end{enumerate}
\end{Theorem}

\begin{Remark}
{\em   The requirement that $u^\mu_{0,z_0}(t) \in D$ for all $t \geq 0$ is necessary because in Lemma \ref{large-init-veloc-causes-exit-lem} we showed that there exist $z_0\in D \times H^{-1}$ such that $u^\mu_{0,z_0}$ leaves $D$ in finite time.  Of course, for these initial conditions, the stochastic processes $u^\mu_{\epsilon,z_0}$ will also exit in finite time for small $\epsilon$.}
\end{Remark}

In \cite{tran} it has been proven that an analogous result to Theorem \ref{m-t81} holds for equation \eqref{heat-eq}.
If we denote by $u_{\epsilon,u_0}$  the mild solutions of equation \eqref{heat-eq}, with initial condition $u_0 \in\,H$, we define  the exit time
\[
  \tau^\epsilon_{u_0} = \inf \left\{t>0: u_{\epsilon,u_0}(t) \not \in D \right\}.
\]
In \cite{tran} it has been proven that for any $u_0 \in\,D$ such that $u_{0,u_0}(t) \in\,D$, for any $t\geq 0$, it holds
\[
\lim_{\epsilon \to 0} \epsilon \log  \E \left( \tau^{\epsilon}_{u_0} \right)= \inf_{u \in \partial D} V(u).
\]
Similarly, as we would expect, it also holds that
\[\lim_{\e \to 0} \e \log \tau^\e_{u_0} = \inf_{u \in \partial D} V(u),\ \ \ \text{in probability},\]
and if $N \subset \partial D$ is closed and $\inf_{u \in N} V(u)> \inf_{u \in \partial D} V(u)$,
\[\lim_{\e \to 0} \Pro \left(u^\e_{u_0}(\tau^\e_{u_0}) \in N \right) =0.\]
The proof of these facts is analogous to the proof of Theorem \ref{m-t81}.

In view of what we have proven in Sections \ref{sec:7-upper} and \ref{sec:8-lower} and of Theorem \ref{m-t81}, this implies   that
the following Smoluchowski-Kramers approximations holds for the exit time.

\begin{Theorem}\cite[Theorem 9.4]{cs-annals}
\begin{enumerate}
  \item
  For any initial conditions $z_0 = (u_0,v_0)$,
  \[
    \lim_{\mu \to 0} \lim_{\epsilon \to 0} \e\log\,\E \left(\tau^{\mu,\epsilon}_{z_0} \right) = \lim_{\epsilon \to 0} \e\log\,\E \left( \tau^{\epsilon}_{u_0} \right)=\inf_{x \in\,\partial D}V(x).
  \]

  \item  For any $\eta>0$, there exists $\mu_0>0$ such that for $\mu< \mu_0$
  \begin{equation}
  \label{m-fine200}
    \lim_{\epsilon \to 0} \Pro \left( e^{\frac{1}{\epsilon}(\bar{V} -\eta)} \leq \tau^{\mu,\epsilon}_{z_0} \leq e^{\frac{1}{\epsilon}(\bar{V} + \eta)}  \right) =1
  \end{equation}

  \item  For any $N \subset \partial D$ such that $\inf_{x \in N} V(x) < \inf_{x \in \partial D} V(x)$, there exits $\mu_0>0$ such that for all $\mu < \mu_0$,
  \begin{equation*}
    \lim_{\epsilon \to 0} \Pro_{z_0} \left( u^\mu_\epsilon(\tau^{\mu,\epsilon}) \in N \right) = 0.
  \end{equation*}
\end{enumerate}
\end{Theorem}

We recall that in Theorem \ref{m2-10} we have proved that, in the case of gradient systems, for any $\mu>0$
\[\bar{V}_\mu(x)=V(x),\ \ \ \ x \in\,H.\]
This means that in this case for any $z_0=(u_0,v_0) \in\,\H$ and $\mu>0$  such that the unperturbed system $u^\mu_{0,z_0}(t) \in D$ for all $t>0$
\[\lim_{\epsilon \to 0} \e\,\log \E \left(\tau^{\mu,\epsilon}_{z_0} \right) = \lim_{\epsilon \to 0} \e\,\log \E \left( \tau^{\epsilon}_{u_0} \right)=\inf_{x \in\,\partial D}V(x).\]
and \eqref{m-fine200} holds for any $\mu>0$.

\end{document}